\documentclass[a4paper,twoside,10pt]{article}

\usepackage[a4paper,left=3cm,right=3cm, top=3cm, bottom=3cm]{geometry}
\usepackage[latin1]{inputenc}
\usepackage{graphicx}
\usepackage{epstopdf}
\usepackage{amsmath}
\usepackage{amsthm}
\usepackage{amssymb}
\usepackage{hyperref}
\usepackage{subfigure}
\usepackage{float}
\usepackage{cite}
\usepackage{color}
\usepackage[font=footnotesize,labelfont=bf]{caption}
\usepackage{tikz}
\usepackage{algorithmic}
\usepackage{rotating}
\usetikzlibrary{calc,arrows,quotes,angles,shapes, plotmarks,fit}
\tikzstyle{rect} = [draw, rectangle, minimum height = 4em, text width = 6em, text centered]\tikzstyle{sum} = [draw, fill=blue!20, circle, node distance=1cm]
\tikzstyle{line} = [draw, -latex']
\usepackage[author={Lorenzo}]{pdfcomment}
\usepackage{soul,xcolor}

\restylefloat{table}
\theoremstyle{plain}

\theoremstyle{definition}

\theoremstyle{remark}
\newtheorem{remark}{Remark}

\newcommand{\p}{p} % polynomial degree
\newcommand{\F}{F} % generic face
\newcommand{\E}{K} % generic polyhedron
\newcommand{\e}{s} % generic edge
\newcommand{\np}     {n_\p} 	 % dimension polynomial of degree \p
\newcommand{\npmu}{n_{\p-1}} % dimension polynomial of degree \p-1
\newcommand{\npmd}{n_{\p-2}} % dimension polynomial of degree \p-2
\newcommand{\m}{m} % symbol denoting polynomial element WITHOUT index
\newcommand{\malpha}{\m_\alpha} % polynomial subscript \alpha
\newcommand{\mbeta}{\m_\beta} % polynomial subscript \beta
\newcommand{\mbetaF}{\m_\beta^\F} % polynomial subscript \beta^\F
\newcommand{\mgamma}{\m_\gamma} % polynomial subscript \gamma
\newcommand{\mgammaF}{\m_\gamma^\F} % polynomial subscript \gamma^\F
\newcommand{\mdelta}{\m_\delta} % polynomial subscript \delta
\newcommand{\mdeltaF}{\m_\delta^\F} % polynomial subscript \delta^|F
\newcommand{\mbar}{\overline \m} % symbol denoting polynomial element WITHOUT index bar
\newcommand{\mbaralpha}{\overline \m_\alpha} % polynomial subscript \alpha bar
\newcommand{\mbarbeta}{\overline \m_\beta} % polynomial subscript \beta bar
\newcommand{\mbargamma}{\overline \m_\gamma} % polynomial subscript \gamma bar
\newcommand{\mbargammaF}{\overline \m_\gamma^\F} % polynomial subscript \gamma bar^\F
\renewcommand{\span}{\text{span}} % span
\newcommand{\f}{f} % continuous rhs
\renewcommand{\u}{u} % continuous solution
\newcommand{\up}{\u_\p} % discrete solution
\newcommand{\vv}{v} % continuous test function
\newcommand{\vp}{\vv_\p} % discrete test function
\newcommand{\V}{V} % continuous space
\newcommand{\Vp}{\V_\p} % discrete space
\newcommand{\Vtildep}{\widetilde\V_\p} % discrete space tilde
\renewcommand{\a}{a} % global continuous bilinear form
\newcommand{\ap}{\a_\p} % global discrete bilinear form
\newcommand{\aE}{\a^\E} % local continuous bilinear form
\newcommand{\SE}{S^\E} % local stabilization
\newcommand{\aF}{\a^\F} % local continuous bilinear form (face)
\newcommand{\apE}{\ap^\E} % local discrete bilinear form
\newcommand{\taun}{\mathcal T_n} % polyhedral mesh
\newcommand{\Pinablap}{\Pi^{\nabla}_\p} % operator Pinabla in the bulk
\newcommand{\PinablapF}{\Pi^{\nabla,\F}_\p} % operator Pinabla on the face \F
\newcommand{\PinablaF}{\PinablapF} % operator Pinabla on the face \F
\newcommand{\PinablabarF}{\overline\Pi^{\nabla,\F}_\p} % operator Pinabla bart on the face \F
\newcommand{\Pizp}{\Pi^0_\p} % L2 projection
\newcommand{\PizpF}{\Pi^{0,\F}_\p} % L^2 projection on face \F
\newcommand{\q}{q} % generic polynomial
\newcommand{\qp}{\q_\p} % generic polynomial of degree \p
\newcommand{\NV}{N_V} % number of vertices
\newcommand{\NS}{N_S} % number of (local) skeletal dofs
\newcommand{\NF}{N_F} % number of (local) face dofs
\newcommand{\NB}{N_B} % number of (local) bulk dofs
\newcommand{\h}{h} % mesh size function
\newcommand{\hF}{\h_\F} % diameter of face \F
\newcommand{\hE}{\h_\E} % diameter of polyhedron \E
\newcommand{\dof}{\text{dof}} % dof
\newcommand{\dofbar}{\overline{\text{dof}}} % dof bar
\newcommand{\xbold}{\mathbf x} % x bold
\newcommand{\xFbold}{\mathbf{x_\F}} % barycenter face \F
\newcommand{\xF}{x_\F} % barycenter face \F, coordinate x
\newcommand{\yF}{y_\F} % barycenter face \F, coordinate y
\newcommand{\xEbold}{\mathbf{x_\E}} % barycenter polyhedron \E
\newcommand{\xE}{x_\E} % barycenter polyhedron \E, coordinate x
\newcommand{\yE}{y_\E} % barycenter polyhedron \E, coordinate y
\newcommand{\zE}{z_\E} % barycenter polyhedron \E, coordinate z
\newcommand{\boldalpha}{\boldsymbol{\alpha}} % \alpha bold
\newcommand{\boldbeta}{\boldsymbol{\beta}} % \alpha bold
\newcommand{\boldgamma}{\boldsymbol{\gamma}} % \alpha bold
\newcommand{\standard}{standard } % standard
\newcommand{\orthogonal}{orthogonal } % orthogonal
\newcommand{\hybrid}{hybrid } % hybrid
\newcommand{\G}{\mathbf G} % VEM matrix G
\newcommand{\Gtilde}{\widetilde{\mathbf G}} % VEM matrix Gtilde
\newcommand{\B}{\mathbf B} % VEM matrix B
\newcommand{\D}{\mathbf D} % VEM matrix D
\renewcommand{\H}{\mathbf H} % VEM matrix H
\newcommand{\C}{\mathbf C} % VEM matrix C
\renewcommand{\L}{\mathbf L} % VEM matrix C
\newcommand{\Gbar}{\overline{\mathbf G}} % VEM matrix G bar
\newcommand{\Gtildebar}{\overline{\widetilde{\mathbf G}}} % VEM matrix Gtilde bar
\newcommand{\Bbar}{\overline{\mathbf B}} % VEM matrix B bar
\newcommand{\Dbar}{\overline{\mathbf D}} % VEM matrix D bar
\newcommand{\Hbar}{\overline{\mathbf H}} % VEM matrix H bar
\newcommand{\Cbar}{\overline{\mathbf C}} % VEM matrix C bar
\newcommand{\Lbar}{\overline{\mathbf L}} % VEM matrix L bar
\newcommand{\Mbar}{\overline{\mathbf M}} % VEM matrix M bar
\newcommand{\Lambdabar}{\overline{\mathbf \Lambda}} % VEM matrix \Lambda bar
\newcommand{\Ndof}{N_{dof}} % number of dofs
\newcommand{\GS}{\mathbf {GS}} % Gram Schmidt matrix
\newcommand{\GSF}{\mathbf {GS}^\F} % Gram Schmidt matrix ^\F
\newcommand{\NVE}{\NV} % number of vertices of \E
\newcommand{\Id}{\textbf{Id}} % identity
\newcommand{\phibar}{\overline \varphi} % phi bar
\newcommand{\n}{\mathbf n} % normal derivative
\newcommand{\lambdabar}{\overline \lambda} % lambdabar
\newcommand{\chibar}{\overline \chi} % \chi bar
\newcommand{\mubar}{\overline \mu} % \mu bar
\newcommand{\hyb}{{hyb}} % hyb
\newcommand{\upi}{\u_\pi} % piecewise discontinuous polynomial approximating \u
\newcommand{\uI}{\u_I} % VE function approximating \u
\providecommand{\keywords}[1]{\textbf{\textit{Keywords: }} #1}

\begin{tiny}
\author{
\normalsize{
F. Dassi \thanks{Dip. di Matematica e Applicazioni,  Universit\`a degli Studi di Milano-Bicocca, E-mail: {\tt franco.dassi@unimib.it}} \quad 
L. Mascotto
\thanks{Dip. di Matematica,  Universit\`a degli Studi di Milano, E-mail: {\tt lorenzo.mascotto@unimi.it}}
\thanks{Inst. f\"ur Mathematik, C. von Ossietzky Universit\"at Oldenburg, E-mail: {\tt lorenzo.mascotto@uni-oldenburg.de}}
}}
\end{tiny}

\date{}

\title{\textbf{\normalsize{Exploring High-order three dimensional Virtual Elements: bases and stabilizations}}}

\begin{document}
%%%%%%%%%%%%%%%%%%%%%%%%%%%%%%%%%%%%%
\maketitle

\begin{abstract}
We present numerical tests of the Virtual Element Method (VEM) tailored for the discretization of a three dimensional Poisson problem with high-order ``polynomial'' degree (up to $\p=10$).
Besides, we discuss possible reasons for which the method could return suboptimal-wrong error convergence curves.
Among these motivations, we highlight ill-conditioning of the stiffness matrix and not particularly ``clever'' choices of the stabilizations.
We propose variants of the definition of face/bulk degrees of freedom, as well as of stabilizations, which lead to methods that are much more robust in terms of numerical performances.
\end{abstract}

\keywords{Virtual Element Method, Polyhedral meshes, high-order methods, ill-conditioning}

\medskip

\begin{flushright}
{\footnotesize\textsl{\indent Darest thou now, O Soul,\\
Walk out with me toward the Unknown Region,\\
Where neither ground is for the feet, nor any path to follow?\\}}

%\medskip
%{\footnotesize \textsl{No map, there, nor guide,\\
%Nor voice sounding, nor touch of human hand,\\
%Nor face with blooming flesh, nor lips, nor eyes, are in that land\dots}\\}
\medskip

{\footnotesize{\textrm{Walt Whitman, Leaves of Grass, 1855.}}}
\end{flushright}

% --------------------------------------------------------------------------------------------------------------------------------------------------------------------------------------------------------------------------------------------------------------------
\section {Introduction} \label{section introduction}
% --------------------------------------------------------------------------------------------------------------------------------------------------------------------------------------------------------------------------------------------------------------------
The Virtual Element Method (VEM) is a generalization of the Finite Element Method (FEM) that allows for general polytopal meshes,
thus including non-convex elements and hanging nodes.

Approximation spaces in VEM contain locally polynomials and, more in general, consist of functions which solve local problems mimicking the original ones and, consequently, are not known in a closed form (hence the name virtual).
For this reason, the operators involved in the discretization of the problem are not computed exactly;
rather, the construction of the method is based on two ingredients: proper projectors onto piecewise discontinuous polynomial spaces and stabilizing bilinear forms mimicking their continuous counterparts.
Both ingredients can be computed exactly only with the aid of the degrees of freedom.

Although the VEM technology is very recent, it has been applied to a large number of two dimensional problems;
a short of list of them is: \cite{equivalentprojectorsforVEM, bbmr_VEM_generalsecondorderelliptic,BLV_StokesVEMdivergencefree,VEMelasticity,cangianigeorgulispryersutton_VEMaposteriori,serendipityVEM,Helmholtz-VEM};
in particular, high-order VEM are investigated in \cite{hpVEMbasic, preprint_hpVEMcorner, pVEMmultigrid, preprint_HarmonicVEM, fetishVEM}.

The literature dealing with three dimensional problems is much less broad, see \cite{gain2014virtual, Topology-VEM, chi2017some}.
The only attempt, at the best of our knowledge, to increase the order of VEM in 3D is \cite{preprint_VEM3Dbasic}, where the highest order achieved in numerical tests is $\p=5$.

\medskip

In the present work, we have a double aim.
Firstly, we present numerical tests for three dimensional VEM of order higher than $5$, thus %enriching what done in \cite{preprint_VEM3Dbasic} up to order $\p=10$.
inviting the reader to raise the anchor from the safe port \cite{preprint_VEM3Dbasic} and to ``\emph{\dots Walk out with us toward the Unknown Region\dots}'', reaching in fact the pinnacle of degree of accuracy $\p=10$.
Secondly, we numerically investigate the reasons of possible suboptimal/wrong behaviour in the error convergence curves,
highlighting two among them: the ill-conditioning of the linear system stemming from the method and the choice of the stabilization.

We tackle the (possible) issue of suboptimality of VEM, when considering its $\h$ and $\p$ versions as well as when it is applied to meshes with elements having collapsing bulk,
by proposing two novel approaches of the definition of the face/bulk degrees of freedom and proposing three different stabilizations.

The outline of the paper follows. In Section \ref{section VEM}, we review the construction of three dimensional VEM,
emphasizing in particular different stabilizations and face/bulk degrees of freedom.
Next, in Section \ref{section numerical results}, we provide a number of numerical results comparing the effects on the method of the above-mentioned stabilizations and degrees of freedom;
more precisely, we study the $\h$ and the $\p$ versions of the method, paying attention also to VEM applied to meshes with degenerate elements.
Concluding remarks are stated in Section \ref{section conclusions}.
Finally, in Appendix \ref{appendix hitchhikers}, we give some hints regarding the implementation of the method with the novel canonical basis functions.

\medskip 

\paragraph*{Notation.} %We fix here some notations.
By $\mathbb P_\p(\F)$ and $\mathbb P_\p(\E)$, $\p \in \mathbb N$, we denote the spaces of two and three dimensional polynomials of degree $\p$ over a polygon $\F$ and a polyhedron $\E$, respectively;
if $\p=-1$, then we set $\mathbb P_{-1}(\F) = \mathbb P_{-1}(\E) = \emptyset$. Moreover, we fix:
\begin{equation} \label{dimension polynomial spaces}
\np^\F = \dim(\mathbb P_\p(\F)),\quad \quad \np = \dim (\mathbb P_\p(\E)) \quad \forall\, \p \in \mathbb N.
\end{equation}
Assume now that we are given $\{\malpha\}_{\alpha=1}^{\np}$, $\p\in \mathbb N$, a basis of $\mathbb P_\p(\E)$ such that:
\begin{equation} \label{assumption splitting polynomials}
\span \left(\{\malpha\}_{\alpha=1}^{\npmd} \right) = \mathbb P_{\p-2}(\E) \quad \text{and } \quad \span \left(\{\malpha\}_{\alpha=1}^{\npmu} \right) = \mathbb P_{\p-1}(\E).
\end{equation}
It will be convenient to split the polynomial basis into:
\begin{equation} \label{polynomial basis splitting}
\{\malpha\}_{\alpha=1}^{\np} = \{\malpha\}_{\alpha=1}^{\npmd} \cup \{\malpha\}_{\alpha=\npmd+1}^{\npmu} \cup \{\malpha\}_{\alpha=\npmu+1}^{\np}.
\end{equation}
%and to fix the following notation:
%\begin{equation} \label{splitting cut polynomials}
%\mathbb P_{\p-1}^\cut (\E) = \span (\{\malpha\}_{\alpha=\npmd+1}^{\npmu}),\quad \mathbb P_\p^\cut(\E) = \span(\{\malpha\}_{\alpha=\npmu+1}^{\np}).
%\end{equation}
We assume that the polygonal counterpart of \eqref{assumption splitting polynomials} holds true; consequently, we can consider a splitting analogous to the one in \eqref{polynomial basis splitting}
on \[ \mathbb P_\p (\F) = \span \left( \{\malpha^\F\}_{\alpha=1}^{\np^\F} \right),\] the space of polynomial of degree $\p$ over polygon $\F$.

% --------------------------------------------------------------------------------------------------------------------------------------------------------------------------------------------------------------------------------------------------------------------
\section {VEM: definition, stabilizations and bases} \label{section VEM}
% --------------------------------------------------------------------------------------------------------------------------------------------------------------------------------------------------------------------------------------------------------------------
In this section, we introduce a family of VEM tailored for the approximation of the following Poisson problem in three dimensions with (for simplicity) homogeneous boundary conditions.
Given $\Omega \subset \mathbb R^3$ a polyhedral domain and $\f \in L^2(\Omega)$:
\begin{equation} \label{weak formulation Poisson problem}
\begin{cases}
\text{find } \u \in \V \text{ s. th.}\\
\a(\u,\vv) = (\f,\vv) \quad \forall \, \vv \in \V
\end{cases},
\end{equation}
where:
\begin{equation} \label{fixing simple notation}
\V = H^1_0(\Omega) , \quad \quad \a(\cdot, \cdot) = (\nabla \cdot, \nabla \cdot)_{0,\Omega}.
\end{equation}
In Section \ref{subsection family VEM}, we briefly recall from \cite{equivalentprojectorsforVEM, preprint_VEM3Dbasic, hitchhikersguideVEM}
the construction of three dimensional VEM for the approximation of the solution of problem \eqref{weak formulation Poisson problem},
keeping yet at a very general level the definition of the degrees of freedom and of the stabilization of the method, typical of the VEM framework.
Various choices of stabilizations as well as of face/bulk degrees of freedom are investigated in Sections \ref{subsection choice stabilizations} and \ref{subsection choice of dofs}, respectively.

%%%%%%%%%%%%%%%%%%%%%%%%%%%%%%%%%%%%%%%%%%%%%%%%%
\subsection{A family of VEM} \label{subsection family VEM}
In this section, we introduce, following \cite{equivalentprojectorsforVEM, preprint_VEM3Dbasic}, a family of VEM in three dimensions for the approximation of problem \eqref{weak formulation Poisson problem}.

The VEM in three dimensions is based on conforming sequences $\taun$ of polyhedra partitioning the physical domain $\Omega$ of the PDE of interest.
By conforming sequence, we mean that, given $\mathcal F_n$, $\mathcal E_n$ and $\mathcal V_n$ the sets of all faces, edges and vertices of the polyhedra in $\taun$, respectively, then %, for all $\E_1$ and $\E_2$ in $\taun$, it can only happen that
%$\overline{\E_1} \cap \overline{\E_2}$ either is empty, a face $\F$ in $\mathcal F_n$, an edge $\e$ in $\mathcal E_n$ or a vertex $\nu$ in $\mathcal V_n$.
all the internal faces $\F \in \mathcal F_n$ must belong to the intersection of two polyhedra.

We observe that, since the aim of the present paper is to test the robustness of the method to mesh-distortion and to increasing ``polynomial degrees'',
no particular geometrical assumptions on the mesh are demanded.

\medskip

Let now $\p\in \mathbb N$; such $\p$ denotes the ``polynomial degree'' of the method.
We begin by defining the local spaces on each face $\F \in \mathcal F_n$:
\begin{equation} \label{local auxiliary space face}
\Vtildep(\F) = \{\vp \in H^1(\F) \mid \Delta \vp \in \mathbb P_{\p}(\F), \, \vp |_{\partial \F} \in \mathbb B_\p(\partial \F)\},
\end{equation}
where:
\begin{equation} \label{polynomials on the skeleton}
\mathbb B_\p (\partial \F) = \{\vp \in \mathcal C^0(\partial \F) \mid \vp|_\e \in \mathbb P_\p(\e) \text{ for all edges } \e \text{ of } \F\}.
\end{equation}

Given \emph{any} polynomial basis $\{\malpha^\F\}_{\alpha=1}^{\npmd^\F}$ of $\mathbb P_{\p-2}(\F)$ satisfying the face counterpart of \eqref{assumption splitting polynomials},
we can endow space \eqref{local auxiliary space face} with the following set of linear functionals. For every $\vp \in \Vtildep(\F)$:
\begin{itemize}
\item the values of $\vp$ at the vertices of $\F$;
\item the values of the $\p-1$ internal Gau\ss-Lobatto nodes on each edge $\e$ of face $\F$;
\item the (scaled) face moments:
\begin{equation} \label{internal moments faces}
\frac{1}{\vert \F \vert} \int_\F \malpha^\F \, \vp \quad \quad \forall \, \alpha =1, \dots, \npmd^\F.
\end{equation}
\end{itemize}
It can be proven, see \cite{VEMvolley}, that it is possible to compute via such linear functionals the energy projector $\PinablapF : \Vtildep(\F) \rightarrow \mathbb P_\p(\F)$ defined as:
\begin{equation} \label{H1 projection face}
\begin{cases}
\medskip
(\nabla \qp^\F, \nabla (\vp - \PinablapF \vp)) _{0,\F}=:\aF(\qp^\F, \vp - \PinablapF \vp)  = 0\\
\medskip
\begin{cases}
\sum_{i=1}^{\NV^\F} (\vp - \PinablapF \vp) (\nu_i^\F) = 0 & \text{if } \p=1\\
\int_\F \vp - \PinablapF \vp = 0 & \text{if } \p \ge 2\\
\end{cases}
\end{cases}\quad \forall \, \qp^\F\in \mathbb P_\p(\F),\, \forall \, \vp \in \Vtildep(\F),\\
\end{equation}
where $\NV^\F$ and $\{\nu_i^\F\}_{i=1}^{\NV^\F}$ denotes the number and the set of vertices of face $\F$, respectively.

Following now \cite{equivalentprojectorsforVEM}, we restrict space $\Vtildep(\F)$ defined in \eqref{local auxiliary space face} so that one is able to compute an $L^2$ projector onto $\mathbb P_\p(\F)$
via the set of linear functionals introduced above.

Such a space, which goes under the name of ``enhanced VE planar space'', is defined as:
\begin{equation} \label{local space face}
\Vp(\F) = \left\{\vp \in \Vtildep(\F) \,\left | \, \int_\F (\vp - \PinablapF \vp) \malpha^\F,\, \text{for all } \alpha=\npmd^\F+1,\dots, \np^\F \right. \right\},
\end{equation}
see \eqref{polynomial basis splitting} for the splitting of the polynomial basis.

As already stressed, it is possible to compute on such space, in addition the the energy projector defined in \eqref{H1 projection face}, the $L^2$ projector $\PizpF:\Vp(\F) \rightarrow \mathbb P_\p(\F)$ defined as:
\begin{equation} \label{L2 projector face rich}
(\qp^\F, \vp - \PizpF \vp)_{0,\F} = 0 \quad \forall \, \qp^\F \in \mathbb P_\p(\F),\, \forall \, \vp \in \Vp(\F).
\end{equation}
\begin{remark}
We observe that the definition of space \eqref{local space face} and of the $H^1$ and $L^2$ orthogonal projectors
are independent of the choice of the polynomial basis employed in the definition of face moments~\eqref{internal moments faces}.
\end{remark}

\medskip

At this point, we are in business for defining local VE spaces on polyhedra.
We begin also in this case by introducing an auxiliary space:
\begin{equation} \label{local auxiliary space bulk}
\Vtildep(\E) = \left\{ \vp \in H^1(\E) \mid \Delta \vp \in \mathbb P_\p(\E),\, \vp|_\F \in \Vp(\F) \text{ for all faces $\F$ of } \E   \right\}.
\end{equation}
Given any polynomial bases $\{ \malpha^\F \}_{\alpha=1}^{\np^\F}$ satisfying the face counterpart of \eqref{assumption splitting polynomials} on all faces $\F$ of $\E$
and any polynomial basis $\{\malpha \}_{\alpha=1}^{\np}$ satisfying \eqref{assumption splitting polynomials},
space \eqref{local auxiliary space bulk} can be endowed with the following set of linear functionals:
\begin{itemize}
\item the values of $\vp$ at the vertices of $\E$;
\item the values of the $\p-1$ internal Gau\ss-Lobatto nodes on each edge $\e$ of polyhedron $\E$;
\item for all faces $\F$ of polyhedron $\E$ the (scaled) face moments:
\begin{equation} \label{internal moments faces polyhedron}
\frac{1}{\vert \F \vert} \int_\F \malpha^\F \, \vp \quad \quad \forall \, \alpha =1, \dots, \npmd^\F;
\end{equation}
\item the (scaled) bulk moments:
\begin{equation} \label{internal moments bulk}
\frac{1}{\vert \E \vert} \int_\E \malpha \, \vp \quad \forall \, \alpha =1, \dots, \npmd.
\end{equation}
\end{itemize}
Such functionals allow to construct the energy projector $\Pinablap : \Vtildep(\E) \rightarrow \mathbb P_\p (\E)$ defined as:
\begin{equation} \label{H1 projector bulk}
\begin{cases}
\medskip
(\nabla \qp, \nabla(\vp - \Pinablap \vp)) =:\aE(\qp, \vp - \Pinablap \vp) = 0\\
\medskip
\begin{cases}
\sum_{i=1}^{\NV} (\vp - \Pinablap \vp)(\nu_i) = 0 & \text{if } \p=1 \\
\int_\E(\vp - \Pinablap \vp) = 0 & \text{if } \p \ge 2\\
\end{cases}\\
\end{cases}\quad \forall \, \qp \in \mathbb P_\p(\E),\, \forall\, \vp \in \Vtildep(\E),
\end{equation}
where $\NV$ and $\{\nu_i\}_{i=1}^{\NV}$ denote the number and the set of vertices of polyhedron $\E$, respectively.

Similarly to the two dimensional case, one can restrict space $\Vtildep(\E)$ defined in \eqref{local auxiliary space bulk} so that it is possible to compute the $L^2$ projection on the space $\mathbb P_\p(\E)$;
such space, which goes under the name of ``enhanced VE bulk space'', reads:
\begin{equation} \label{local space bulk}
\Vp(\E) = \left\{ \vp \in \Vtildep (\E) \,   \left|  \, \int_\E (\vp - \Pinablap \vp) \malpha=0 \text{ for all } \alpha=\npmd+1, \dots, \np         \right.  \right\}.
\end{equation}
Importantly, $\mathbb P_\p(\E) \subseteq \Vp(\E)$; this inclusion guarantees good approximation properties by functions in space \eqref{local space bulk}.

It is possible to compute the $L^2$ projector $\Pizp : \Vp(\E) \rightarrow \mathbb P_\p(\E)$ defined as:
\begin{equation} \label{L2 projector bulk rich}
(\qp, \vp - \Pizp \vp)_{0,\E} = 0 \quad \forall \, \qp \in \mathbb P_\p (\E), \forall \, \vp \in \Vp(\E).
\end{equation}

\begin{remark}
As in the case of face spaces and projectors, the definitions of space \eqref{local space bulk} and of the bulk $H^1$ and $L^2$ orthogonal projectors defined in \eqref{H1 projector bulk} and \eqref{L2 projector bulk rich}, respectively, are independent of
both the face and bulk polynomial basis choices.
\end{remark}

The aforementioned linear functionals forms a unisolvent set of degrees of freedom for space $\Vp(\E)$ defined in \eqref{local space bulk}.
In particular, one has: skeletal dofs given by evaluation at the vertices and ($\p-1$) internal Gau\ss-Lobatto nodes on each edge,
face dofs given by (scaled) face moments \eqref{internal moments faces polyhedron} and (scaled) bulk moments \eqref{internal moments bulk}.

We denote henceforth by:
\begin{equation} \label{dofs and canonical basis}
\left\{ \dof_i \right\}_{i=1}^{\dim(\Vp(\E))}  \quad \text{ and } \quad \left\{ \varphi_i  \right\}_{i=1}^{\dim(\Vp(\E))}
\end{equation}
the set of local degrees of freedom and the local canonical basis, respectively.

\begin{remark} 
Both the face and bulk moments are taken with respect to rather general polynomial bases, which, so far, are required to fulfill assumption \eqref{assumption splitting polynomials} only.
Explicit choices of such polynomial bases are the topic of Section \ref{subsection choice of dofs}.
\end{remark}

\medskip \medskip

The global VE space is obtained by a standard conforming dof coupling and by imposing homogeneous boundary conditions:
\begin{equation} \label{global VE space}
\Vp = \{\vp \in H^1_0(\Omega) \cap \mathcal C^0(\overline \Omega) \mid \vp|_{\overline \E} \in \Vp(\E) \}.
\end{equation}

\medskip

For what concerns the definition of the discrete bilinear form, we follow the VEM gospel and we split it into a sum of local terms:
\begin{equation} \label{global discrete bilinear form}
\ap(\up,\vp) = \sum_{\E \in \taun} \apE(\up,\vp) \quad \forall \, \up,\, \vp \in \Vp,
\end{equation}
which are spit in turn into a sum of two terms, known in the VEM literature as \emph{consistency} and \emph{stabilization} local terms:
\begin{equation} \label{local discrete bilinear form}
\apE(\up,\vp) = \ap(\Pinablap \up, \Pinablap \vp) + \SE((I-\Pinablap ) \up, (I-\Pinablap ) \vp) \quad \forall \, \up,\,\vp \in \Vp(\E).
\end{equation}
Here, $\SE$ is \emph{any} bilinear form satisfying:
\begin{equation} \label{local stabilization}
c_*(\p) \vert \vp \vert^2_{1,\E} \lesssim \SE(\vp, \vp) \lesssim c^*(\p) \vert \vp \vert^2_{1,\E} \quad \forall\, \vp \in \ker (\Pinablap),
\end{equation}
where $c_*(\p)$ and $c^*(\p)$ are two positive constants possibly depending on $\p$;
for an analysis regarding the dependence of $c_*(\p)$ and $c^*(\p)$ on $\p$, we refer to \cite{hpVEMbasic, preprint_hpVEMcorner}.
At the present stage, no explicit bounds are available for the three dimensional case.

\medskip\medskip

We now introduce a family of VEM based on arbitrary stabilizations:
\begin{equation} \label{family VEM}
\begin{cases}
\text{find } \up \in \Vp \text{ s. th.}\\
\ap(\up, \vp) = (\f, \Pizp \vp)_{0,\E} \quad \forall \, \vp \in \Vp
\end{cases}.
\end{equation}
After having introduced the $H^1$ broken Sobolev seminorm associated with polyhedral decomposition $\taun$:
\begin{equation} \label{broken Sobolev seminorm}
\vert \cdot \vert_{1,\taun} = \sqrt{\sum_{\E \in \taun} \vert \cdot \vert^2_ {1,\E}},
\end{equation}
we recall the following abstract error result from \cite{preprint_hpVEMcorner}.
Given $\u$ and $\up$ the solutions of \eqref{weak formulation Poisson problem} and \eqref{family VEM}, respectively,
for any $\upi$ piecewise in $\mathbb P_\p(\E)$ and for any $\uI \in \Vp$, one gets:
\begin{equation} \label{VEM Strang}
\vert \u - \up \vert_{1,\Omega} \lesssim \alpha(\p) \left\{ \vert \u - \upi\vert_{1,\taun} + \vert \u - \uI\vert_{1,\Omega} + \Vert \f - \Pizp \f\Vert_{0,\Omega} \right\},
\end{equation}
where $\alpha(\p)$, the so-called pollution factor in \eqref{VEM Strang}, reads:
\begin{equation} \label{stabilization factor}
\alpha(\p) = \frac{\max(1,c^*(\p))}{\min(1,c_*(\p))}.
\end{equation}

\begin{remark} \label{remark on pollution factor}
Importantly, $\alpha(\p)$ may depend on $\p$, polluting thus the convergence rate of the $\p$ version of VEM. Moreover, $\alpha(\p)$ depends both on the choice of the stabilization and on the choice of the degrees of freedom.
\end{remark}

By requiring proper geometric regularity assumptions on $\taun$, one may prove from \eqref{VEM Strang} $\h$ convergence estimates, see \cite{preprint_VEM3Dbasic}.
Instead, error estimates in terms of $\p$ in three dimensional VEM are not available yet and will be the object of future studies.

%%%%%%%%%%%%%%%%%%%%%%%%%%%%%%%%%%%%%%%%%%%%%%%%%
\subsection{Stabilizations} \label{subsection choice stabilizations}
In this section, we provide a short list of possible local stabilization terms \eqref{local stabilization},
on which we will perform numerical comparisons in Section \ref{section numerical results}.

\begin{enumerate}
\item The first stabilization that we present is somehow the standard one in VEM literature, since it is employed in the pioneering works \cite{VEMvolley, hitchhikersguideVEM} as well as in the majority of VEM works. It reads:
\begin{equation} \label{standard choice stabilization}
\SE_1(\up,\vp) = \hE \sum_{i=1}^{\dim(\Vp(\E))} \dof_i(\up) \, \dof_i(\vp) \quad \forall \, \up,\, \vp \in \ker (\Pinablap).
\end{equation}
We highlight that the presence of factor $\hE$ is used in order to have a stabilization $\SE_1$ which scales like the $H^1$ seminorm for arbitrary diameter $\hE$.

\item
It was observed in \cite{preprint_VEM3Dbasic} that in three dimensions this choice may lead to suboptimal convergence results when employing (moderately) high degrees of accuracy.
This effect was avoided by employing another stabilization, also known as ``\emph{D-recipe}'' stabilization, which is defined as follows.
After that one observes that $\SE(\varphi_i, \varphi_j) = \hE\, \delta_{i,j}$, $\delta_{i,j}$ being the Kronecker delta, completely defines stabilization $\SE_1$, one sets:
\begin{equation} \label{D recipe stabilization}
\SE_2(\varphi_i, \varphi_j) = \max(\hE, \aE(\Pinablap \varphi_i, \Pinablap \varphi_j)) \, \delta_{i,j}.
\end{equation}
It is not hard to understand why stabilization $\SE_2$ is preferable to stabilization $\SE_1$. Indeed, 
it may occur that, if for some reason the energy of $\Pinablap \varphi_i$ is extremely high for most of the basis elements $\varphi_i$,
then the effects of the stabilization $\SE_1$ are negligible in practice; contrarily, by picking stabilization $\SE_2$, one levels off the importance of the consistency and stabilization contributions.

\item
An additional stabilization is obtained by applying the \emph{D-recipe} to stabilization $\SE_1$ only on boundary (i.e. skeleton and face) dofs, neglecting the bulk ones. More precisely, we set:
\begin{equation} \label{D recipe boundary stabilization}
\SE_3(\varphi_i, \varphi_j) =
\begin{cases}
\max(\hE, \aE(\Pinablap \varphi_i, \Pinablap \varphi_j)) \, \delta_{i,j} & \text{if } \varphi_{i} \text{ is a boundary dof}\\
0 & \text{otherwise}
\end{cases}.
\end{equation}

\end{enumerate}

%%%%%%%%%%%%%%%%%%%%%%%%%%%%%%%%%%%%%%%%%%%%%%%%%
\subsection{Polynomial and canonical VEM bases} \label{subsection choice of dofs}
In this section, we discuss some choices for what concerns the polynomial spaces employed in the definition of the face and 
scaled moments introduced in \eqref{internal moments faces polyhedron} and \eqref{internal moments bulk}, respectively, generalizing what done for the two dimensional case in \cite{fetishVEM}.
Importantly, the definition of the face and bulk polynomial spaces are utterly independent.
For the sake of simplicity, we define one type of polynomial basis on all faces and one type of polynomial basis in every polyhedron.

We extensively use the two natural bijections $\mathbb N ^2 \leftrightarrow \mathbb N$ and $\mathbb N ^3 \leftrightarrow \mathbb N$ defined as:
\begin{equation} \label{bijection 2D}
(0,0) \leftrightarrow 1, \; (1,0) \leftrightarrow 2, \; (0,1) \leftrightarrow 3 \dots
\end{equation}
and
\begin{equation} \label{bijection 3D}
(0,0,0) \leftrightarrow 1, \; (1,0,0) \leftrightarrow 2, \; (0,1,0) \leftrightarrow 3, \; (0,0,1) \leftrightarrow 4 \dots
\end{equation}

\medskip

We start with polynomial bases on the faces. Given $\F \in \mathcal F_n$ a face, we define by $\xFbold = (\xF, \yF)$ and $\hF$ its barycenter and diameter, respectively.
Note that $\xFbold$ is written with respect to the local coordinates system on face $\F$.
Our first choice reads:
\begin{equation} \label{monomials face}
\begin{aligned}
\malpha^\F (\xbold) 	& = \left( \frac{\xbold - \xFbold}{\hF} \right)^{\boldalpha} \\
				& = \left( \frac{x-\xF}{\hF} \right)^{\alpha_1}  \left( \frac{y-\yF}{\hF} \right)^{\alpha_2} \; \forall\, \boldalpha = (\alpha_1,\alpha_2) \in \mathbb N^2,\, \vert \boldalpha \vert=0,\dots, \p-2.\\
\end{aligned}
\end{equation}
A second choice is given by $\{\mbaralpha^\F\}_{\alpha=1}^{\npmd^\F}$ which can be obtained from \eqref{monomials face} via a stable $L^2(\F)$ orthonormalizing Gram-Schmidt process, e.g. the one presented in \cite{BassiBottiColomboDipietroTesini}.
This choice was in fact already performed on polygons in \cite{fetishVEM}.

Next, we introduce polynomial bases in the bulk. Given $\E \in \taun$ polyhedron, we define $\xEbold$ and $\hE$ its barycenter and diameter, respectively.
Note that $\xEbold$ is written with respect to the global coordinate system of $\mathbb R^3$.
Our first choice reads:
\begin{equation} \label{monomials bulk}
\begin{aligned}
\malpha (\xbold) 	 & = \left( \frac{\xbold - \xEbold}{\hE} \right)^{\boldalpha} 	\\
			& = \left( \frac{x-\xE}{\hE} \right)^{\alpha_1}  \left( \frac{y-\yE}{\hE} \right)^{\alpha_2} \left( \frac{z-\zE}{\hE} \right)^{\alpha_3} \; \forall\, \boldalpha = (\alpha_1,\alpha_2, \alpha_3) \in \mathbb N^3,\,
																													\vert \boldalpha \vert=0,\dots, \p-2.\\
\end{aligned}
\end{equation}
A second choice is given by $\{\mbaralpha\}_{\alpha=1}^{\npmd}$ which can again be obtained from \eqref{monomials bulk} via a stable $L^2(\E)$ orthonormalizing Gram-Schmidt process, see \cite{BassiBottiColomboDipietroTesini}.

In the numerical tests presented in the forthcoming Section \ref{section numerical results}, we employ the following combinations of polynomial bases:
\begin{itemize}
\item \textbf{``\standard choice''}: monomials $\{\malpha^\F\}_{\alpha=1}^{\npmd^\F}$ on all faces $\F$ of $\E$ and monomials $\{ \malpha \}_{\alpha=1}^{\npmd}$ in the bulk;
this choice is the standard one in three dimensional VEM, see e.g. \cite{preprint_VEM3Dbasic, equivalentprojectorsforVEM};
the implementation details when employing such basis are already known from \cite{hitchhikersguideVEM};
\item \textbf{``\orthogonal choice''}: $L^2(\F)$ orthonormal polynomials $\{\mbaralpha^\F\}_{\alpha=1}^{\npmd^\F}$ on all faces $\F$ of $\E$ and $L^2(\E)$ orthonormal polynomials $\{ \mbaralpha \}_{\alpha=1}^{\npmd}$ in the bulk;
some implementation details of this new approach are discussed in Appendix \ref{subappendix pure fetish};
\item \textbf{``\hybrid choice''}: monomials $\{\malpha^\F\}_{\alpha=1}^{\npmd^\F}$ on all faces $\F$ of $\E$ and $L^2(\F)$ orthonormal polynomials $\{ \mbaralpha \}_{\alpha=1}^{\npmd}$ in the bulk;
%this choice, which on ``non-exotic''  meshes performs slightly worse than the ``orthonormal choice'', seems to be particularly advisable when employing polyhedral meshes with small faces;
some implementation details of this new approach are the topic of Appendix \ref{subappendix hybrid fetish}.
\end{itemize}

As a byproduct we remark that in principle one could use a sort of mix between the ``\orthogonal choice'' and the ``\hybrid choice'', by picking the ``orthonormal'' one on some faces only.
For the sake of an easy implementation of the method and also for the sake of a more straightforward presentation, we stick to the case of uniform choice on all faces.

% --------------------------------------------------------------------------------------------------------------------------------------------------------------------------------------------------------------------------------------------------------------------
\section {Numerical results} \label{section numerical results}
% --------------------------------------------------------------------------------------------------------------------------------------------------------------------------------------------------------------------------------------------------------------------
This section is devoted to numerically compare the choices of the stabilizations introduced in Section \ref{subsection choice stabilizations} and of the face/bulk moments introduced in Section \ref{subsection choice of dofs}.

When studying the error of the method, owing to the fact that functions in the VE spaces are not known explicitly neither on the faces nor in the bulk of each element but only on the skeleton of the mesh, we compute the following couple of relative errors:
\begin{equation} \label{computed errors}
\frac{\vert \u - \Pinablap \up \vert_{1,\taun}}{\vert \u \vert_{1,\Omega}}, \quad \quad \quad \quad \frac{\Vert \u - \Pizp \up \Vert_{0,\Omega}}{\Vert \u \Vert_{0,\Omega}},
\end{equation}
where $\vert \cdot \vert_{1,\taun}$ is the $H^1$ broken Sobolev seminorm introduced in \eqref{broken Sobolev seminorm}, $\Pinablap$ is defined in \eqref{H1 projector bulk} and
$\Pizp$ is defined in \eqref{L2 projector bulk rich}, $\u$ is the exact solution of problem \eqref{weak formulation Poisson problem} and $\up$ is the solution of VEM \eqref{family VEM}.

In the forthcoming sections, we perform a number of numerical tests
by taking the standard unit cube $\Omega:=[0,\,1]^3$ as physical domain and 
by considering as solutions of problem \eqref{weak formulation Poisson problem} the two test cases defined as:
\begin{subequations}
\begin{equation} \label{solution 1}
\u_1(x,y,z) = \sin(\pi\, x) \sin(\pi \, y) \sin(\pi \, z),
\end{equation}
\begin{equation} \label{solution 2}
\u_2(x,y,z)  = 1+x+y+z.
\end{equation}
\end{subequations}
We underline that $u_1$ is analytic while $u_2$ is a polynomial of degree $1$.
Hence, the method should return up to machine precision the polynomial solution $\u_2$, see \cite{hitchhikersguideVEM,preprint_VEM3Dbasic}. 

We carry out numerical tests by employing three different types of polyhedral decomposition:
\begin{itemize}
\item meshes made of structured cubes; we refer in the following to a mesh of this sort as ``\texttt{cube mesh}'';
\item meshes obtained by a Voronoi tessellation of sets of points randomly chosen inside $\Omega$ optimized via Lloyd's algorithm~\cite{dufabergunzburgerVoronoi};
we refer in the following to a mesh of this sort as ``\texttt{Voronoi mesh}'';
\item meshes obtained by a Voronoi tessellation of sets of points randomly chosen inside $\Omega$; we refer in the following to a mesh of this sort as  ``\texttt{rand mesh}''.
\end{itemize}
We point out that meshes of type ``\texttt{rand mesh}'' contain distorted elements which  are instrumental for severely testing the robustness of VEM with respect to mesh distortion.

%; thanks to the particular construction of the discrete bilinear form, see \eqref{global discrete bilinear form} and \eqref{local discrete bilinear form},
%this method should return up to machine precision the polynomial solution. Importantly, this leads to the following consideration: the error computed in the approximation of the patch test $\u_2$ is the ``litmus test''
%to understand what is the impact of the condition number when solving the linear system; in fact, the error on the patch test should be zero in principle, but in machine precision it is not and it grows as the condition number of the stiffness matrix grows.

\medskip

The outline of this section follows.
Firstly, in Section \ref{subsection nr: reasons suboptimality}, we highlight some of the possible reasons for which the method can return unexpected/wrong error convergence slopes.
In Section \ref{subsection nr: stabilizations}, we fix the choice of the polynomial bases dual to face \eqref{internal moments faces polyhedron} and bulk \eqref{internal moments bulk} moments
to the ``standard choice'' presented in Section \ref{subsection choice of dofs} and we compare the effects of the three stabilizations presented in Section \ref{subsection choice stabilizations} of the $\p$ version of VEM.
Next, in Section \ref{subsection nr: p version}, we consider again the $\p$ version of VEM and we investigate its sensibility by fixing the stabilizations and by varying
the choice of the polynomials bases dual to face/bulk  moments.
Instead, in Section \ref{subsection nr: collapsing}, we study the effects due to choice of the stabilizations and the face/bulk moments when employing polyhedral meshes with extremely degenerate elements.
Finally, in Section \ref{subsection nr: h VEM}, we perform some tests on the $\h$ version of the method comparing again the effects of the choice of face/bulk moments.

We add at the end of each section a condensed summary, highlighting therein in short the conclusions of each set of numerical experiments.

\paragraph*{Notation employed in Section \ref{section numerical results}.}
In order to manage the large (and somehow cumbersome) amount of data and not to jeopardize the understanding of the reader, we fix here once and for all some notations.

We will test the method with solutions $\u_1$ and $\u_2$ defined in \eqref{solution 1} and \eqref{solution 2},
employing the ``standard'', the ``orthogonal'' and the ``hybrid choices'' described in Section \ref{subsection choice of dofs}
and employing stabilizations $\SE_1$, $\SE_2$ and $\SE_3$ defined in \eqref{standard choice stabilization}, \eqref{D recipe stabilization} and \eqref{D recipe boundary stabilization}, respectively.
Moreover, by $H^1$ and $L^2$ error, we denote those introduced in \eqref{computed errors}.
\par

%%%%%%%%%%%%%%%%%%%%%%%%%%%%%%%%%%%%%%%%%%%%%%%%%
\subsection{Possible reasons for suboptimality in the error convergence curves} \label{subsection nr: reasons suboptimality}
As sometimes happens in scientific computing and more specifically in the numerical approximation of PDEs,
it may occur in the VEM framework that the numerical %tests do not return the error convergence curves that are expected from the abstract error result \eqref{VEM Strang}
%and from best approximation results.
performances of the method suffer a lack of accuracy/convergence.
This has been already observed for the $\p$ version of two and three dimensional VEM in \cite{fetishVEM} and \cite{preprint_VEM3Dbasic}, respectively.

\begin{figure}[h] 
\begin{center}
\begin{tikzpicture}[auto, node distance=4cm,>=latex']
	\node[rect](input){\footnotesize{choice dofs}};
	\node[rect,  above of = input](cond){\footnotesize{condition number}};
	\path[line] (input) -- (cond);
	\node[rect,  right of = cond](stab){\footnotesize{choice stabilization}};
	\path[line] (stab) -- (cond);
	\node[rect,  below of = stab](pollution){\footnotesize{pollution factor \eqref{stabilization factor}}};
	\path[line] (stab) -- (pollution);
	\node[rect,  right of = pollution](error){\footnotesize{error slopes}};
	\path[line] (pollution) -- (error);
	\path[line] (input) -- (pollution);
	\path[line] (cond) -- (error);
\end{tikzpicture}
\end{center}
\caption{Scheme summarizing the reasons of possible wrong behaviours of error slopes.} \label{scheme for reasons suboptimality}
\end{figure}
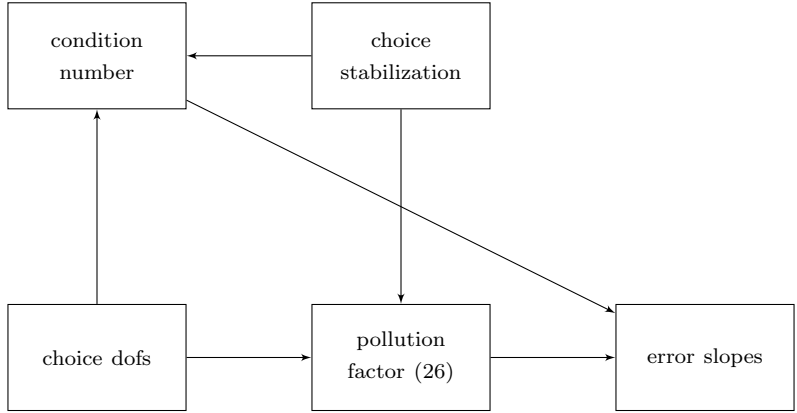

In the forthcoming sections, we will behold additional suboptimal/wrong results when considering ``not shrewd'' stabilizations and face/bulk moments.
%We wonder which are the reasons for such unexpected behaviour. For this purpose, we list here two of them.
We highlight here two among the possible reasons for such unexpected behaviours and, in Figure \ref{scheme for reasons suboptimality},
we depict a scheme summarizing what are their effects on the performances of the method.

\begin{itemize}
\item
The first one is the \textbf{condition number} of the stiffness matrix. A possible way to understand its impact on the method
is to study the error slopes on the so-called patch test, i.e. on an exact solution which is polynomial on the complete computational domain.
Indeed, owing to the particular construction of the discrete bilinear form \eqref{global discrete bilinear form} and \eqref{local discrete bilinear form}, the method should return in this case the exact solution up to machine precision.
In practice, the error on the patch test grows as the condition number of the stiffness matrix grows.
We emphasize that such condition number depends both on the choice of the canonical basis, see Section \ref{subsection choice stabilizations}, and on the choice of the stabilization, see Section \ref{subsection choice of dofs}.
\item
The \textbf{stabilization} of the method, which is the second reason for possible suboptimal/wrong behaviour of the method, has effects also on the error estimates though the pollution factor $\alpha(\p)$ defined in \eqref{stabilization factor}
which pops up in the theoretical convergence analysis of the method \eqref{VEM Strang}.
We underline that the behavior of the pollution factor depends also on the choice of the canonical basis  of the VE space as explained e.g. in \cite[Section 6.4]{hpVEMbasic}.
In particular, this means that ``clever'' choices of the stabilization do not automatically entails optimal convergence slopes;
one has to pick a ``clever'' choice of the canonical basis as well.
\end{itemize}

%In Figure \ref{scheme for reasons suboptimality}, we depict a scheme summarizing the reasons of possible unexpected behaviours of error slopes.

%%%%%%%%%%%%%%%%%%%%%%%%%%%%%%%%%%%%%%%%%%%%%%%%%
\subsection{Numerical results: the effects of the stabilization} \label{subsection nr: stabilizations}
In this section, we investigate the effects of the choice of the stabilization in the $\p$ version of three dimensional VEM keeping fixed the choice of the degrees of freedom.
This is a step-forward with respect to what was presented in \cite{preprint_VEM3Dbasic}.

\begin{figure}  [h]
\centering
\subfigure {\includegraphics [angle=0, width=0.49\textwidth]{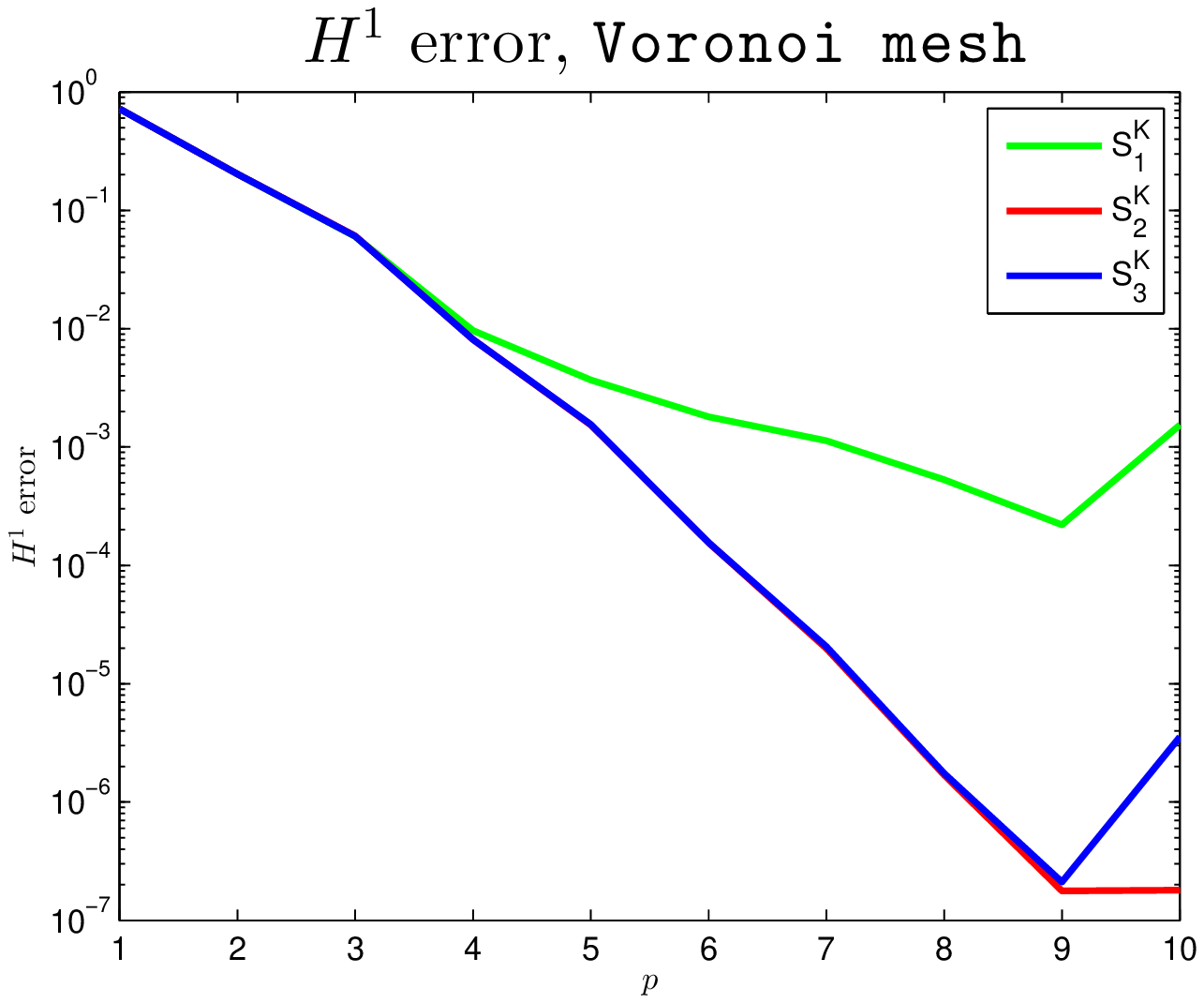}}
\subfigure {\includegraphics [angle=0, width=0.49\textwidth]{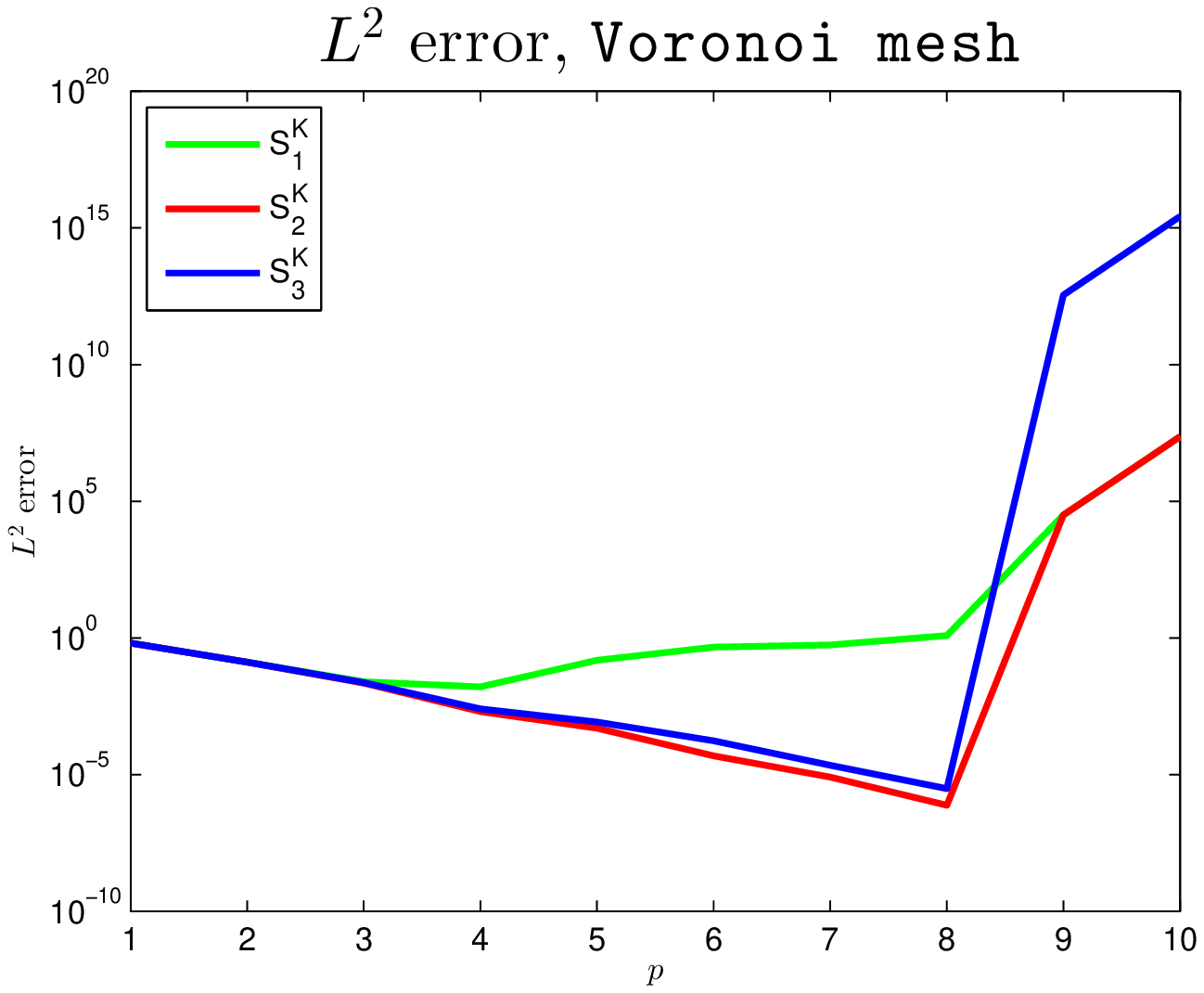}}
\subfigure {\includegraphics [angle=0, width=0.49\textwidth]{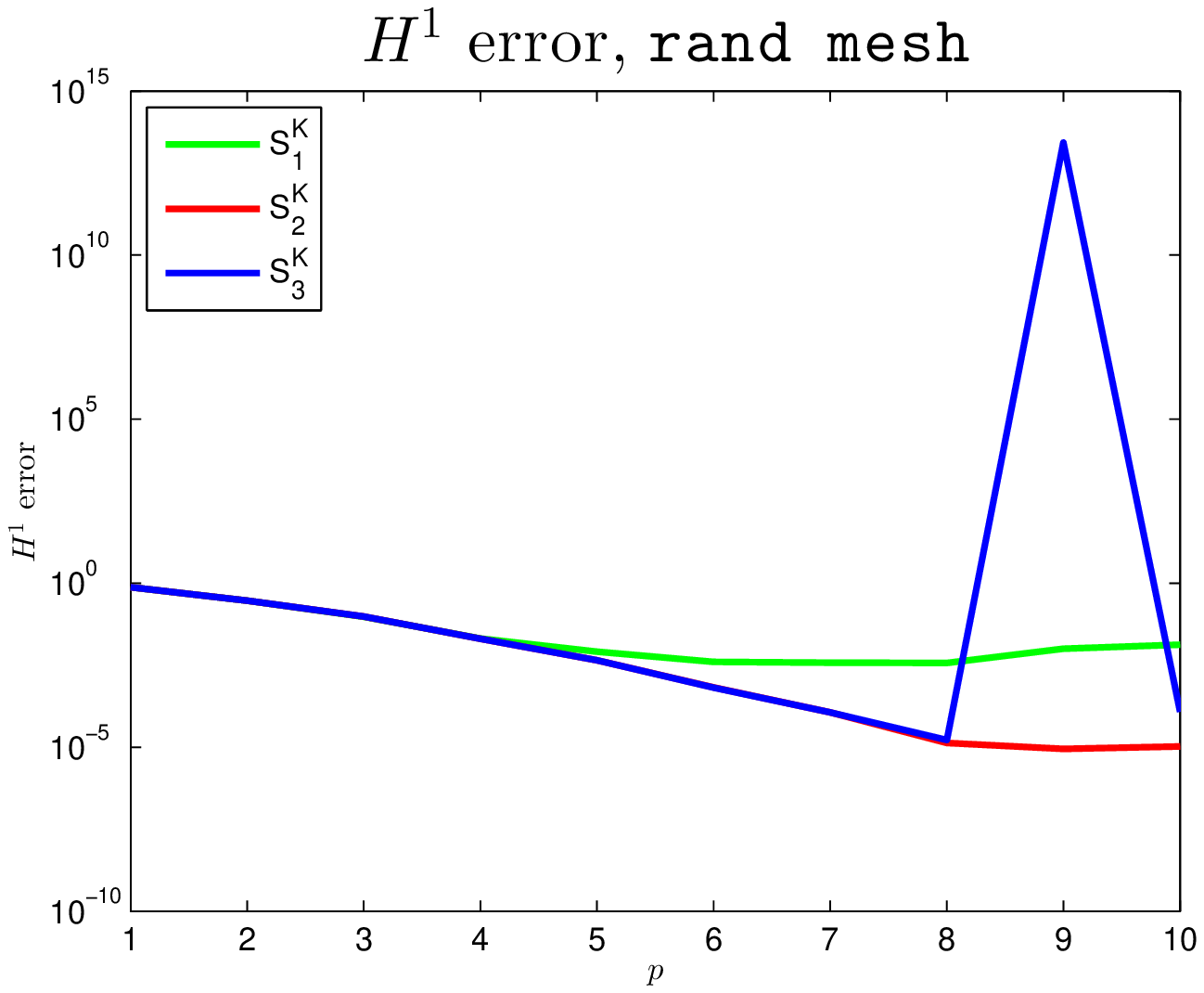}}
\subfigure {\includegraphics [angle=0, width=0.49\textwidth]{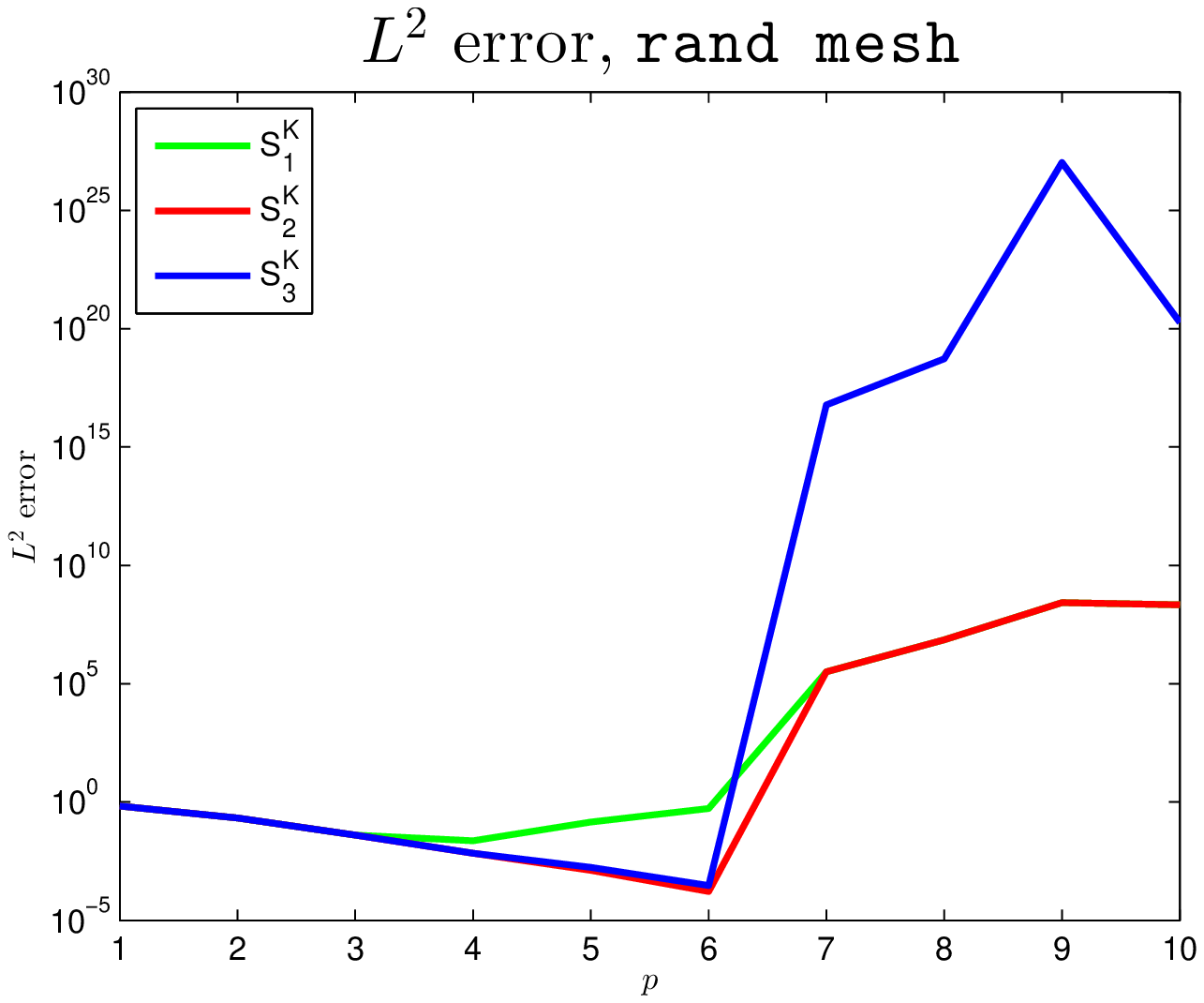}}
\caption{$\p$ version of VEM. Exact solution $\u_1$. Stabilizations employed: $\SE_1$, $\SE_2$ and $\SE_3$.
Face/bulk moments employed: ``standard choice''.
Up-left: $H^1$ error on a \texttt{Voronoi mesh}. Up-right: $L^2$ error on a \texttt{Voronoi mesh}.
Down-left: $H^1$ error on a \texttt{rand mesh}. Down-right: $L^2$ error on a \texttt{rand mesh}.}
\label{figure comparison stabilizations 2 meshes}
\end{figure}

In particular, we consider the ``standard choice'' 
and stabilizations $\SE_1$, $\SE_2$ and $\SE_3$. As a test case, we consider the analytic function $\u_1$
and we consider two meshes, namely a \texttt{Voronoi mesh} and a \texttt{rand mesh}.
In Figure \ref{figure comparison stabilizations 2 meshes}, we show the convergence of the $H^1$ and $L^2$ errors defined in \eqref{computed errors} on a \texttt{Voronoi mesh} and on a \texttt{rand mesh}.

What we observe is that with no doubts the best performances are those related to stabilization $\SE_2$. The standard stabilization leads to suboptimal convergence even for moderately low degrees of accuracy;
nonetheless, it prevents the error slopes to suddenly blow up as it happens for stabilization $\SE_3$.
Having said this, in the forthcoming sections, we present a number of numerical tests dropping stabilization $\SE_3$.

\paragraph*{Summary:} in $\p$ VEM, employ stabilization $\SE_2$, avoid stabilizations $\SE_3$; better to avoid $\SE_1$.

%%%%%%%%%%%%%%%%%%%%%%%%%%%%%%%%%%%%%%%%%%%%%%%%%
\subsection{Numerical results: the effects of the face/bulk moments} \label{subsection nr: p version}
This section is devoted to understand the impact that the choice of face \eqref{internal moments faces polyhedron} and bulk \eqref{internal moments bulk} moments 
has on the convergence of the $\p$ version of VEM.

We assume to use stabilizations $\SE_1$ and $\SE_2$ and we compare the effects on the performances of the method employing the ``standard'', ``orthogonal'' and ``hybrid choices''.

\begin{figure}  [h]
\centering
\subfigure {\includegraphics [angle=0, width=0.49\textwidth]{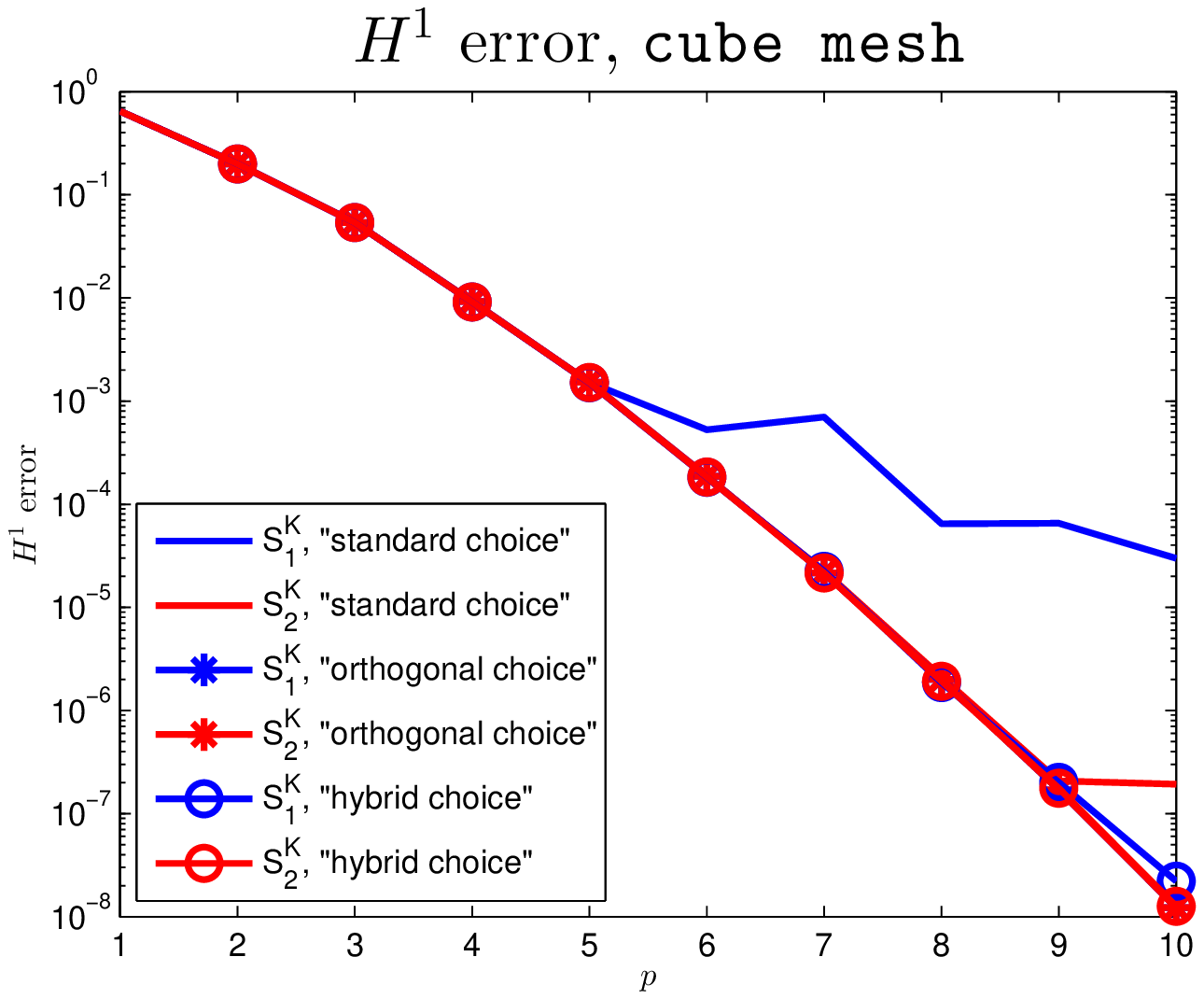}}
\subfigure {\includegraphics [angle=0, width=0.49\textwidth]{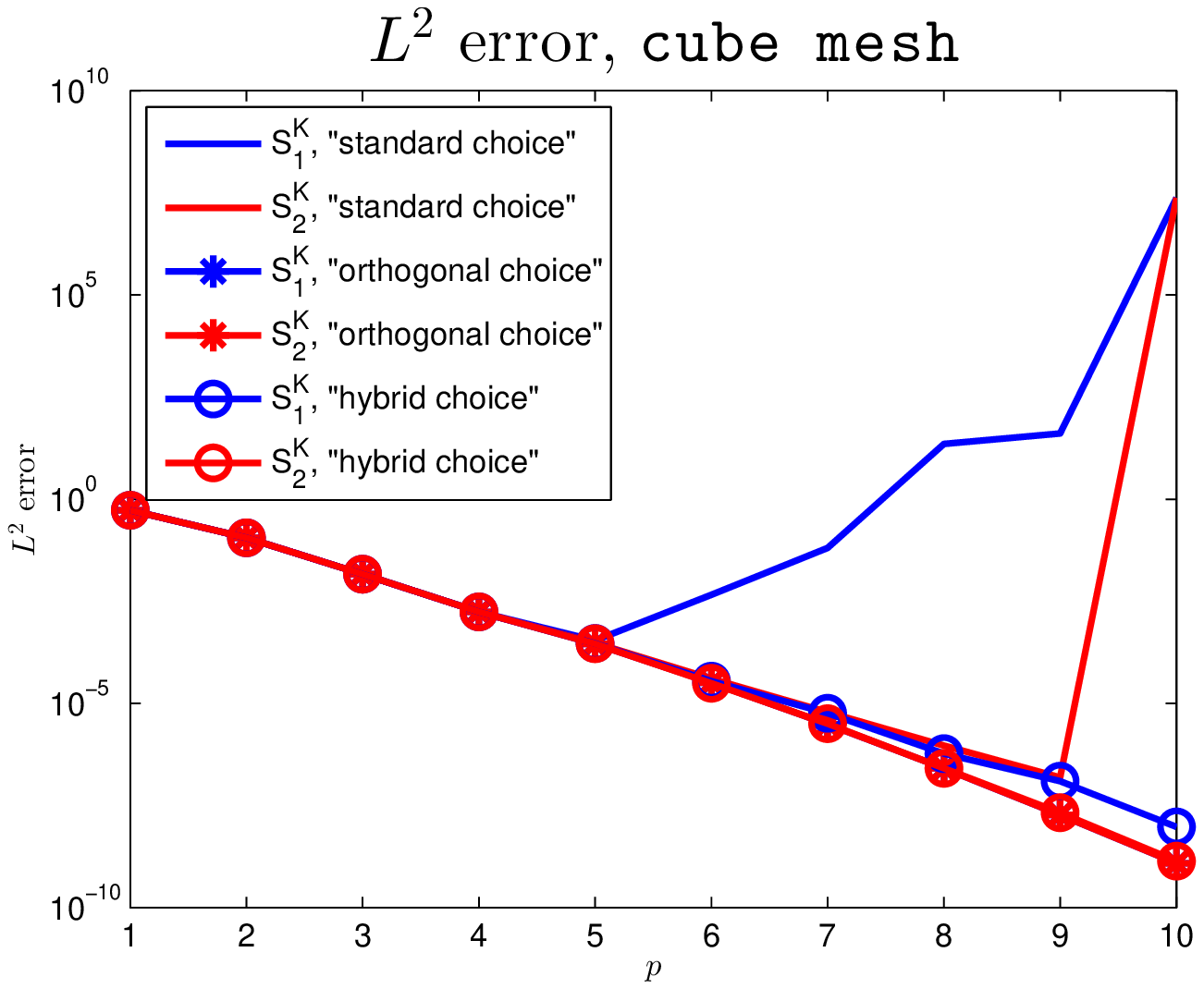}}
\caption{$\p$ version of VEM. Exact solution $\u_1$. Stabilizations employed: $\SE_1$ and $\SE_2$.
Face/bulk moments employed: ``standard choice'', ``orthogonal choice'' and ``hybrid choice''.
\texttt{Cube mesh}. Left: $H^1$ error. Right: $L^2$ error.}
\label{figure changing polynomial bases analytic solution cube}
\end{figure}

\begin{figure}  [h]
\centering
\subfigure {\includegraphics [angle=0, width=0.49\textwidth]{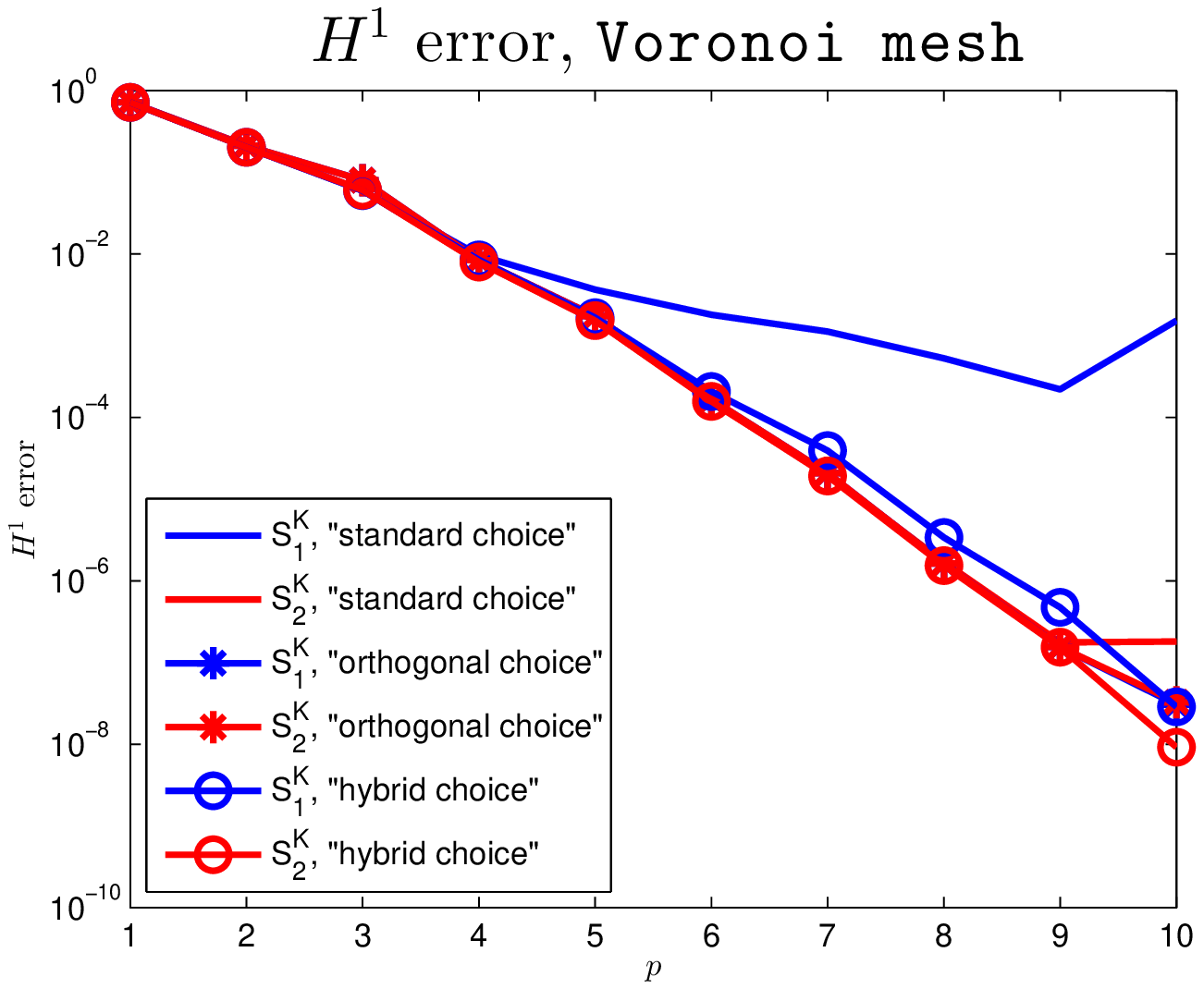}}
\subfigure {\includegraphics [angle=0, width=0.49\textwidth]{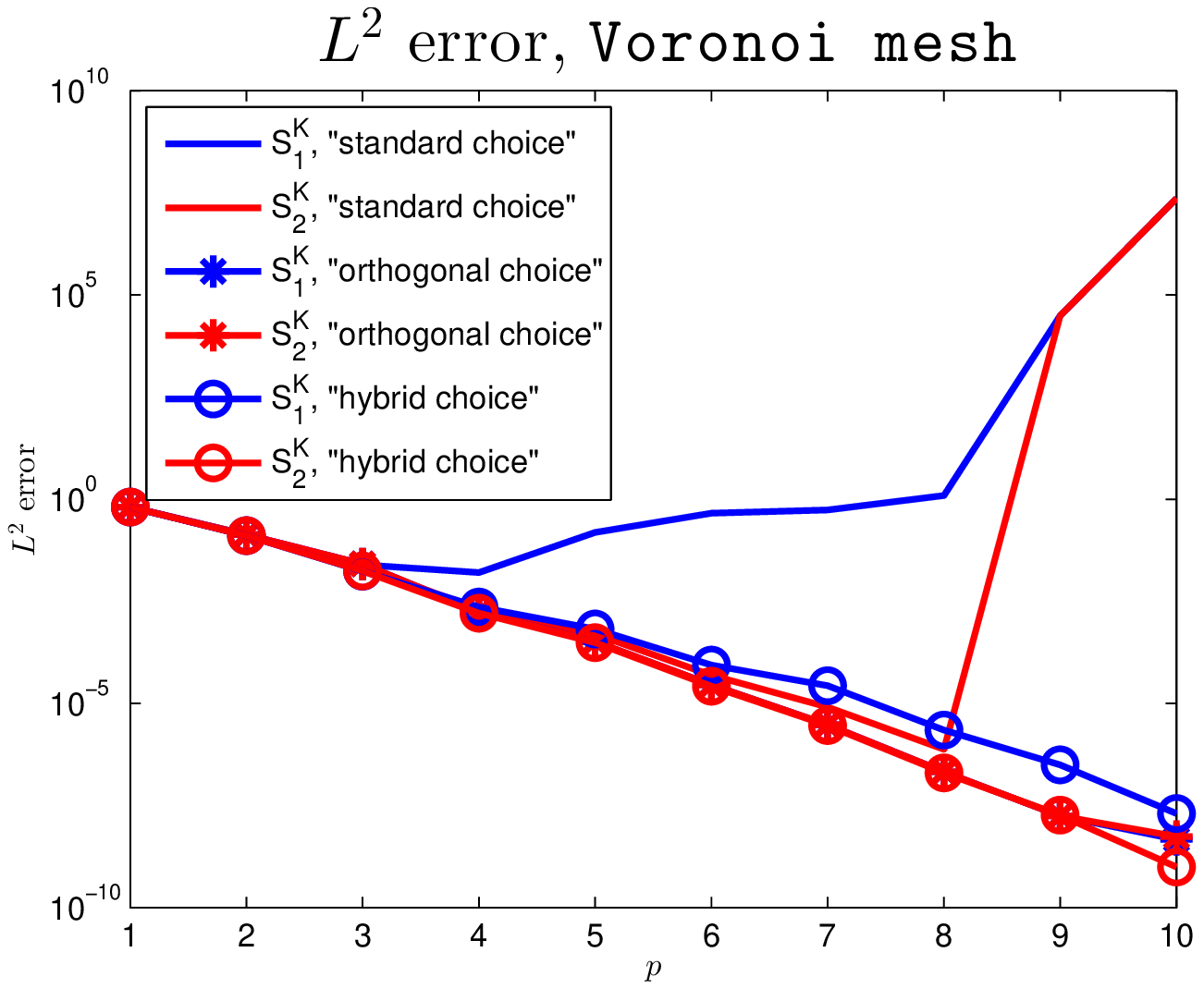}}
\caption{$\p$ version of VEM. Exact solution $\u_1$. Stabilizations employed: $\SE_1$ and $\SE_2$.
Face/bulk moments employed: ``standard choice'', ``orthogonal choice'' and ``hybrid choice''.
\texttt{Voronoi mesh}. Left: $H^1$ error. Right: $L^2$ error.}
\label{figure changing polynomial bases analytic solution Voronoi}
\end{figure}

\begin{figure}  [h]
\centering
\subfigure {\includegraphics [angle=0, width=0.49\textwidth]{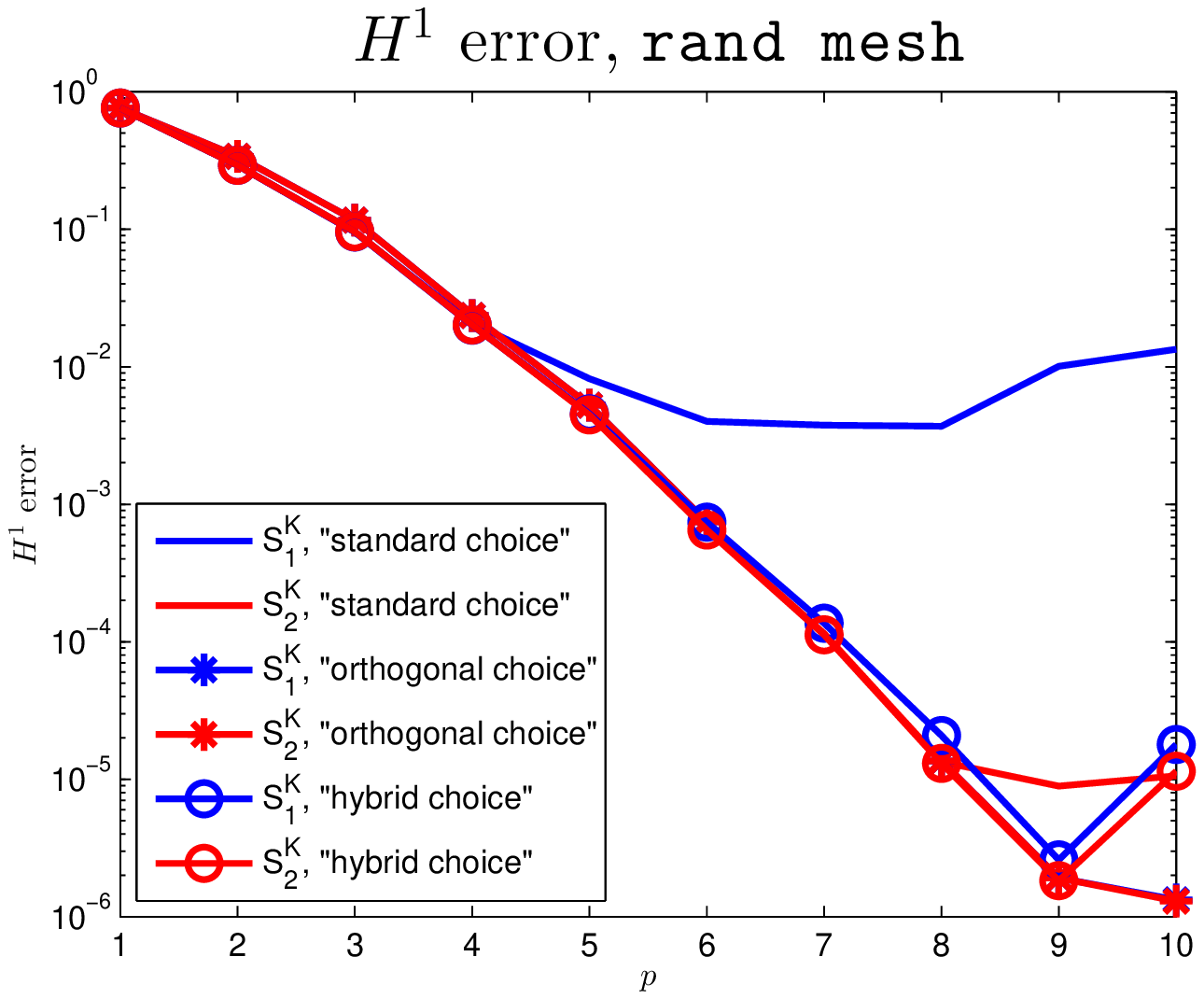}}
\subfigure {\includegraphics [angle=0, width=0.49\textwidth]{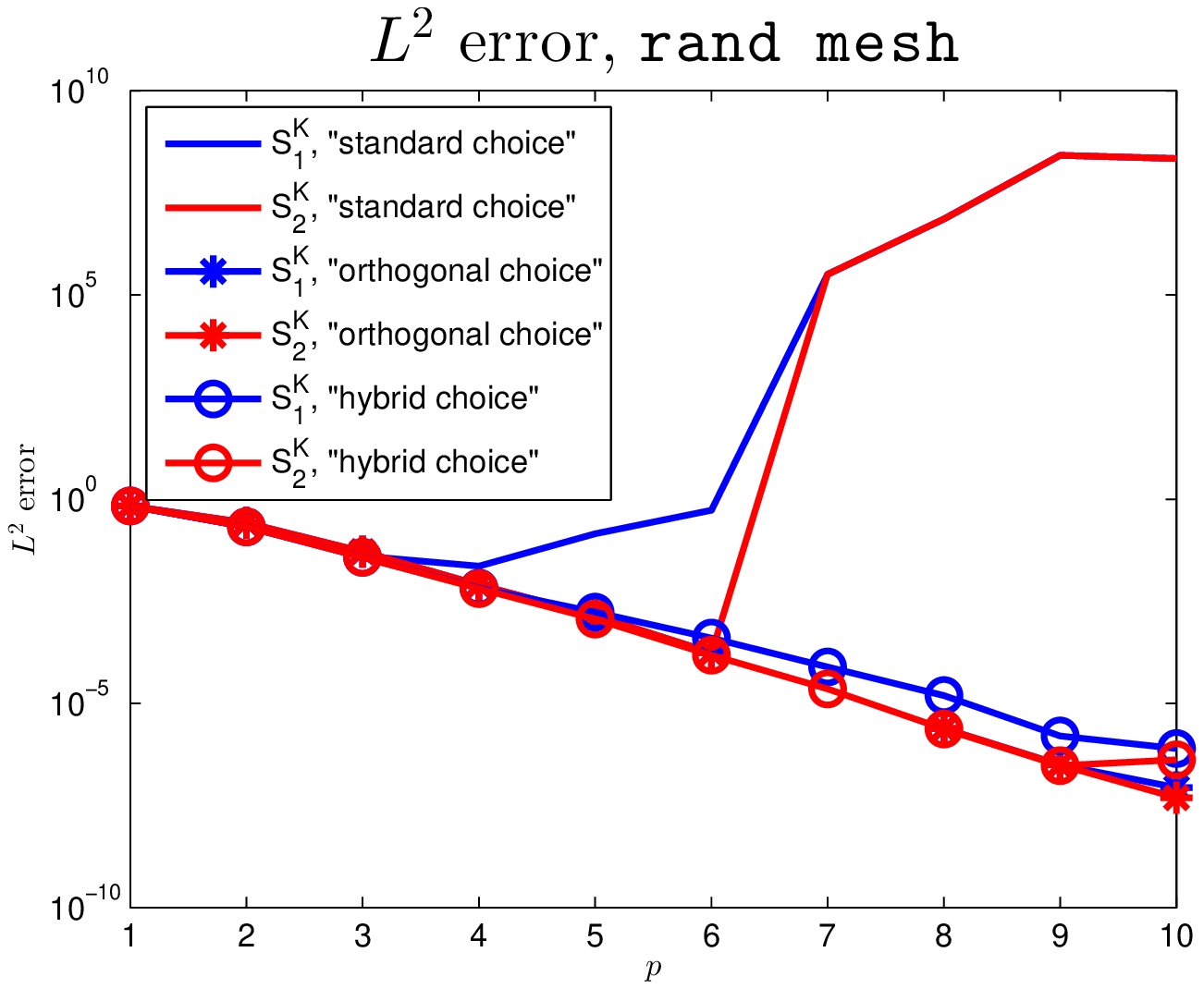}}
\caption{$\p$ version of VEM. Exact solution $\u_1$. Stabilizations employed: $\SE_1$ and $\SE_2$.
Face/bulk moments employed: ``standard choice'', ``orthogonal choice'' and ``hybrid choice''.
\texttt{Rand mesh}. Left: $H^1$ error. Right: $L^2$ error.}
\label{figure changing polynomial bases analytic solution rand}
\end{figure}

We begin by applying the method on the model problem with exact solution $\u_1$ and by employing a \texttt{cube mesh} in Figure \ref{figure changing polynomial bases analytic solution cube},
a \texttt{Voronoi mesh} in Figure \ref{figure changing polynomial bases analytic solution Voronoi} and a \texttt{rand mesh} in Figure \ref{figure changing polynomial bases analytic solution rand}.
We check both the $H^1$ and the $L^2$ errors.

What we observe in Figures \ref{figure changing polynomial bases analytic solution cube}, \ref{figure changing polynomial bases analytic solution Voronoi} and \ref{figure changing polynomial bases analytic solution rand}
is that when employing a very regular mesh (\texttt{cube mesh}) the ``standard choice'' suffers a lack of convergence even employing stabilization $\SE_2$, which is the most robust among those we presented;
on the other hand, by employing both the ``orthogonal choice'' and the ``hybrid choice'', the method converges without any loss.
Analogous comments hold true when employing a \texttt{Voronoi mesh}, which is less regular than the cube one,
although it seems that the $H^1$ and $L^2$ errors on the \texttt{Voronoi mesh} employing the ``orthogonal choice'' and stabilization $\SE_2$ starts to grow for $\p=10$.

Interestingly, analogous results are valid also when employing a much less regular mesh (with small faces, small edges) such as the \texttt{rand mesh}.
The error slopes when employing the ``orthogonal choice'' and the ``hybrid choice'' are identical in practice up to $\p=9$;
for $\p=10$, the ``orthogonal choice'' performs slightly better.

At the end of the day, we can say that in the $\p$ version of VEM the ``orthogonal'' and the ``hybrid choices'' are comparable and show extremely satisfactory results.

\medskip

At this point, we wonder the reasons for which the ``standard choice'' leads to suboptimal convergence.
As discussed in Section \ref{subsection nr: reasons suboptimality}, one possible reason is the choice of the stabilization along with its effect on the pollution factor $\alpha(\p)$ defined in \eqref{stabilization factor};
another one is the condition number of the stiffness matrix.

%Among all the possible causes we highlight two of them which we deem could influence more the behaviour of the error slopes:
%\begin{itemize}
%\item the first one is that the stabilization employed behaves worse in this case than when employing for instance the ``orthogonal choice''
%(one has to take into account that the pollution factor defined in \eqref{stabilization factor} depends on the choice of the stabilization but also on the choice of the canonical basis,
%which in turn depends on the choice of the face/bulk moments); this can be observed in Section \ref{subsection nr: stabilizations} and in \cite{preprint_VEM3Dbasic} as well;
%\item the second one is the condition number of the stiffness matrix.
%\end{itemize}

\begin{figure}  [h]
\centering
\subfigure {\includegraphics [angle=0, width=0.49\textwidth]{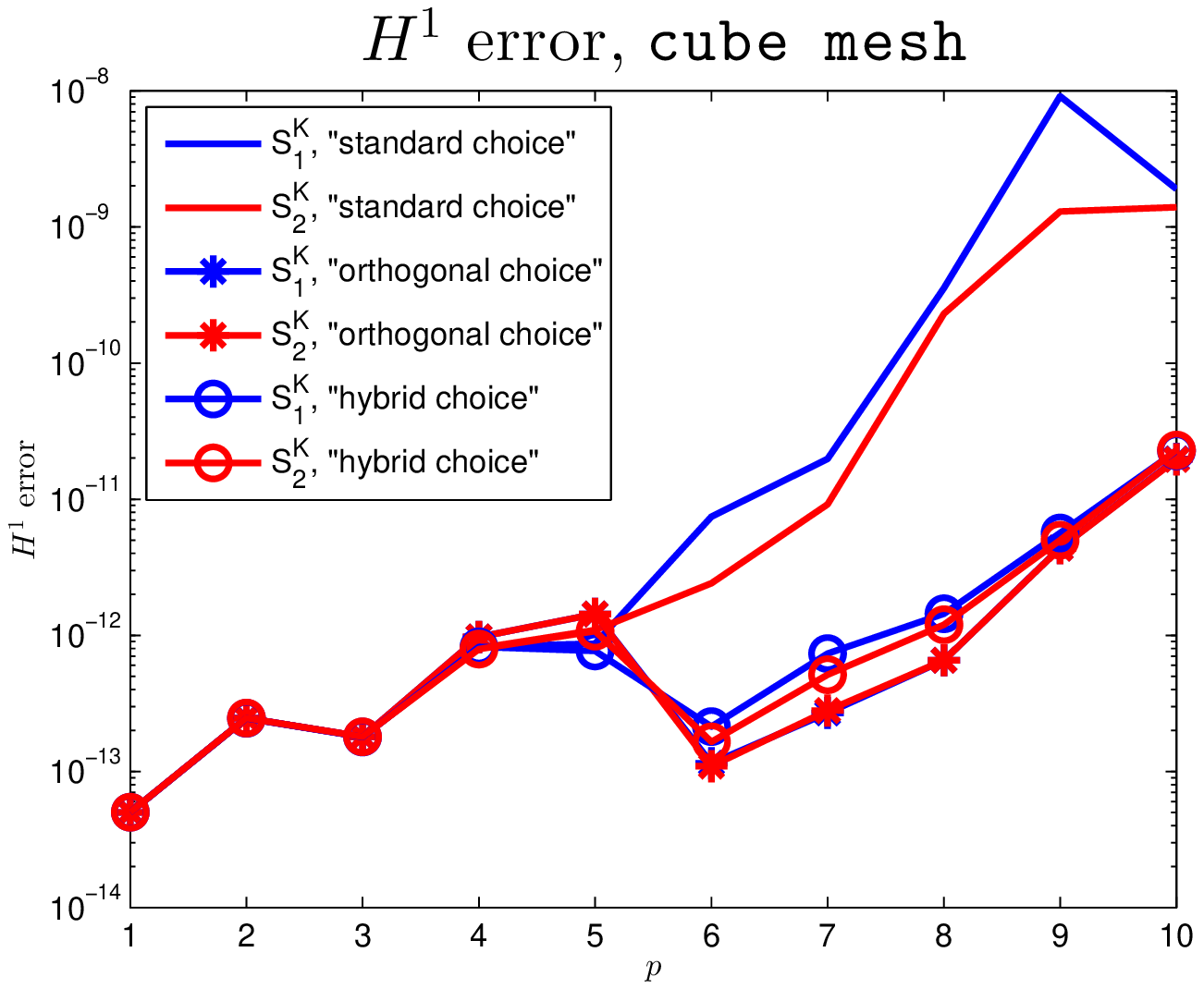}}
\subfigure {\includegraphics [angle=0, width=0.49\textwidth]{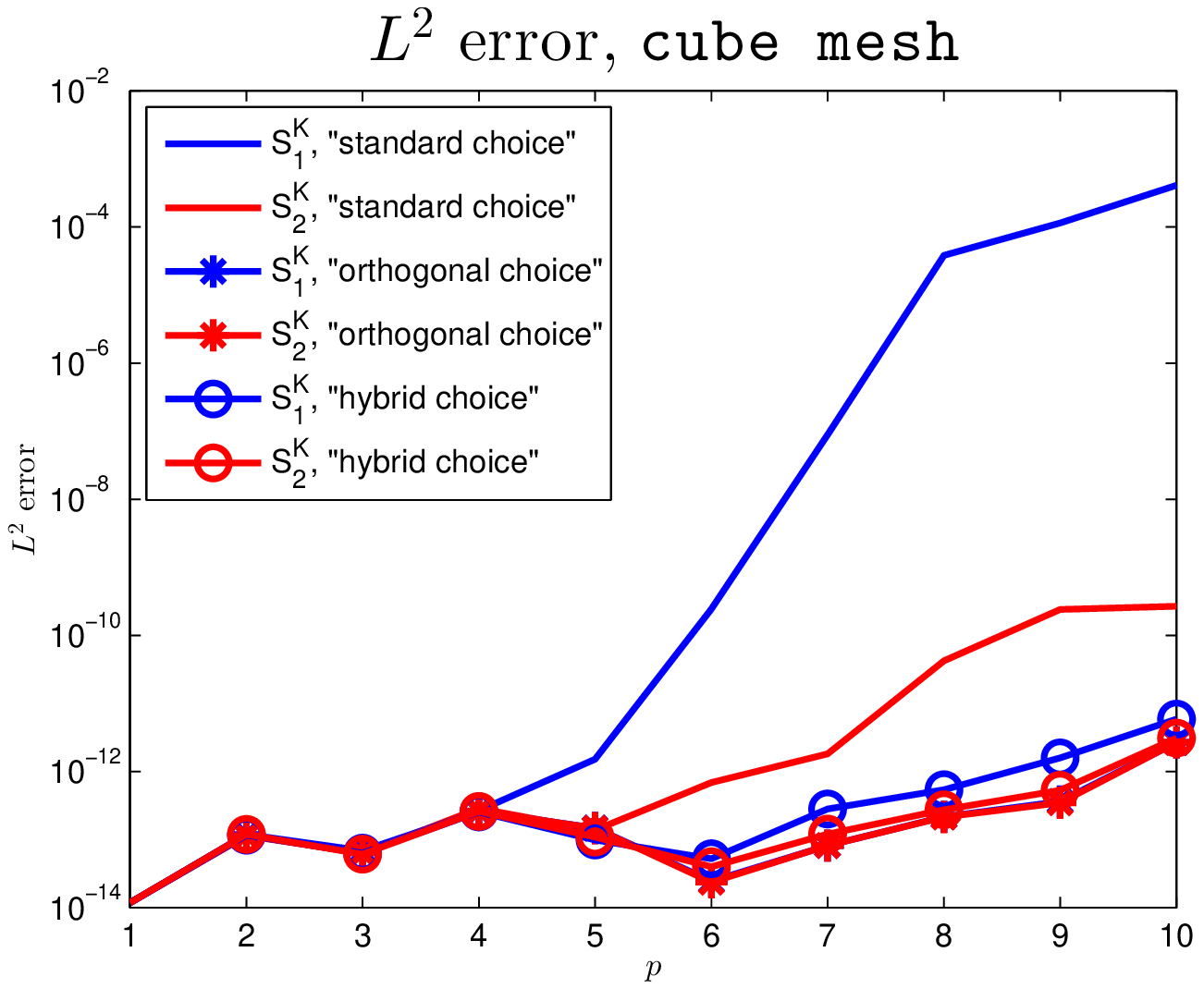}}
\caption{$\p$ version of VEM. Exact solution $\u_2$. Stabilizations employed: $\SE_1$ and $\SE_2$.
Face/bulk moments employed: ``standard choice'', ``orthogonal choice'' and ``hybrid choice''.
\texttt{Cube mesh}. Left: $H^1$ error. Right: $L^2$ error.}
\label{figure changing polynomial bases patch test cube}
\end{figure}

\begin{figure}  [h]
\centering
\subfigure {\includegraphics [angle=0, width=0.49\textwidth]{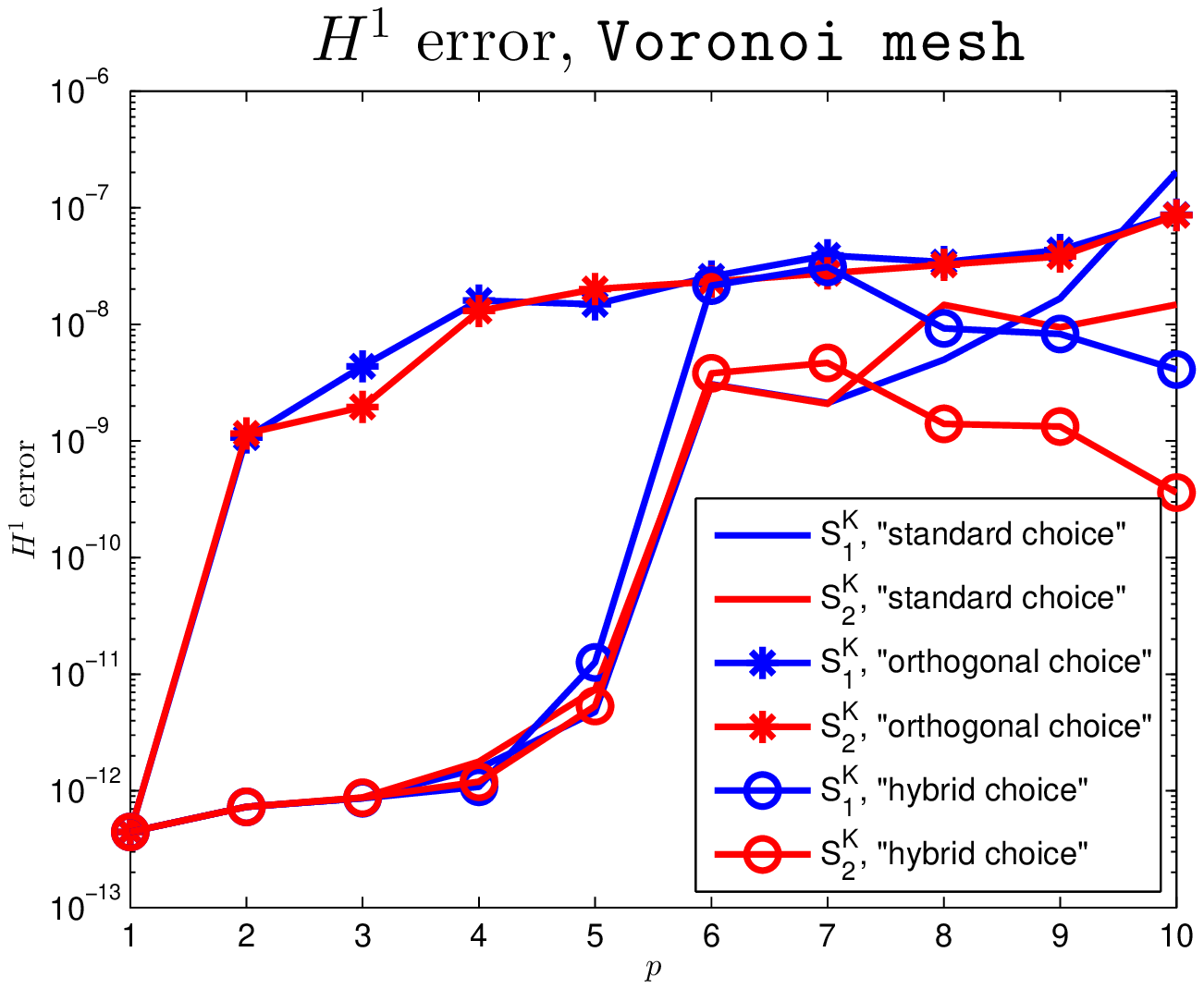}}
\subfigure {\includegraphics [angle=0, width=0.49\textwidth]{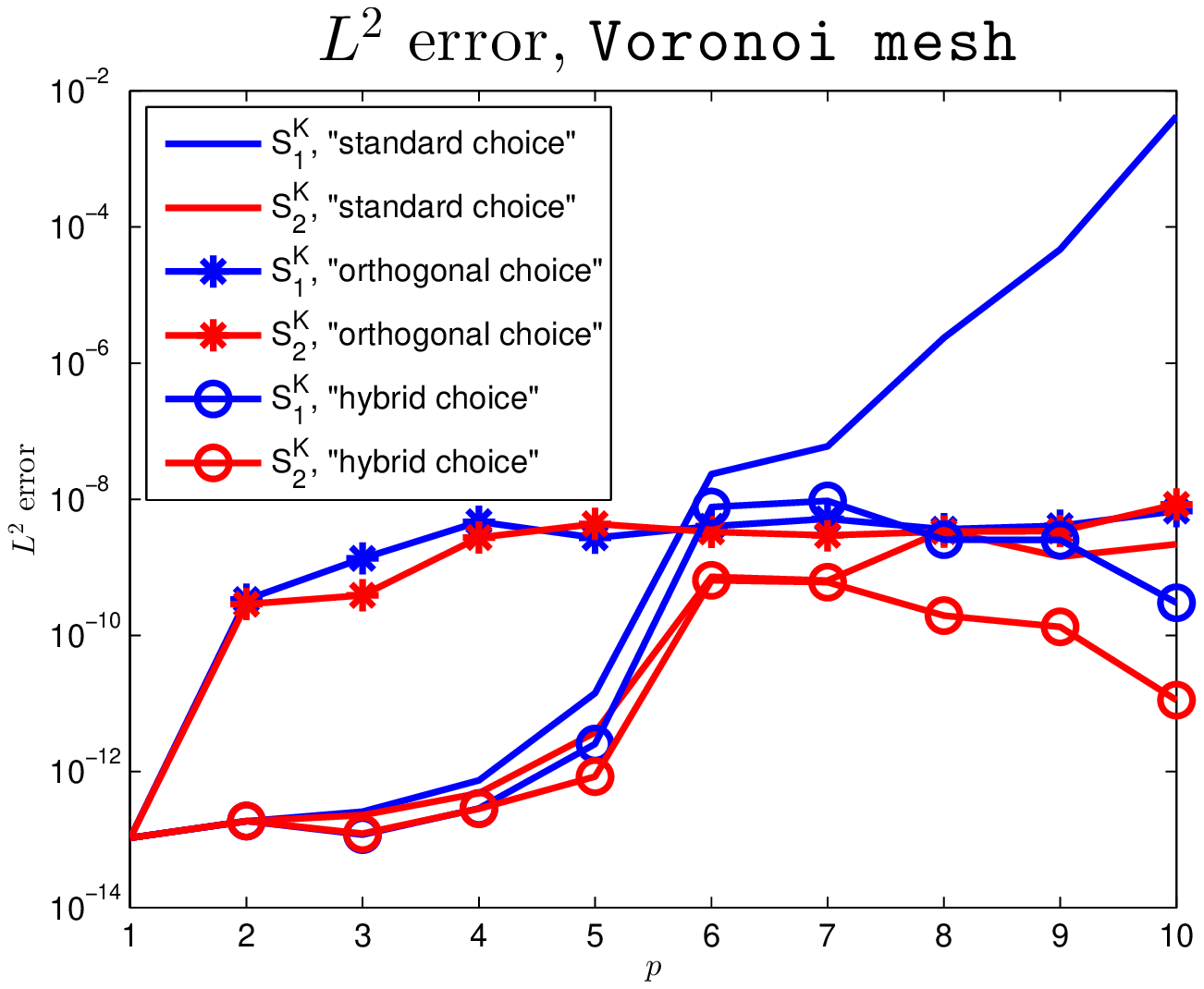}}
\caption{$\p$ version of VEM. Exact solution $\u_2$. Stabilizations employed: $\SE_1$ and $\SE_2$.
Face/bulk moments employed: ``standard choice'', ``orthogonal choice'' and ``hybrid choice''.
\texttt{Voronoi mesh}. Left: $H^1$ error. Right: $L^2$ error.}
\label{figure changing polynomial bases patch test Voronoi}
\end{figure}

\begin{figure}  [h]
\centering
\subfigure {\includegraphics [angle=0, width=0.49\textwidth]{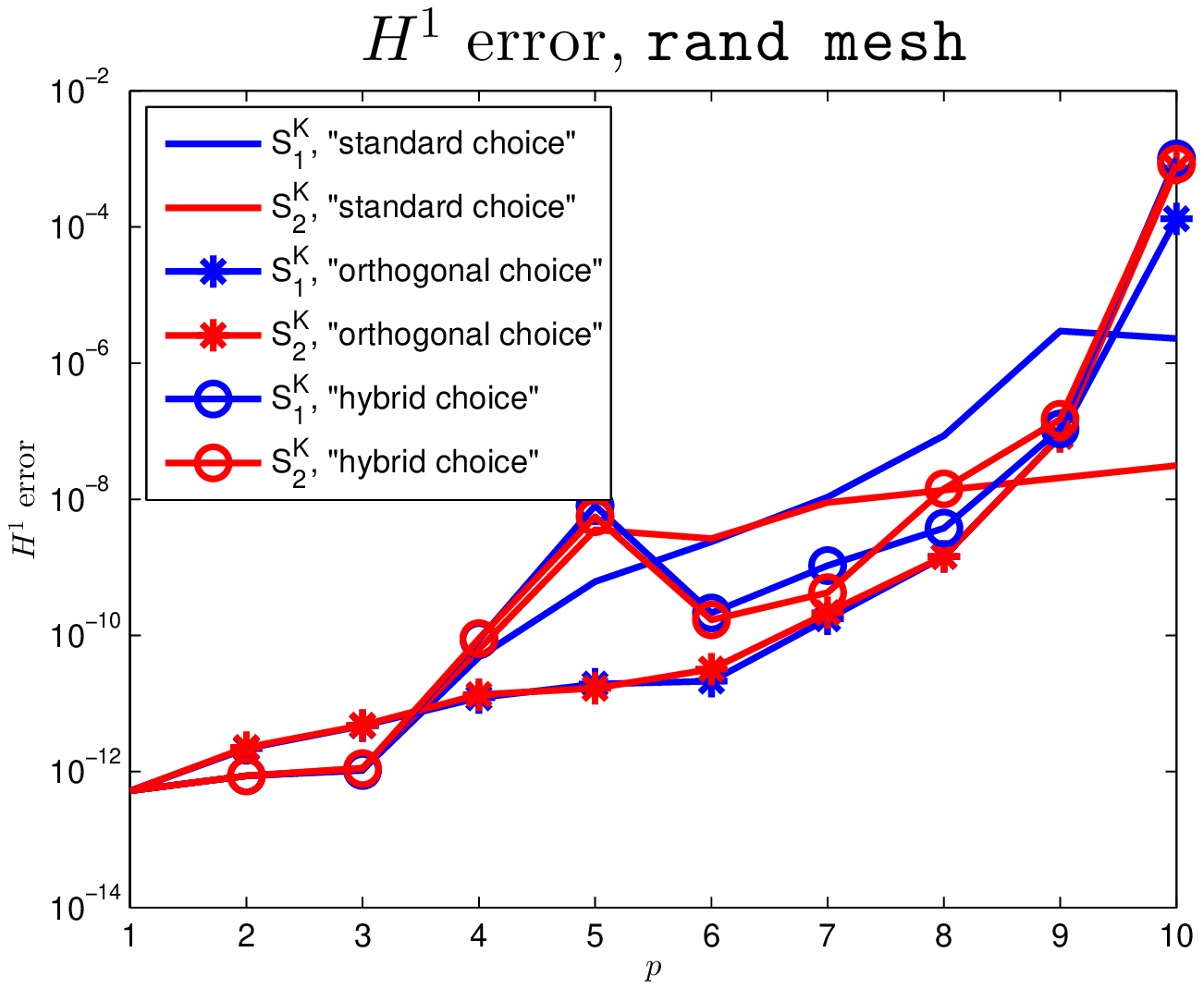}}
\subfigure {\includegraphics [angle=0, width=0.49\textwidth]{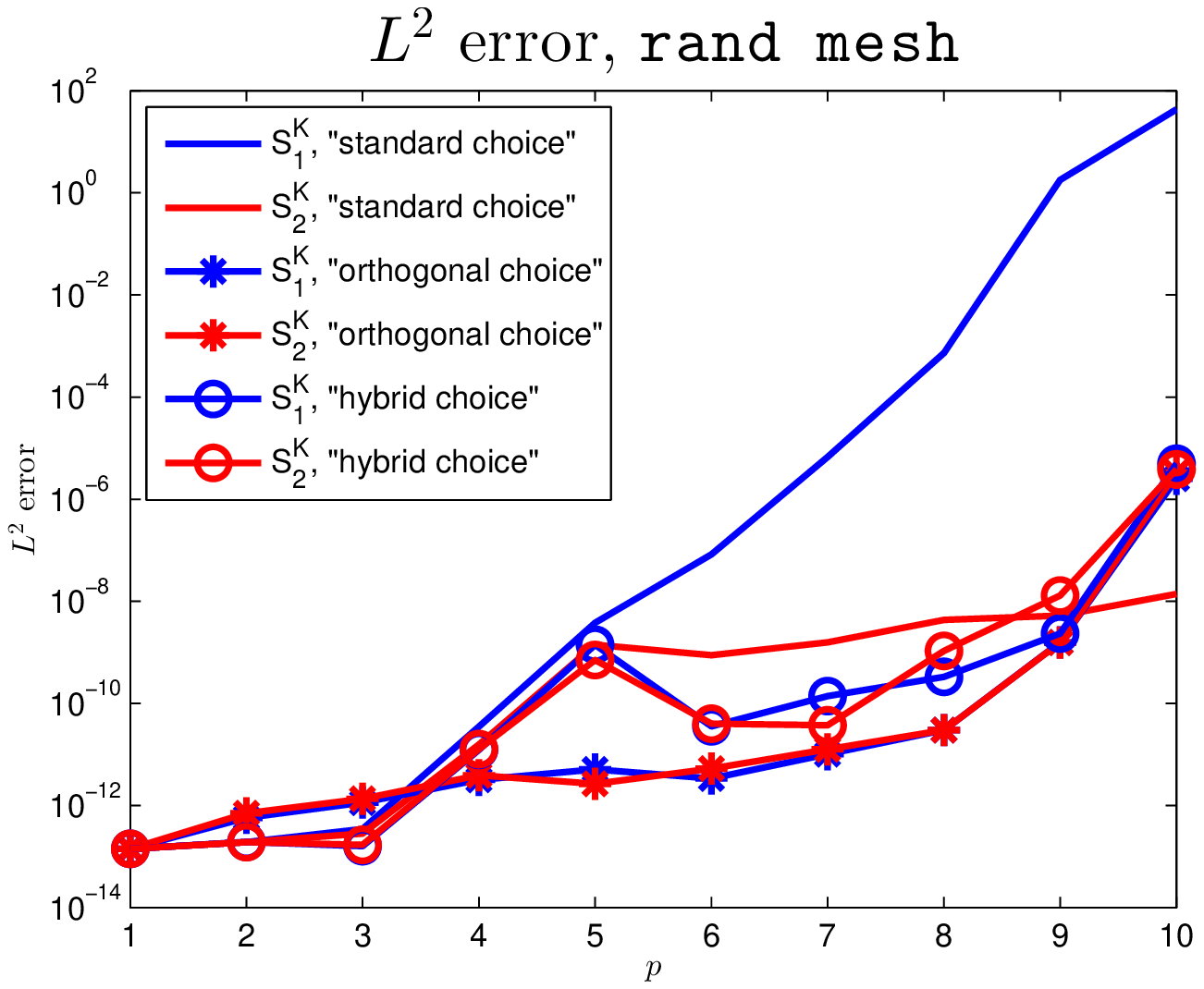}}
\caption{$\p$ version of VEM. Exact solution $\u_2$. Stabilizations employed: $\SE_1$ and $\SE_2$.
Face/bulk moments employed: ``standard choice'', ``orthogonal choice'' and ``hybrid choice''.
\texttt{Rand mesh}. Left: $H^1$ error. Right: $L^2$ error.}
\label{figure changing polynomial bases patch test rand}
\end{figure}

As already observed, in order to see if the condition number is the reason for the suboptimal convergence, we have to test the method on a polynomial solution;
for this reason, we consider in Figures \ref{figure changing polynomial bases patch test cube}, \ref{figure changing polynomial bases patch test Voronoi} and \ref{figure changing polynomial bases patch test rand} the same set of numerical tests
exhibited in Figures \ref{figure changing polynomial bases analytic solution cube}, \ref{figure changing polynomial bases analytic solution Voronoi} and \ref{figure changing polynomial bases analytic solution rand}
applied now to the patch test solution $\u_2$.
%In fact, we know that VEM returns up to machine precision, the exact solution whenever it is a polynomial of degree $\p$;
%therefore, the error that the method returns is due to the effect of the condition number of the stiffness matrix only.

%In Figures \ref{figure changing polynomial bases patch test cube}, \ref{figure changing polynomial bases patch test Voronoi} and \ref{figure changing polynomial bases patch test rand}, we present analogous results
%to those of Figure \ref{figure changing polynomial bases analytic solution cube}, \ref{figure changing polynomial bases analytic solution Voronoi} and \ref{figure changing polynomial bases analytic solution rand}
%employing now, instead of exact solution $\u_1$ defined in \eqref{solution 1}, solution $\u_2$ defined in \eqref{solution 2}.

We deduce that the problem of unexpected decay of the error when employing the ``standard choice'' is not \emph{only} ill-conditioning. Let us focus for instance on the \texttt{Voronoi mesh}.
Employing the ``standard choice'' on the patch test, the $L^2$ error with $\p=10$ is around $10^{-8}$ when using stabilization $\SE_2$ and around $10^{-2}$ when using stabilization $\SE_1$,
while, when testing analytic solution $\u_1$, the $L^2$ error with both stabilization is $10^7$!

Such mismatch between results approximating solutions $\u_1$ and $\u_2$ are not observed with the ``orthogonal'' and ``hybrid choice''.
As a consequence, the two novel choices, i.e. the ``orthogonal'' and the ``hybrid choice'', make the method even more robust when defining the stabilization.

\paragraph*{Summary:}
in $\p$ VEM, employ the ``orthogonal'' and the ``hybrid choices'' (robustness in the choice of stabilizations!); avoid the ``standard choice''.
\par

%%%%%%%%%%%%%%%%%%%%%%%%%%%%%%%%%%%%%%%%%%%%%%%%%

\subsection{Numerical results: collapsing elements} \label{subsection nr: collapsing}
Having observed in Section \ref{subsection nr: p version} that the ``orthogonal'' and the ``hybrid choices'' entail more robustness with respect to the choice of the stabilization than the ``standard choice'',
at least in the $\p$ version of the method,
we want to show here that such two choices perform better also when considering meshes with ``collapsing polyhedra''.

\begin{figure}[!htb]
\begin{center}
\includegraphics[width=0.7\textwidth]{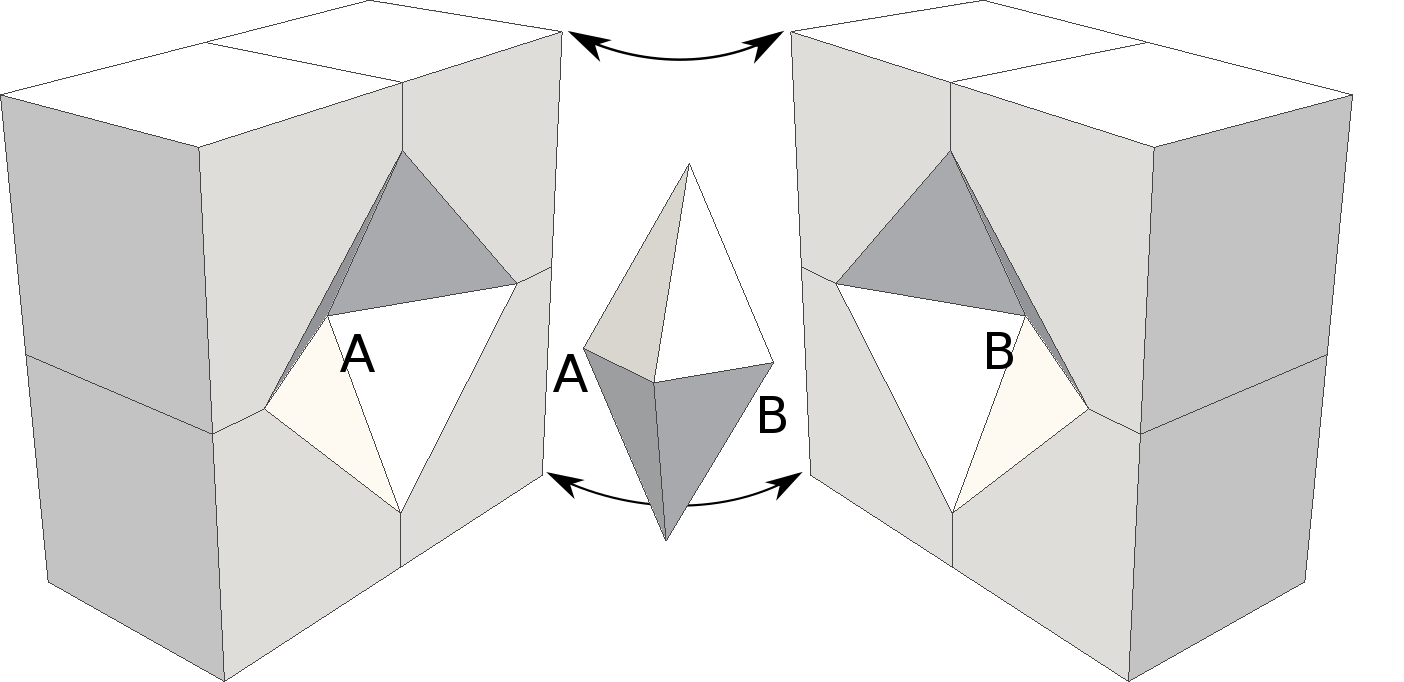} 
\end{center}
\caption{The mesh considered for the ``collapsing test''. 
%We show an exploded representation of this mesh to better understand how it is built and 
We highlight the points $A$ and $B$ that we shift towards the center of mass of the cube, in order to build the sequence of meshes under considerations.}
\label{fig:meshColl}
\end{figure}

For the purpose, we consider a sequence of meshes as that depicted in Figure~\ref{fig:meshColl}.
Such meshes are obtained by splitting the unit cube $[0,\,1]^3$ into four ``boundary'' polyhedra and an internal octahedron.
The sequence is built by shifting two vertices of the octahedron towards $\mathbf x_\Omega = (0.5,\,0.5,\,0.5)$, the center of mass of the cube,
for example we shift the points $A$ and $B$ represented in Figure~\ref{fig:meshColl} towards $\mathbf x_\Omega$.
In this way, the volume of the octahedron collapses.

More precisely, we consider a sequence of five meshes. The first one is built by taking points $A$ and $B$ having the following coordinates:
\[
A=(0.25,\, 0.5,\, 0.5), \quad \quad B=(0.75,\, 0.5,\, 0.5).
\]
The other meshes are obtained by moving these two points towards $\mathbf x_\Omega$ by halving at each step their distance.

\begin{figure}  [h]
\centering
\subfigure {\includegraphics [angle=0, width=0.49\textwidth]{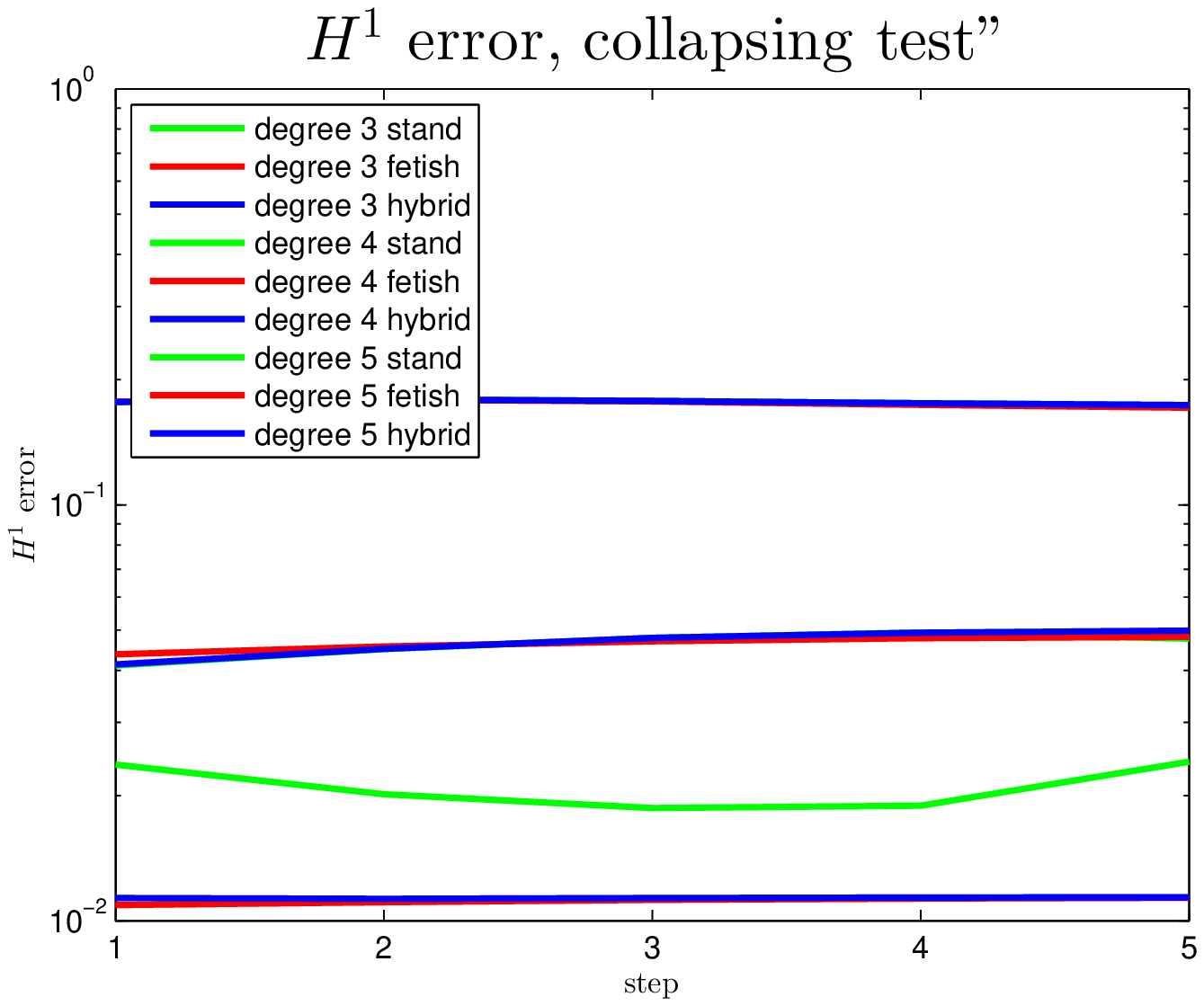}}
\subfigure {\includegraphics [angle=0, width=0.49\textwidth]{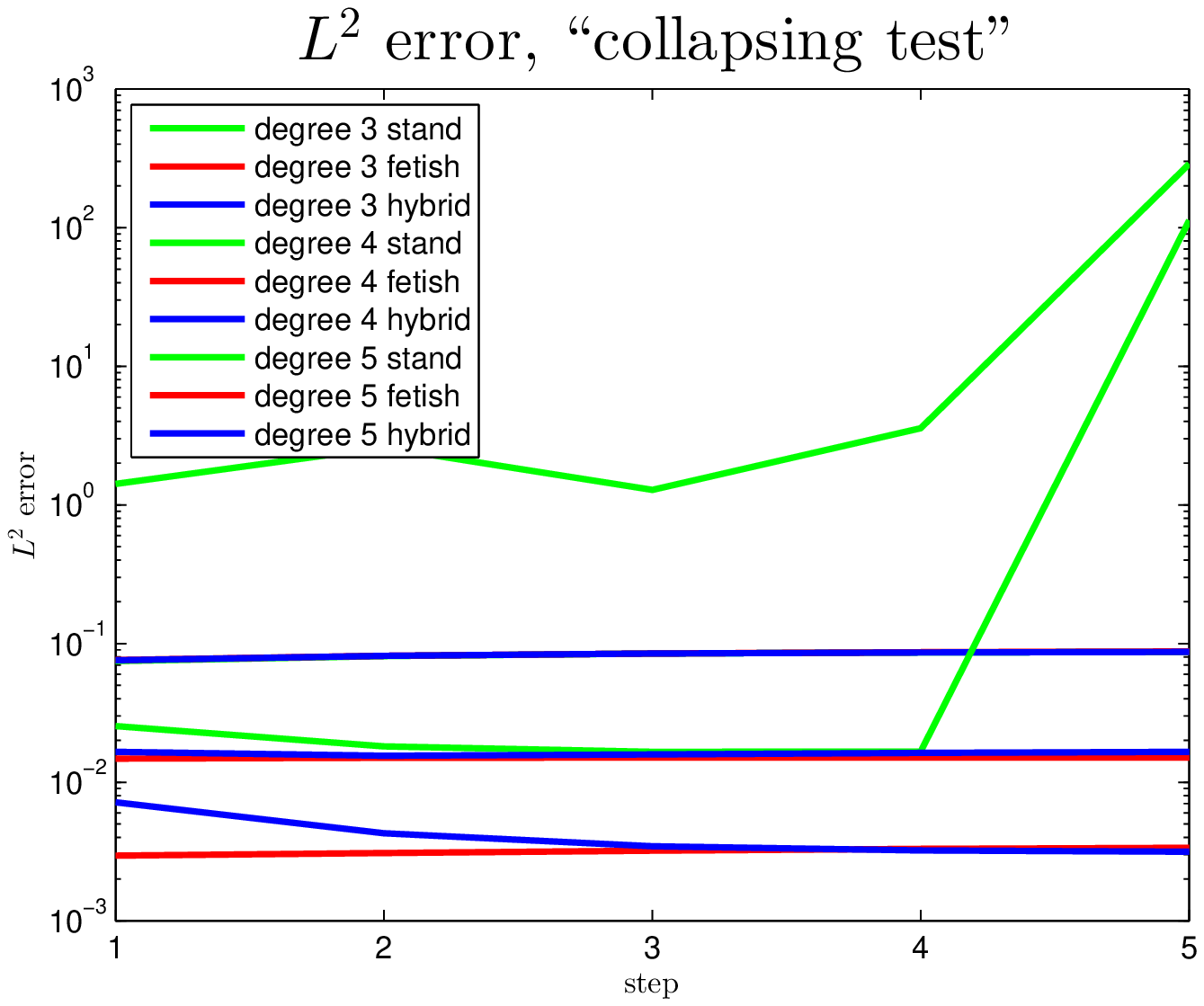}}
\caption{VEM applied to ``collapsing sequences'' of meshes, see Figure \eqref{fig:meshColl}.
Exact solution $\u_1$. Stabilization employed: $\SE_1$.
Face/bulk moments employed: ``standard choice'', ``orthogonal choice'' and ``hybrid choice''.
Degrees of accuracy: $\p=3$, $4$, $5$. Left: $H^1$ error. Right: $L^2$ error.}
\label{figure collapsing}
\end{figure}

We numerically investigate in Figure \ref{figure collapsing} what happens to the $H^1$ and the $L^2$ errors defined in \eqref{computed errors}, when considering as an exact solution $\u_1$ defined in \eqref{solution 1}
and employing stabilization $\SE_2$ defined in \eqref{D recipe stabilization}.
We consider various degrees of accuracy, namely $\p=3,4,5$.
What we observe here is that again the ``standard choice'' suffers after some ``collapsing iterations''.

Importantly, we carried out numerical experiments with stabilization $\SE_1$ and even worse performances of the ``standard choice'' have been observed.

\paragraph*{Summary:}
in VEM on meshes with ``collapsing polyhedra'', employ the ``orthogonal'' and the ``hybrid choices'';
avoid the ``standard choice''.
\par

%%%%%%%%%%%%%%%%%%%%%%%%%%%%%%%%%%%%%%%%%%%%%%%%%
\subsection{Numerical results: the $\h$ version of 3D VEM} \label{subsection nr: h VEM}
In the foregoing sections, we observed that the ``orthogonal'' and the `hybrid choice'' produce similar results when dealing with the $\p$ version of VEM
as well as when considering meshes characterized by elements with collapsing bulk.

In this section, we consider instead the $\h$ version of VEM and we compare the effects on the convergence of the error employing again the ``standard'', the ``orthogonal'' and the ``hybrid choices''.

We aim to approximate the analytic solution $\u_1$.
We consider two sequences of meshes, namely a \texttt{cube mesh} and \texttt{Voronoi mesh}, and
we fix $\SE_2$ as a stabilization.

\begin{figure}  [h]
\centering
\subfigure {\includegraphics [angle=0, width=0.49\textwidth]{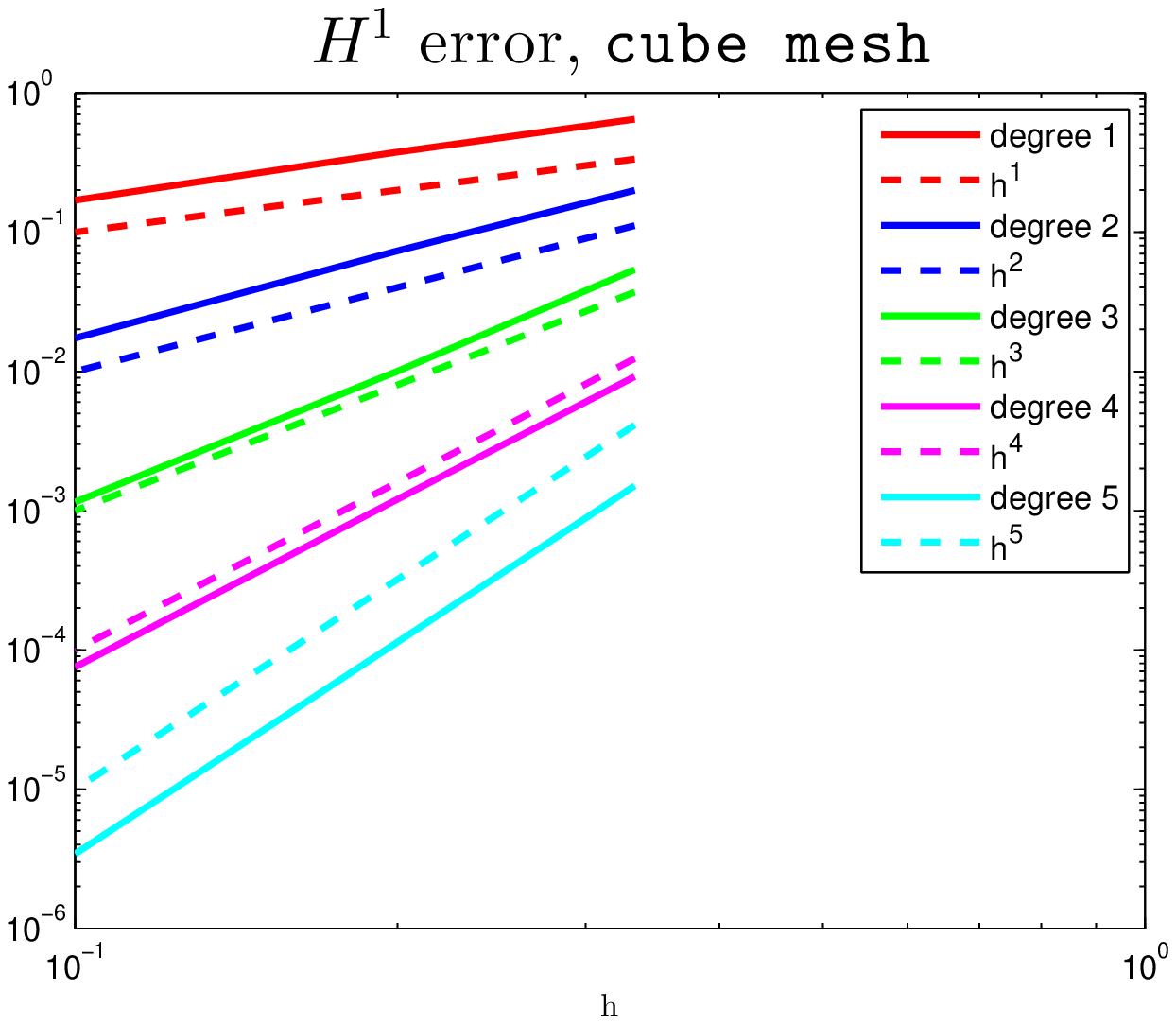}}
\subfigure {\includegraphics [angle=0, width=0.49\textwidth]{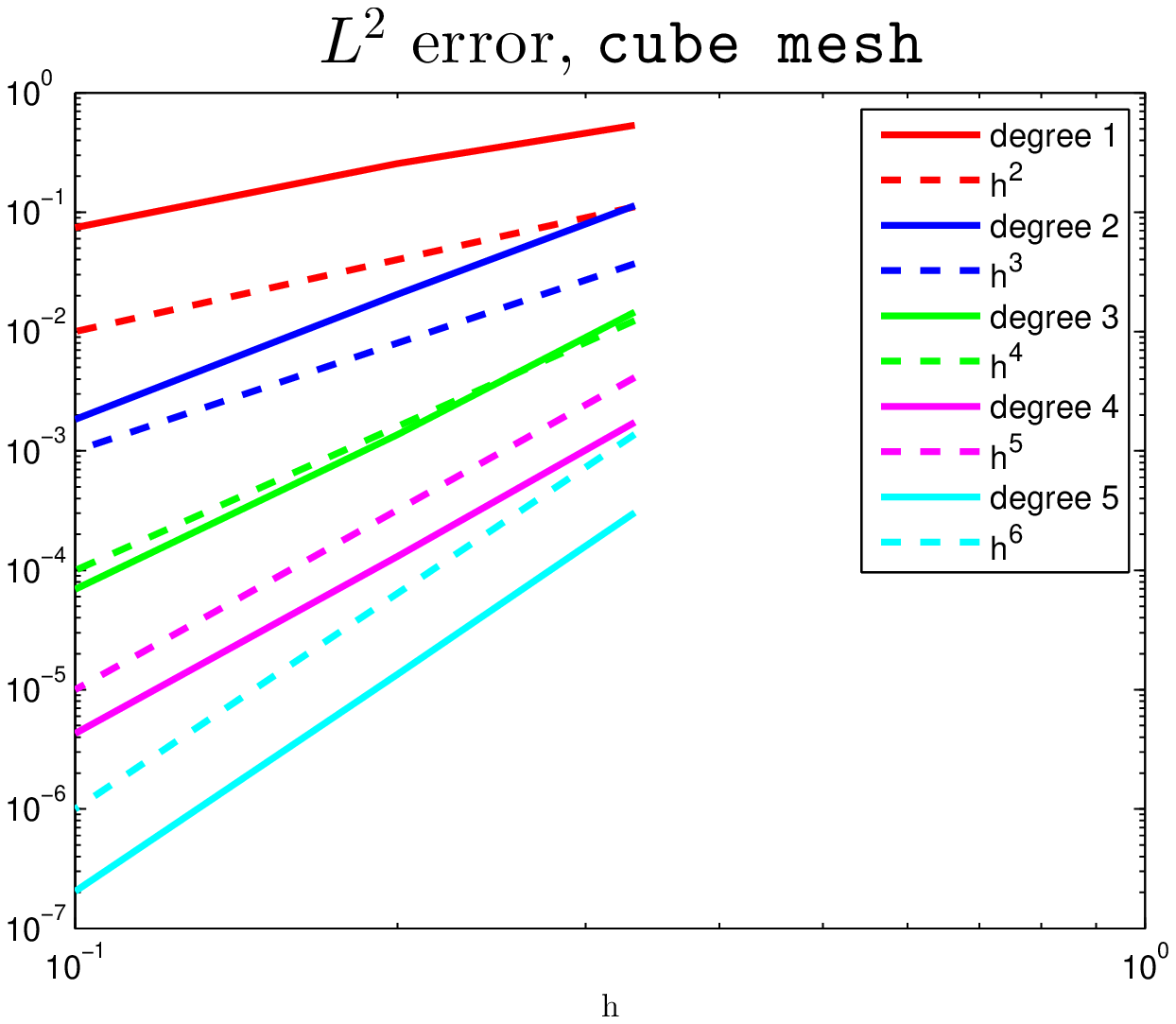}}
\caption{$\h$ version of VEM. Exact solution $\u_1$. Stabilization employed: $\SE_2$. Face/bulk moments employed: ``standard choice''.
\texttt{Cube mesh}. Left: $H^1$ error. Right: $L^2$ error.}
\label{figure h VEM cube analytic solution_stand}
\end{figure}

\begin{figure}  [h]
\centering
\subfigure {\includegraphics [angle=0, width=0.49\textwidth]{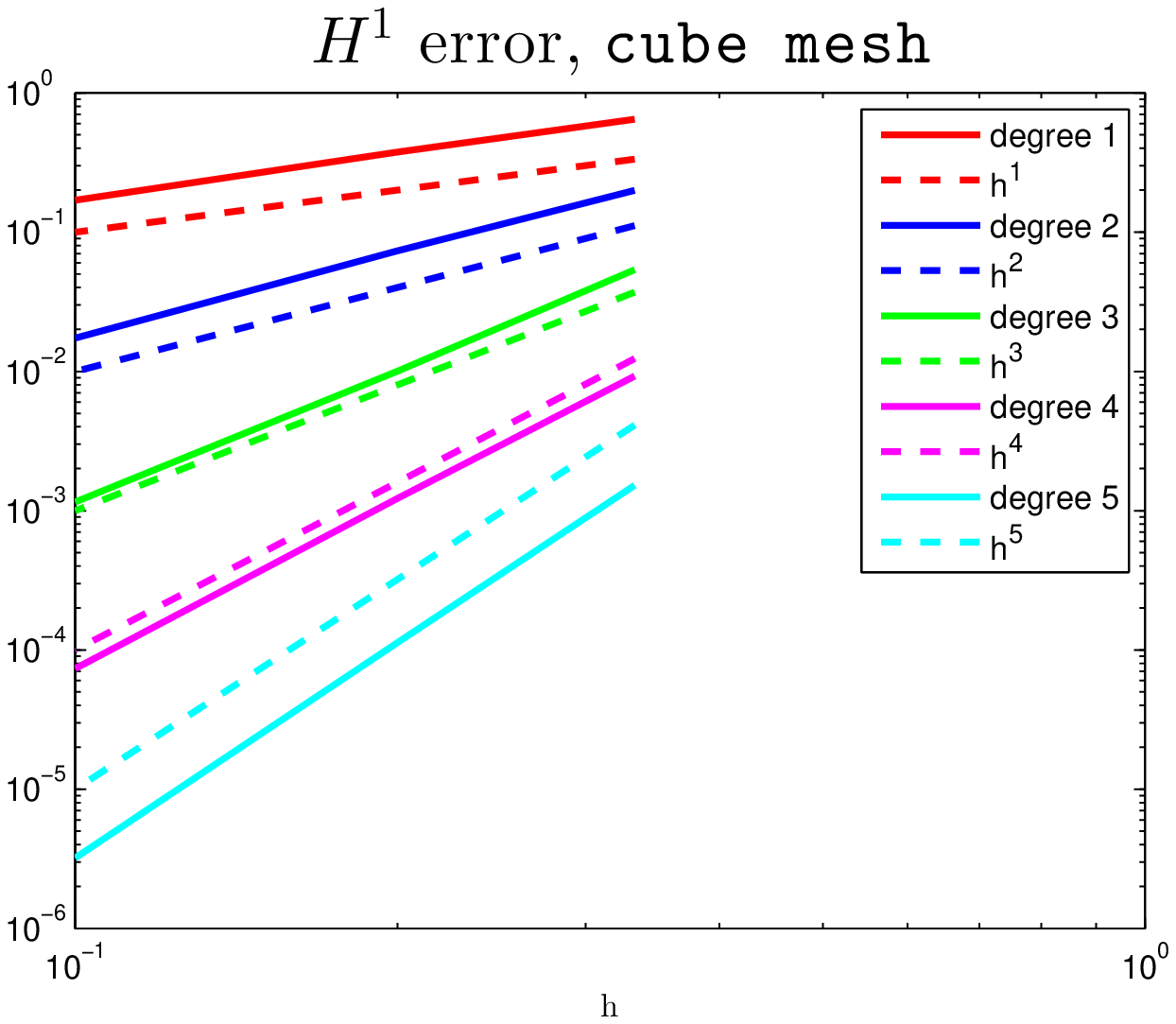}}
\subfigure {\includegraphics [angle=0, width=0.49\textwidth]{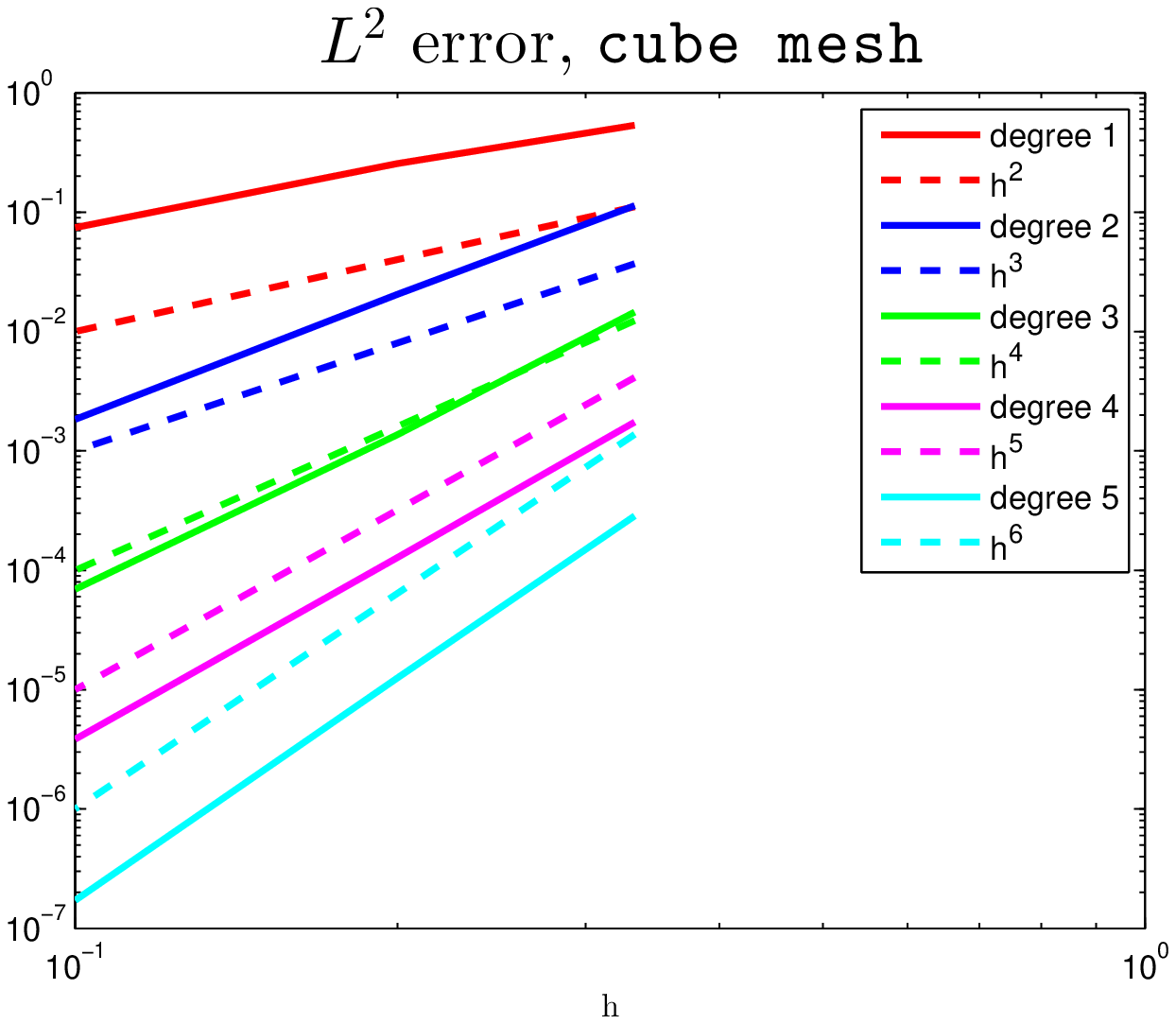}}
\caption{$\h$ version of VEM. Exact solution $\u_1$. Stabilization employed: $\SE_2$. Face/bulk moments employed: ``orthogonal choice''.
\texttt{Cube mesh}. Left: $H^1$ error. Right: $L^2$ error.}
\label{figure h VEM cube analytic solution_fetish}
\end{figure}

\begin{figure}  [h]
\centering
\subfigure {\includegraphics [angle=0, width=0.49\textwidth]{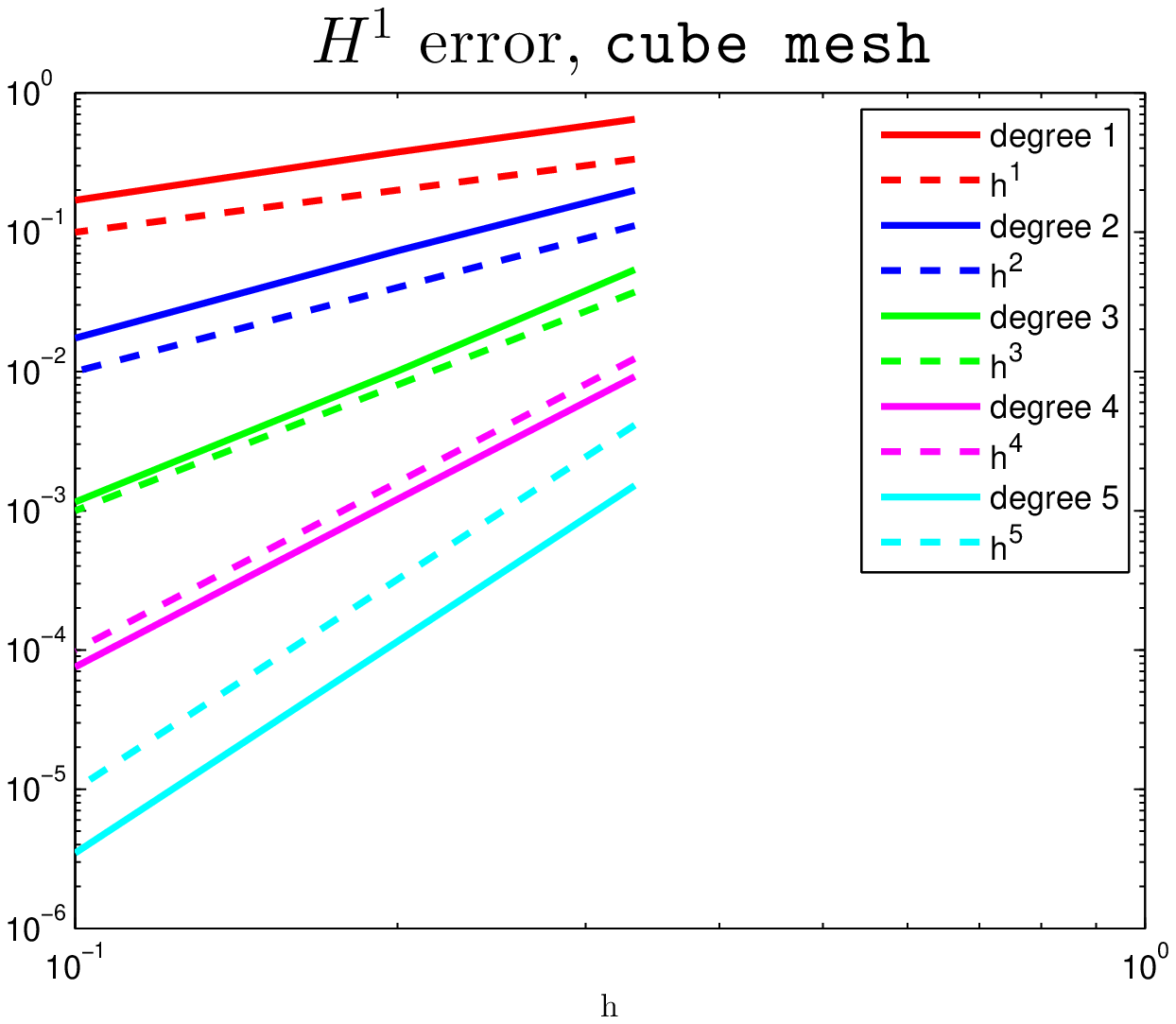}}
\subfigure {\includegraphics [angle=0, width=0.49\textwidth]{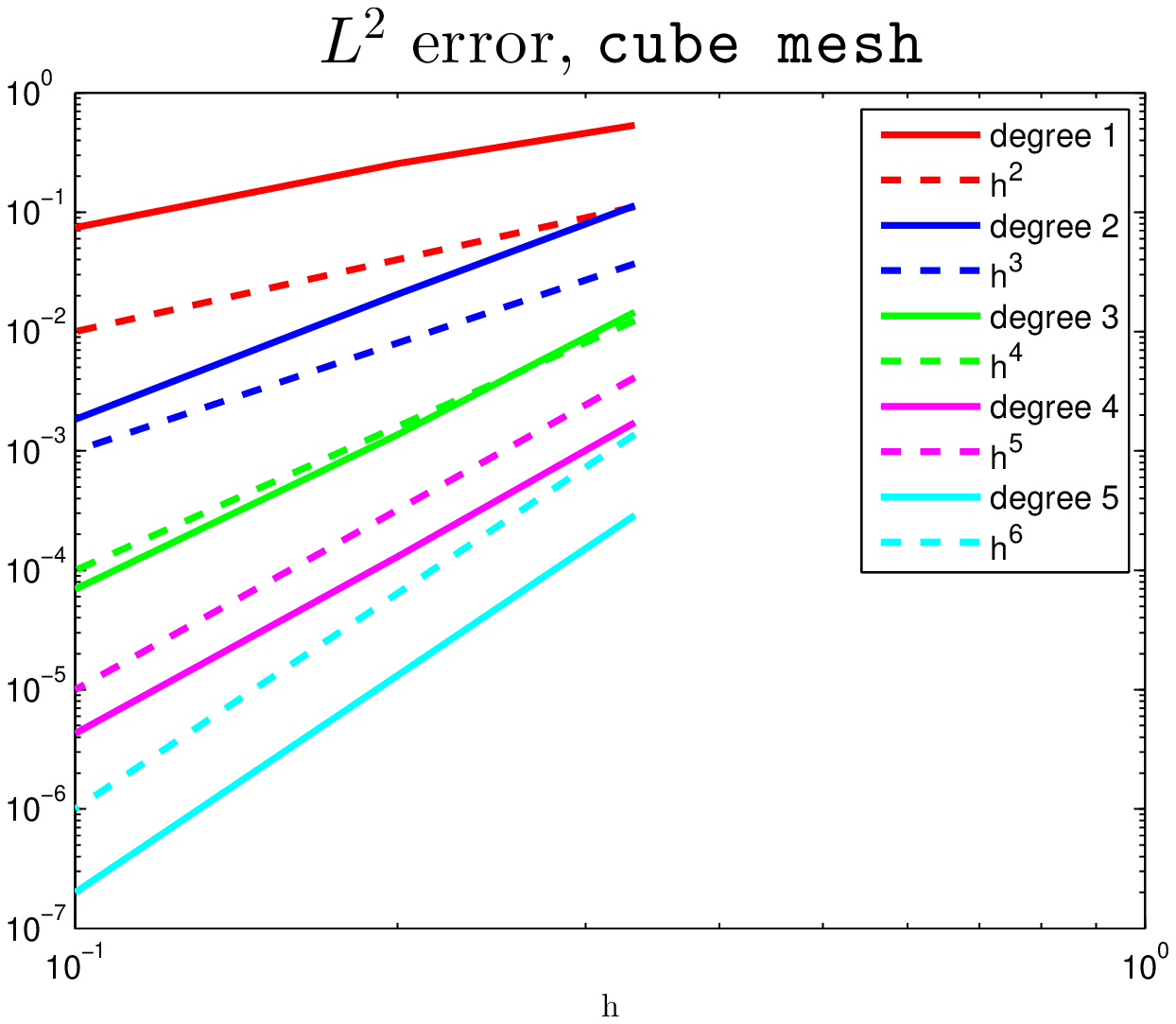}}
\caption{$\h$ version of VEM. Exact solution $\u_1$. Stabilization employed: $\SE_2$. Face/bulk moments employed: ``hybrid choice''.
\texttt{Cube mesh}. Left: $H^1$ error. Right: $L^2$ error.}
\label{figure h VEM cube analytic solution_hybrid}
\end{figure}

\begin{figure}  [h]
\centering
\subfigure {\includegraphics [angle=0, width=0.49\textwidth]{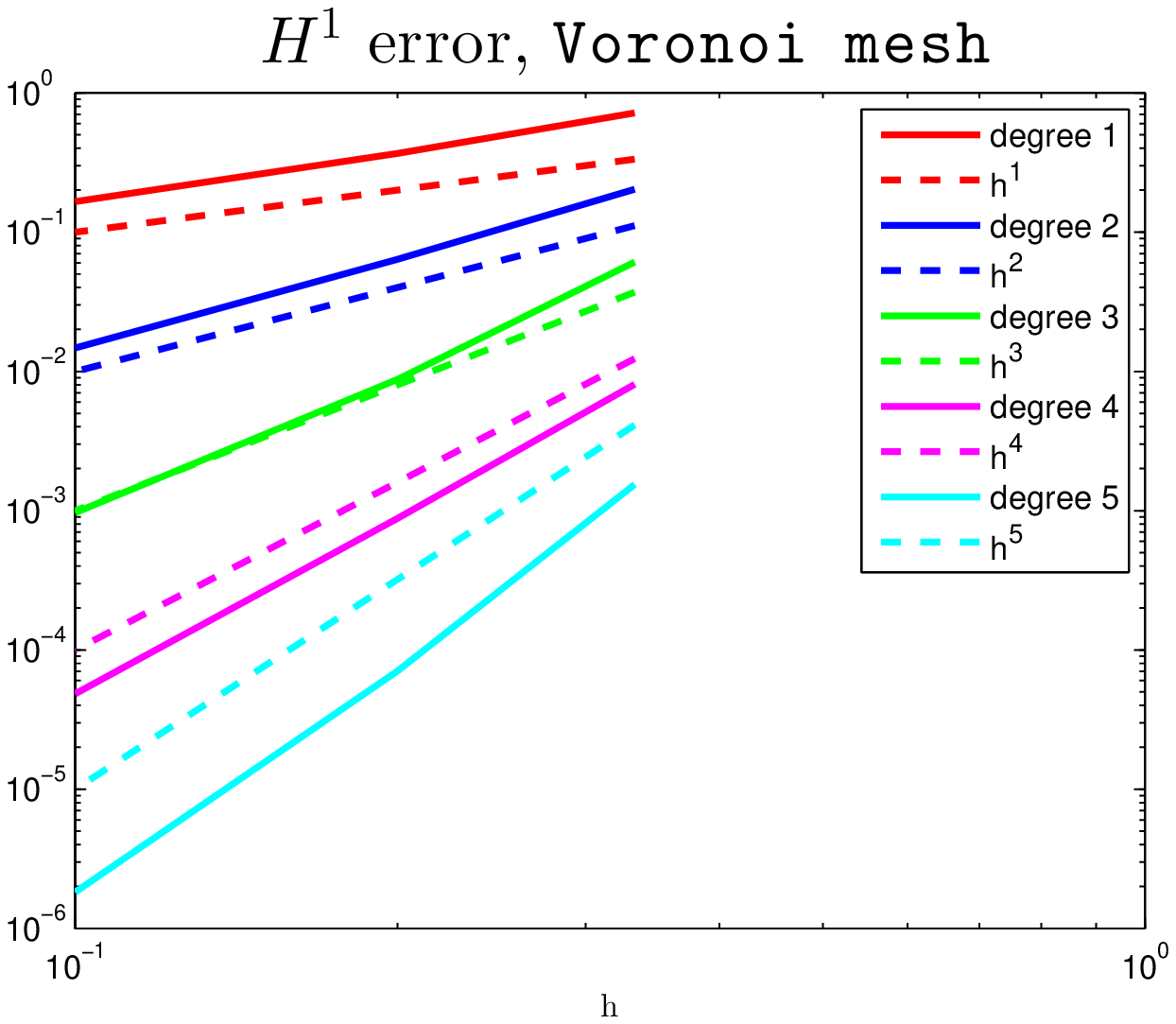}}
\subfigure {\includegraphics [angle=0, width=0.49\textwidth]{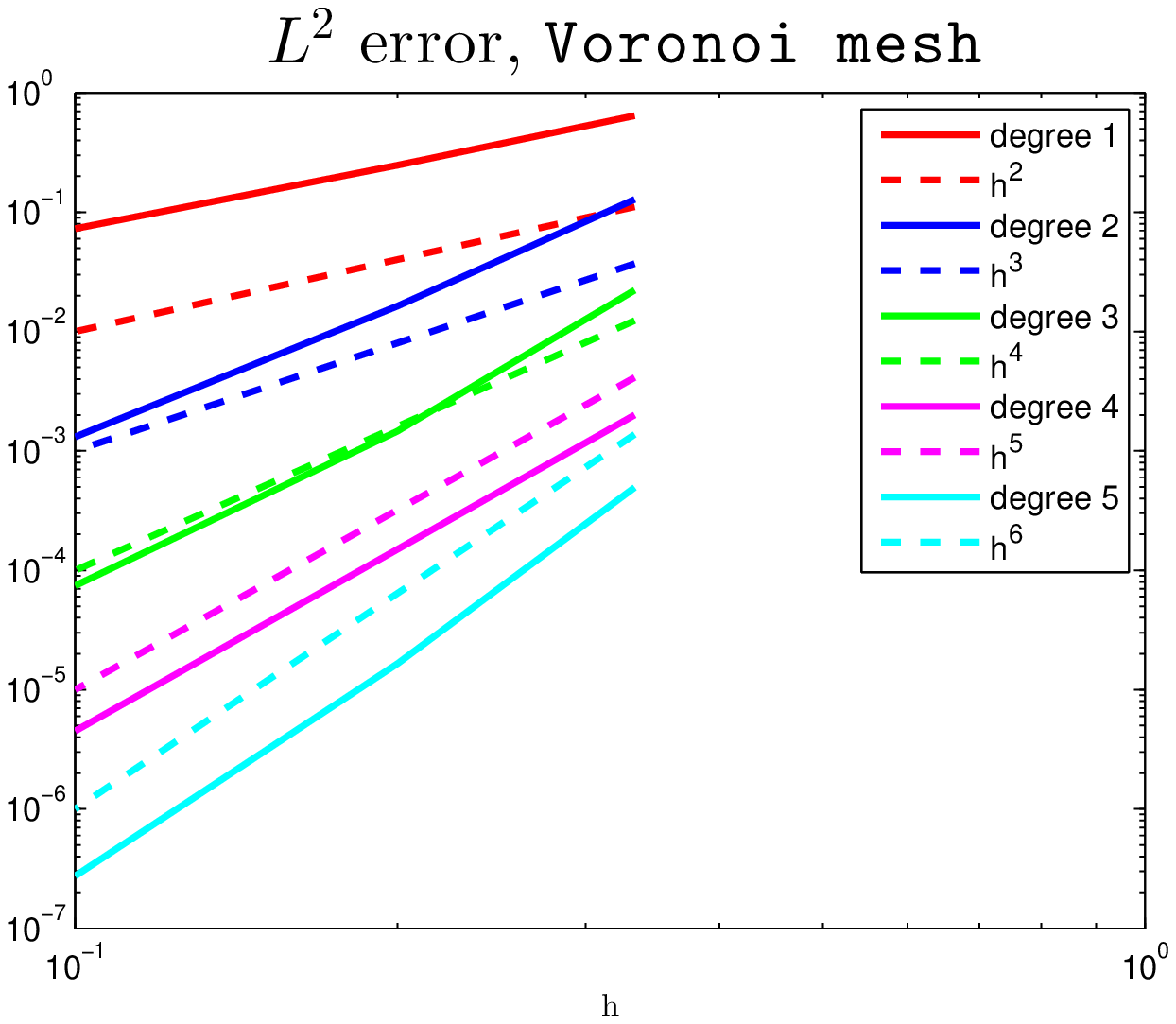}}
\caption{$\h$ version of VEM. Exact solution $\u_1$. Stabilization employed: $\SE_2$.
Face/bulk moments employed: ``standard choice''.
\texttt{Voronoi mesh}. Left: $H^1$ error. Up-right: $L^2$ error.}
\label{figure h VEM Voronoi analytic solution_stand}
\end{figure}

\begin{figure}  [h]
\centering
\subfigure {\includegraphics [angle=0, width=0.49\textwidth]{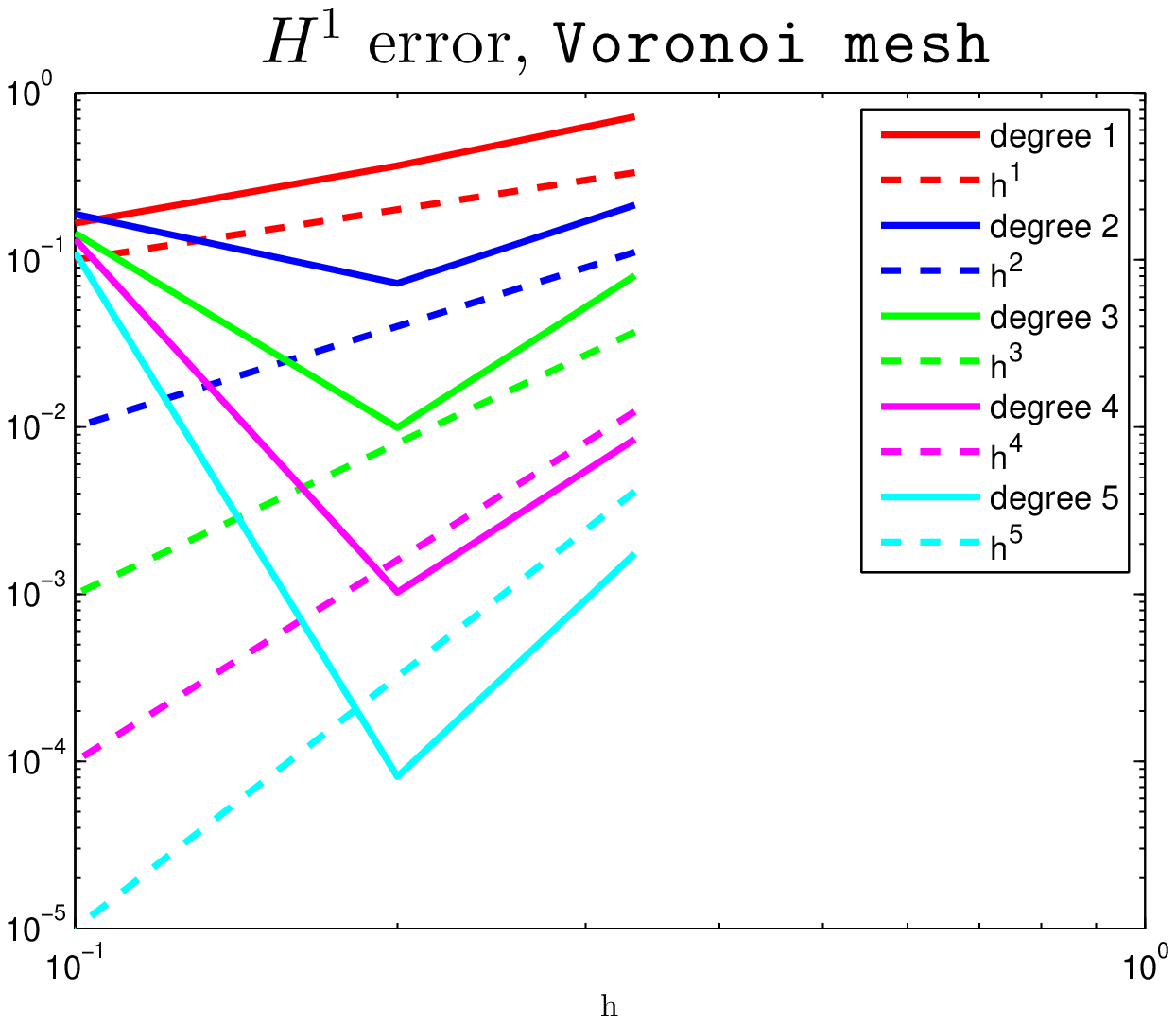}}
\subfigure {\includegraphics [angle=0, width=0.49\textwidth]{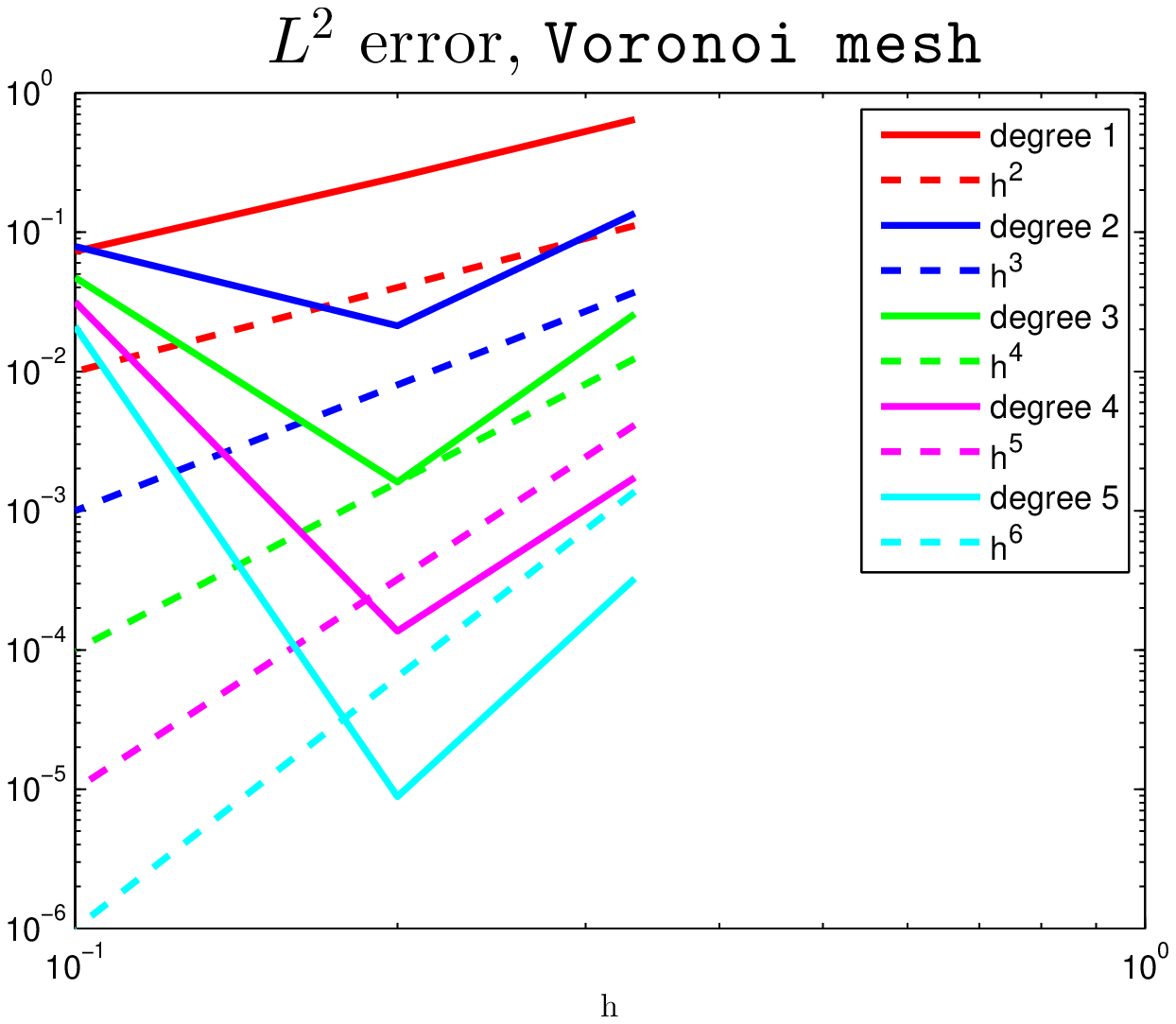}}
\caption{$\h$ version of VEM. Exact solution $\u_1$. Stabilization employed: $\SE_2$.
Face/bulk moments employed: ``orthogonal choice''.
\texttt{Voronoi mesh}. Left: $H^1$ error. Up-right: $L^2$ error.}
\label{figure h VEM Voronoi analytic solution_fetish}
\end{figure}

\begin{figure}  [h]
\centering
\subfigure {\includegraphics [angle=0, width=0.49\textwidth]{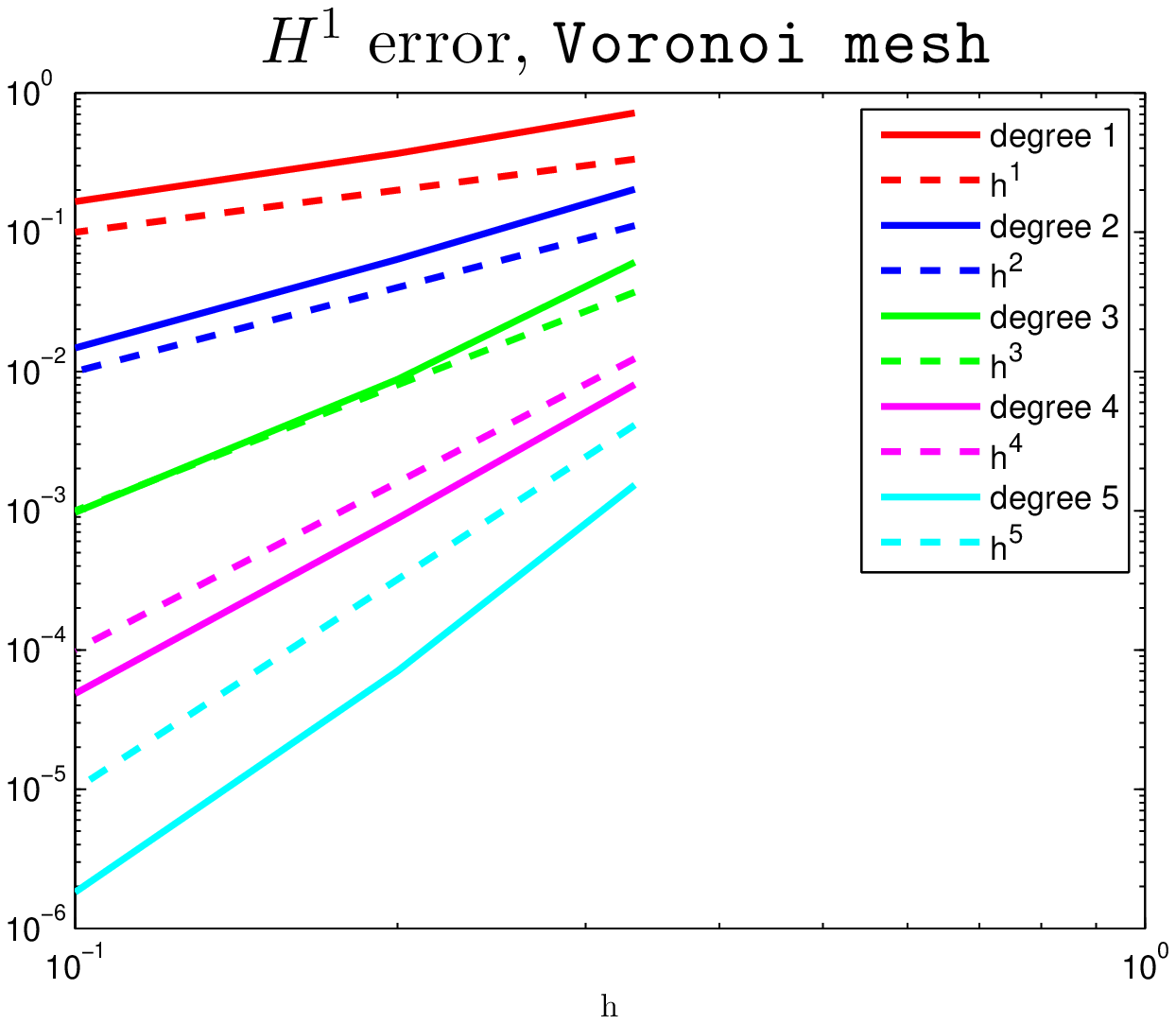}}
\subfigure {\includegraphics [angle=0, width=0.49\textwidth]{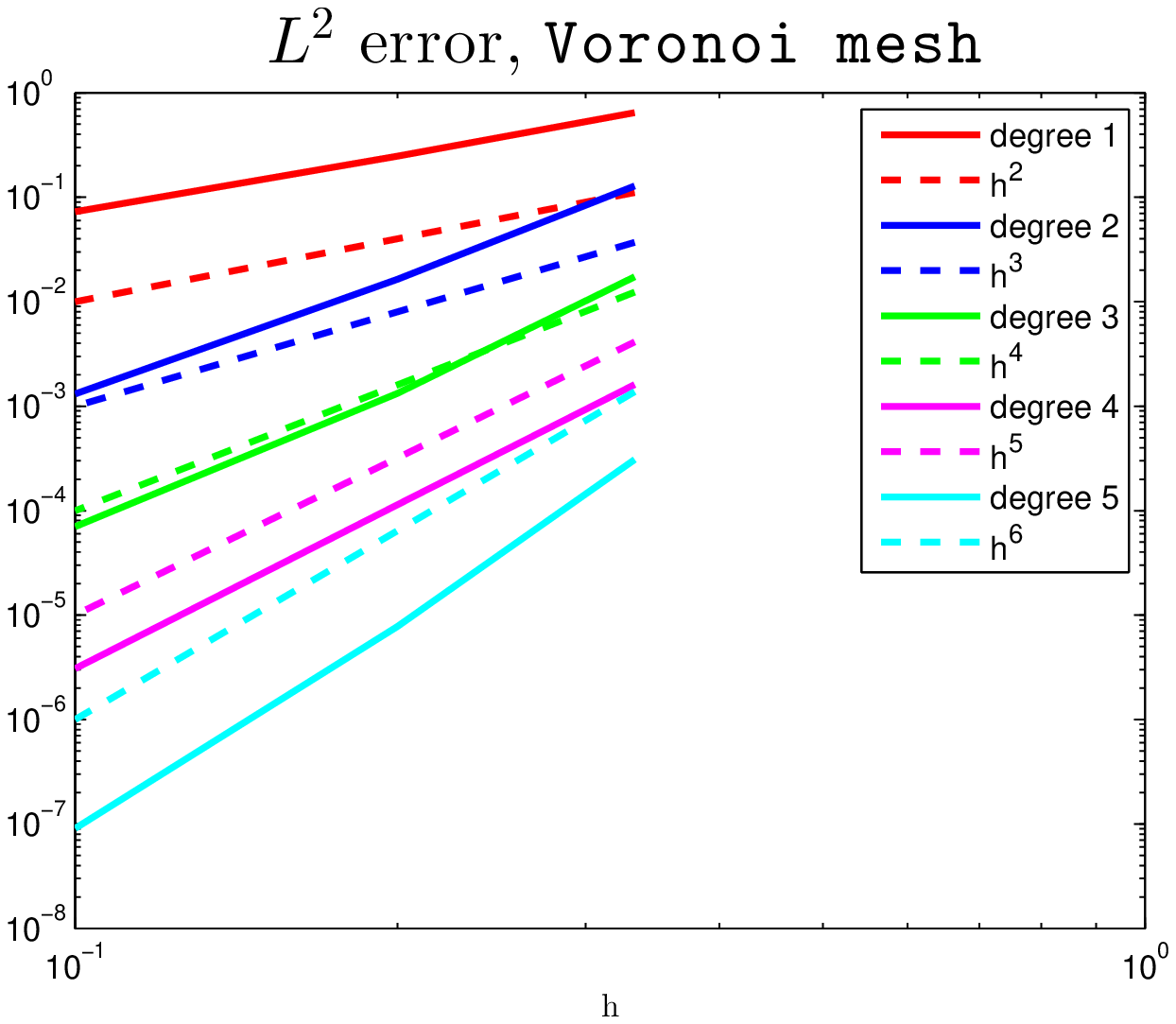}}
\caption{$\h$ version of VEM. Exact solution $\u_1$. Stabilization employed: $\SE_2$.
Face/bulk moments employed: ``hybrid choice''.
\texttt{Voronoi mesh}. Left: $H^1$ error. Up-right: $L^2$ error.}
\label{figure h VEM Voronoi analytic solution_hybrid}
\end{figure}

In Figures \ref{figure h VEM cube analytic solution_stand}, \ref{figure h VEM cube analytic solution_fetish} and \ref{figure h VEM cube analytic solution_hybrid} we perform the tests on the sequence of \texttt{cube meshes} whereas,
in Figures \ref{figure h VEM Voronoi analytic solution_stand}, \ref{figure h VEM Voronoi analytic solution_fetish} and \ref{figure h VEM Voronoi analytic solution_hybrid} we perform the tests on the sequence of \texttt{Voronoi meshes}.


What we deduce is that on sequences of regular polygons all the bases perform rather well;
on the other hand, on sequences of less regular meshes, e.g. \texttt{Voronoi meshes}, the orthogonalization process on faces leads to very bad (even increasing sometimes) error slopes.

\begin{figure}  [h]
\centering
\subfigure {\includegraphics [angle=0, width=0.49\textwidth]{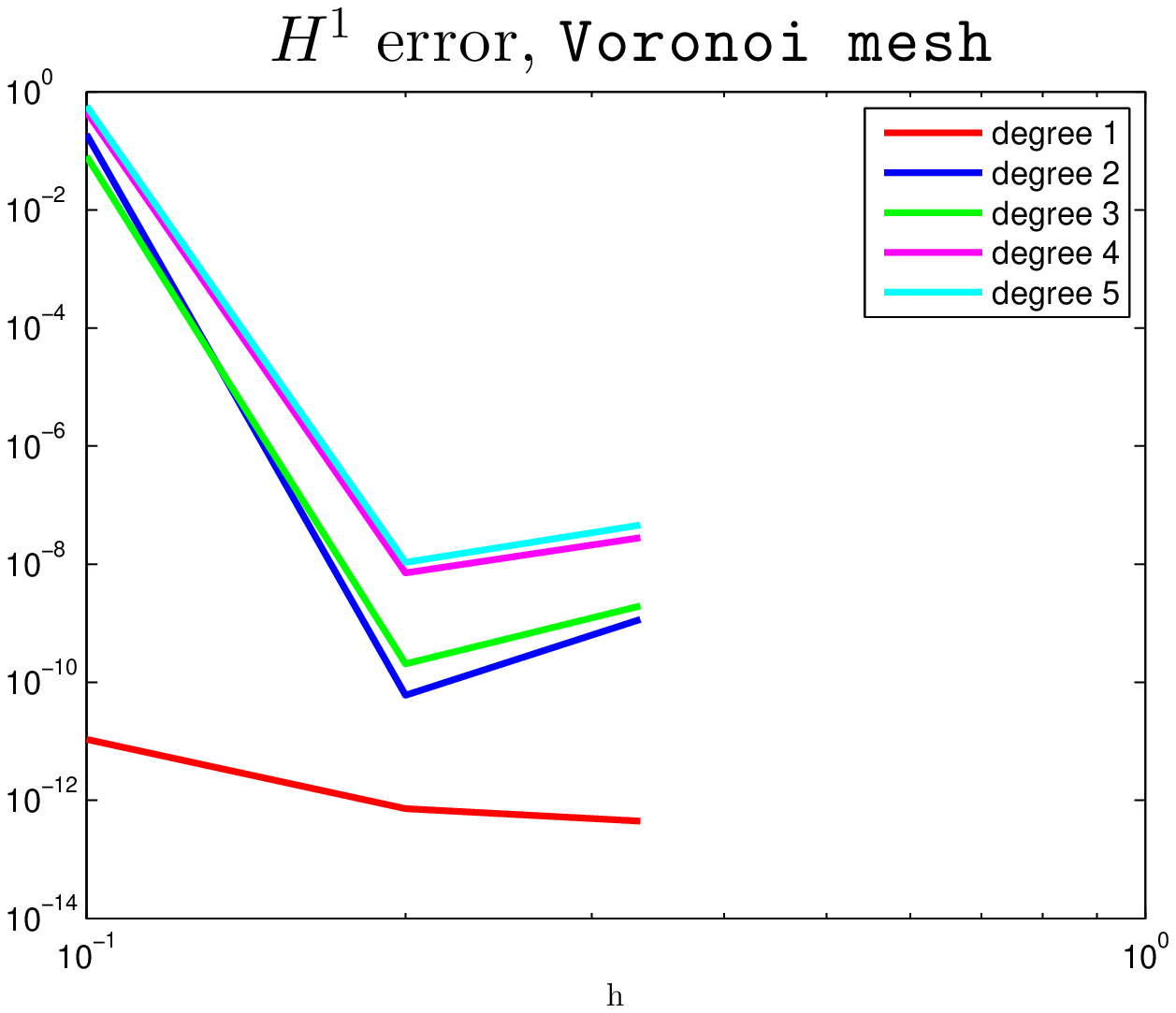}}
\subfigure {\includegraphics [angle=0, width=0.49\textwidth]{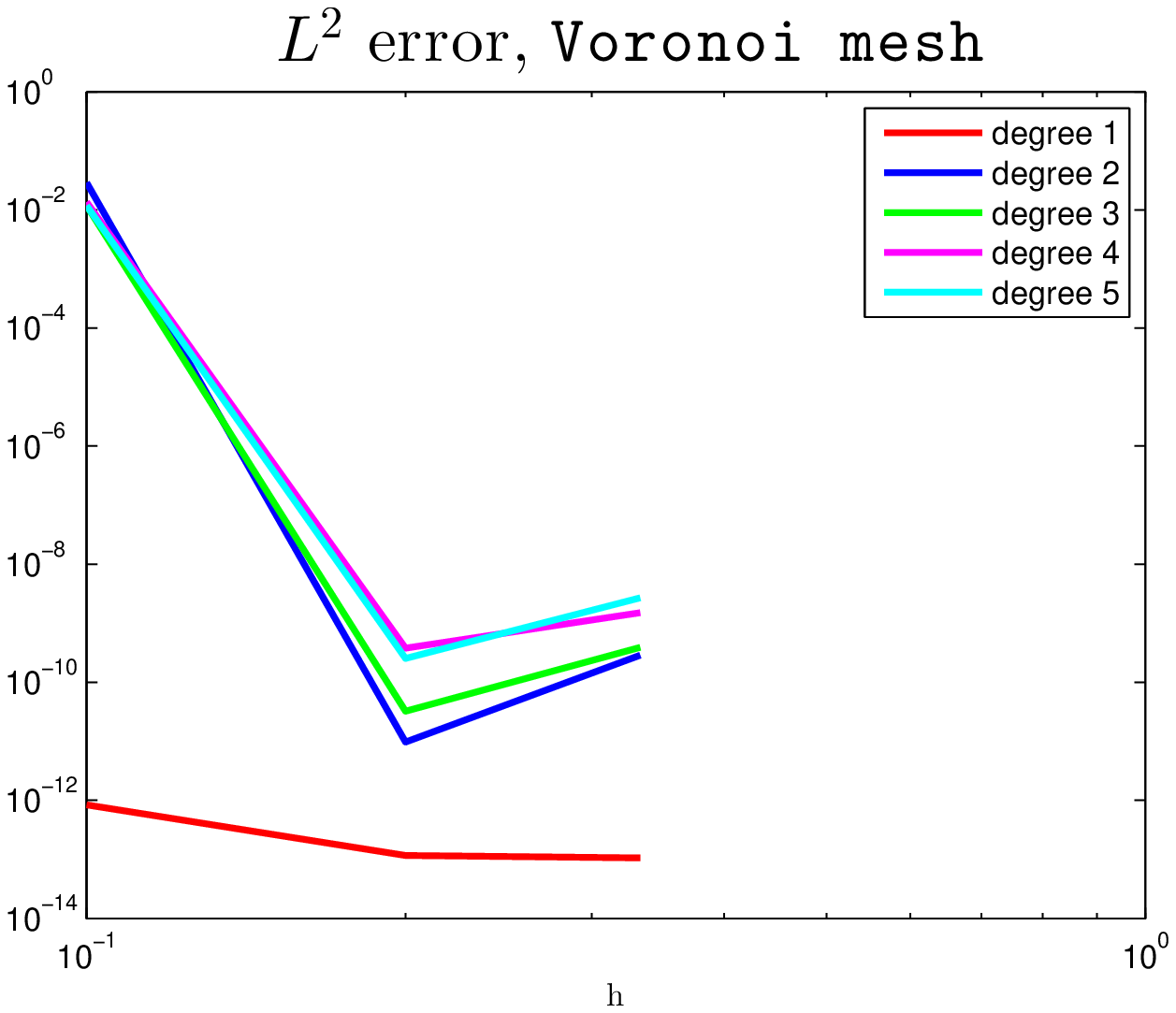}}
\caption{$\h$ version of VEM. Exact solution $\u_2$. Stabilizations employed: $\SE_2$.
Face/bulk moments employed: ``orthogonal choice''.
Left: $H^1$ error ``orthogonal choice''. Right: $L^2$ error ``orthogonal choice''.}
\label{figure h VEM Voronoi patch test}
\end{figure}

The reason for this bad behaviour is entirely ascribable to ill-conditioning. Indeed, in Figure~\ref{figure h VEM Voronoi patch test}, we apply as an example the VEM to the patch test $\u_2$,
taking for instance the \texttt{Voronoi mesh} and the ``orthogonal choice''.
What we observe is that the error grows precisely as it grows when approximating the analytic solution in Figure \ref{figure h VEM Voronoi analytic solution_fetish}.

\paragraph*{Summary:}
in $\h$ VEM, employ ``standard'' and ``hybrid choices''; avoid ``orthogonal choice''.
\par

% --------------------------------------------------------------------------------------------------------------------------------------------------------------------------------------------------------------------------------------------------------------------
\section{Conclusions} \label{section conclusions}
% --------------------------------------------------------------------------------------------------------------------------------------------------------------------------------------------------------------------------------------------------------------------
In this paper, we presented numerical tests dealing with high-order Virtual Element Method for the approximation of a three dimensional Poisson problem, extending thus the numerical analysis of \cite{preprint_VEM3Dbasic, fetishVEM}.

Moreover, we tested the method employing three different stabilizations and three different choices of face/bulk degrees of freedom.

It turned out that the stabilization leading to more performing results is $\SE_2$, which has the merit of ``leveling'' the contribution of the consistency and stability terms in the local splitting \eqref{local discrete bilinear form}.

Regarding the definition of face/bulk degrees of freedom, we have numerical evidence that the ``standard choice'', i.e. moments taken with respect to (scaled and centered) monomials,
leads to satisfactory decay of the errors in the $\h$ version of the method, but
implies suboptimal/wrong numerical results when employing both (moderately) high order degrees of accuracy and when employing meshes with ``bad-shaped'' elements;
this suboptimality can be alleviated by considering either the ``orthogonal choice'' or the ``hybrid choice''.%, i.e. by taking moments with respect to a basis of $L^2$ orthogonal polynomials;
%importantly, the deepest impact of the performances of the method is mainly due to bulk and not to face moments.

However, the ``orthogonal choice'' leads, on some sequences of meshes (characterized by small faces, small edges\dots), to inappropriate error convergence curves when employing the $\h$ version of the method,
whereas the ``hybrid choice'' seems to be more robust with this respect.

For this reason, if we have to suggest one recipe for building a VEM in three dimensions, we recommend the following:
\begin{itemize}
\item employ stabilization $\SE_2$;
\item employ the ``standard'' and the ``hybrid choices'' for the $\h$ version of the method, giving a preference to the ``hybrid choice'' in presence of `` collapsing'' elements;
\item employ the ``orthogonal'' and the ``hybrid choices'' for the $\p$ version of the method.
\end{itemize}

Importantly, we also discussed the possible reasons for which the method could return in some occurrences the wrong error convergence slopes.

The first one is the condition number of the stiffness matrix, whose effect can be determined by checking the errors on the patch test, e.g. on function $\u_2$.

The second one is the choice of the stabilization, which influences the condition number, but plays also a role in approximation estimates through the pollution factor $\alpha(\p)$
defined in \eqref{stabilization factor}.

\begin{appendix}
% --------------------------------------------------------------------------------------------------------------------------------------------------------------------------------------------------------------------------------------------------------------------
\section{A hitchhiker's guide for the ``orthogonal'' and the ``hybrid choice''} \label{appendix hitchhikers}
% --------------------------------------------------------------------------------------------------------------------------------------------------------------------------------------------------------------------------------------------------------------------
This appendix is devoted to discuss some implementation details of the local VEM stiffness matrix by employing specific polynomial bases introduced in Section \ref{subsection choice of dofs}
in the face \eqref{internal moments faces polyhedron} and bulk \eqref{internal moments bulk} moments.

Since the ``standard choice'' has already been investigated in \cite{hitchhikersguideVEM}, we split this appendix into two parts:
in Appendix \ref{subappendix pure fetish}, we discuss the details when employing the ``orthogonal choice'' whereas, in Appendix \ref{subappendix hybrid fetish}, we discuss the details when employing the ``hybrid choice''.

Henceforth we fix some notations; given $\E$ polyhedron in $\taun$, we set $\Ndof$ the dimension of $\Vp(\E)$ and, given $\F$ a face of $\E$, we set $\Ndof^\F$ the dimension of $\Vp(\F)$;
moreover, by $\NS$, $\NF$ and $\NB$ we denote the number of skeletal, face and bulk dofs in local space $\Vp(\E)$ defined in \eqref{local space bulk}, respectively.
In Section \ref{subsection choice of dofs}, we defined polynomial basis up to order $\p-2$. In this appendix, we employ polynomials up to order $\p$.
In particular, we write $\{\malpha^\F\}_{\alpha=1}^{\np^\F}$ and $\{\malpha\}_{\alpha=1}^{\np}$ to denote the monomials of degree $\p$ on face $\F$ and in polyhedron $\E$,
while we write $\{\mbaralpha^\F\}_{\alpha=1}^{\np^\F}$ and $\{\mbaralpha\}_{\alpha=1}^{\np}$ to denote the $L^2$ orthonormal polynomials on face $\F$ and in polyhedron $\E$ 
obtained from $\{\malpha^\F\}_{\alpha=1}^{\np^\F}$ and $\{\malpha\}_{\alpha=1}^{\np}$ via a stable Gram-Schmidt process, respectively.

\medskip

When dealing with calculations between vectors-matrices, we employ the following notation:
\[
\mathbf A(i:j, l:k) \quad \quad \quad \forall \, \mathbf A\in \mathbb R^{n,m},
\]
is the submatrix of $\mathbf A$ from row $i$ to $j$ and from column $\ell$ to column $k$.
If no indications concerning rows-columns are given, then it means that we are considering the full matrix.

We extensively make usage of the two natural bijections \eqref{bijection 2D} and \eqref{bijection 3D}.

The construction of the local stiffness matrix in the ``standard case'' is based on the following matrices, whose implementation has been already discussed in \cite{hitchhikersguideVEM}, both in polyhedron $\E$:
\begin{equation} \label{old VEM local matrices}
\begin{aligned}
& \G_{\alpha,\beta} = 	\begin{cases}
					\begin{cases}
						\sum_{i=1}^{\NV} (\mbeta(\nu_i)) 	& \text{if } \p=1\\
						\int_\E \mbeta 				& \text{otherwise}\\
					\end{cases} & \text{if } \alpha=1\\
					(\nabla \malpha, \nabla \mbeta)_{0,\E} & \text{otherwise}\\
				\end{cases},\\
& \Gtilde _{\alpha, \beta} = (\nabla \malpha, \nabla \mbeta)_{0,\E}, \quad \B_{\alpha,i} = \begin{cases}
															P_0 \varphi_i & \text{if } \alpha=1\\
															(\nabla \malpha, \nabla \varphi_i) _{0,\E} & \text{otherwise}\\
														\end{cases},\\
& \D_{i,\alpha} = \dof_i (\malpha), \quad \H_{\alpha, \beta} =(\malpha, \mbeta)_{0,\E},\quad 	\C_{\alpha,i} = (\malpha, \varphi_i)_{0,\E},\\
\end{aligned}
\end{equation}
for all $\alpha$, $\beta= 1,\dots,\np$ and for all $i=1,\dots, \Ndof$, where we recall that $\NV$ and $\{\nu_i\}_{i=1}^{\NV}$ are the number and the set of vertices of $\E$,
both on face $\F$:
\begin{equation} \label{old VEM local matrices faces}
\begin{aligned}
&\G^\F_{\alpha,\beta} = 	\begin{cases}
					\begin{cases}
						\sum_{i=1}^{\NV^\F} (\mbeta^\F(\nu_i^\F)) 	& \text{if } \p=1\\
						\int_\E \mbeta^\F 				& \text{otherwise}\\
					\end{cases} & \text{if } \alpha=1\\
					(\nabla \malpha^\F, \nabla \mbeta^\F)_{0,\F} & \text{otherwise}\\
				\end{cases},\\
&\Gtilde^\F _{\alpha, \beta} = (\nabla \malpha^\F, \nabla \mbeta^\F)_{0,\E}, \quad \B_{\alpha,i} = \begin{cases}
															P_0 \varphi_i & \text{if } \alpha=1\\
															(\nabla \malpha^\F, \nabla \varphi_i) _{0,\F} & \text{otherwise}\\
														\end{cases},\\
&\D_{i,\alpha} = \dof_i (\malpha^\F), \quad \H_{\alpha, \beta} =(\malpha^\F, \mbeta^\F)_{0,\E},\quad 	\C_{\alpha,i} = (\malpha^\F, \varphi_i)_{0,\F},\\
\end{aligned}
\end{equation}
for all $\alpha$, $\beta= 1,\dots,\np^\F$ and for all $i=1,\dots, \Ndof^\F$, where we recall that $\NV^\F$ and $\{\nu_i^\F\}_{i=1}^{\NV^\F}$ are the number and the set of vertices of $\F$.

\medskip

In the two forthcoming appendices, we explain how to construct the counterparts of the matrices defined in \eqref{old VEM local matrices} and \eqref{old VEM local matrices faces}
with the ``orthogonal'' and the ``hybrid choices''.
Since with these two choices we employ $L^2$ orthonormal polynomial bases, we also need the (lower triangular) matrices containing the orthonormalizing coefficients
with respect to the monomial bases of $\mathbb P_{\p} (\E)$ and $\mathbb P_{\p} (\F)$, respectively.
Such matrices are denoted by $\GS$, matrix belonging to $\mathbb R^{\np \times \np}$ (on polyhedron $\E$), and $\GS^\F$, matrix belonging to $\mathbb R^{\np^\F \times \np^\F}$ (on face $\F$).

In the remainder of this appendix, we denote the local VEM matrices, the local degrees of freedom and the local canonical basis functions with a bar at the top of each of them.

%%%%%%%%%%%%%%%%%%%%%%%%%%%%%%%%%%%%%%%%%%%%%%%%%
\subsection{A hitchhiker's guide for the ``\orthogonal choice''} \label{subappendix pure fetish}
The aim of the present appendix, is to give some details for what concerns the computation of the counterpart of the matrices in \eqref{old VEM local matrices} employing the face/bulk polynomial bases of the so-called ``\orthogonal choice''
discussed in Section \ref{subsection choice of dofs}. In particular, we fix bases made of $L^2$ orthonormal polynomials both on faces and in the bulk.

The assembling of the local stiffness matrix boils down to the construction in \cite{hitchhikersguideVEM} and depends on the choice of the local stabilization, see Section \ref{subsection choice stabilizations}.

\subsubsection{Matrices $\Gbar$ and $\Gtildebar$} \label{subsubsection matrices Gbar and Gtildebar}
We start with matrix $\Gtildebar$ which is defined as:
\[
\Gtildebar = (\nabla \mbaralpha, \nabla \mbarbeta)_{0,\E} \quad \forall \, \alpha,\, \beta= 1,\dots, \np.
\]
One simply has:
\[
\Gtildebar = \GS \cdot \Gtilde \cdot \GS^T.
\]
Next, we consider matrix $\Gbar$ defined as:
\[
\Gtildebar_{\boldalpha, \boldbeta} =
\begin{cases}
P_0(\mbarbeta) & \text{if } \alpha = 1\\
(\nabla \mbaralpha, \nabla \mbarbeta)_{0,\E} & \text{otherwise}
\end{cases}\quad \forall \, \alpha,\, \beta= 1,\dots, \np,
\]
where, recalling that $\{\nu_i\}_{i=1}^{\NVE}$ is the set of vertices of $\E$:
\[
P_0(\cdot) = 
\begin{cases}
\frac{1}{\NVE} \sum_{i=1}^{\NVE} \cdot (\nu_i) & \text{if } \p=1\\
\frac{1}{\vert \E \vert} \int_\E \cdot & \text{otherwise}\\
\end{cases}.
\]
Clearly, we only have to treat the case $\alpha=1$. We distinguish two cases.

If $\p=1$, then we have:
\[
P_0(\mbarbeta) = \frac{1}{\NVE} \mbarbeta(\nu_i) \quad \forall \, \beta= 1,\dots, \np
\]
and therefore
\[
\begin{aligned}
P_0(\mbar_\beta) 	& = \frac{1}{\NVE} \sum _{i=1}^{\NVE} \left(\sum_{\gamma \le \beta} \GS_{\beta, \gamma} \m_{\gamma}(\nu_i) \right) = \sum_{\gamma \le \beta} \GS_{\beta,\gamma} \left( \frac{1}{\NVE} \sum_{i=1}^{\NVE} \m_\gamma (\nu_i)  \right)\\
			& = \sum_{\gamma \le \beta} \GS_{\beta, \gamma} \G_{1,\gamma} = \GS(\beta, 1:\beta) \cdot \G(1,1:\beta)^T \quad \forall \, \beta= 1,\dots, \np.
\end{aligned}
\]
If, on the other hand, $\p \ge 2$:
\[
\begin{aligned}
P_0(\mbar_\beta) 	& = \frac{1}{\vert \E \vert} \int_\E \mbar _\beta = \GS_{1,1}^{-1} \frac{1}{\vert \E \vert} \int _\E \mbar _\beta \GS_{1,1}\\
			& = \GS_{1,1}^{-1} \frac{1}{\vert \E \vert} \int _\E \mbar_1 \mbar_\beta  = 	\begin{cases}
															\frac{1}{\GS_{1,1} \vert \E \vert} & \text{if } \beta=1\\
															0                         & \text{else}\\
															\end{cases}\quad \forall \, \beta= 1,\dots, \np,
\end{aligned}
\]
since basis $\{\mbar_\alpha\}_{\alpha=1}^{\np}$ is $L^2(\E)$ orthonormal by construction.

%-------------------------------------------------------------------------------------------------------------------------------------------
\subsubsection{Matrix $\Dbar$} \label{subsubsection Dbar}
We define matrix $\Dbar$ as:
\[
\Dbar_{i,\alpha} = \dofbar_i(\mbar_\alpha) \quad \forall \, i=1,\dots, \Ndof,\, \forall \, \alpha= 1,\dots, \np.
\]
Let us consider firstly the boundary dofs. For all $i=1,\dots, \NS$ and for all $\alpha= 1,\dots, \np$:
\[
\Dbar_{i,\alpha} = \dofbar_i(\mbar_\alpha) = \mbar_ \alpha(\xi_i) = \sum_{\beta \le \alpha} \GS_{\alpha, \beta}\, \m_\beta(\xi_i) = \sum_{\beta\le \alpha} \GS_{\alpha, \beta} \D_{i,\beta} = \D(i, 1:\alpha)^T \cdot \GS (\alpha, 1: \alpha)^T,
\]
where $\xi_i$ is a proper node on the boundary. In short:
\[
\Dbar (1:\NS, 1: \np) = \D (1:\NS, 1:\np) \cdot \GS^T.
\]

Next, we deal with the face dofs. Assume that the $i$-th $\dofbar$ is associated with polynomial $\mbargamma^\F$.
Then, for all $i=\NS+1, \dots, \NS+\NF$ and for all $\alpha= 1,\dots, \np$:
\[
\dofbar_{i} (\mbaralpha) = \frac{1}{\vert \F \vert} \int_\F \mbargammaF \mbaralpha = \sum_{\beta \le \alpha} \GS_{\alpha,\beta} \frac{1}{\vert \F \vert} \int_\F \mbargammaF \mbeta 
						= \sum_{\beta \le \alpha} \GS_{\alpha, \beta} \left\{ \sum_{\delta \le \gamma} \GS_{\gamma,\delta}^\F \frac{1}{\vert \F \vert} \int_\F \mdeltaF  \mbeta \right\}.
\]
In short:
\[
\Dbar (\NS+1:\NS+\NF+1, 1:\np) = \GS^\F(1:\npmd^\F, 1: \npmd^\F) \cdot \D (\NS+1:\NS+\NF+1, 1 : \npmd) \cdot \GS^T.
\]

Finally, we treat the bulk dofs. Assume that the $i$-th dof is associated with polynomial $\mbargamma$. Then, for all $i=\NS+\NF+1,\dots, \NS+ \NF + \NB$ and for all $\alpha= 1,\dots, \np$:
\[
\dofbar_i(\mbaralpha) = \frac{1}{\vert \E \vert} \int_\E \mbaralpha \mbargamma = \frac{1}{\vert \E \vert} \delta_{\alpha,i}.
\]

%-------------------------------------------------------------------------------------------------------------------------------------------
\subsubsection{Matrix $\Hbar$} \label{subsubsection matrix Hbar}
We define matrix $\Hbar$ as:
\[
\Hbar_{\alpha, \beta} = (\mbaralpha, \mbarbeta)_{0,\E}\quad \forall \, \alpha,\, \beta= 1,\dots, \np.
\]
One directly has:
\[
\Hbar = \GS \cdot \H \cdot \GS ^T = \Id.
\]
As a byproduct, we observe that by verifying this last equality one can check whether $\GS$ has been actually properly computed.

%-------------------------------------------------------------------------------------------------------------------------------------------
\subsubsection{Matrix $\Bbar$} \label{subsubsection matrix Bbar}
We define matrix $\Bbar$ as follows:
\[
\Bbar_{\alpha,i} =
\begin{cases}
P_0(\phibar _i) & \text{if }\alpha=1\\
(\nabla \mbaralpha, \nabla \phibar_i) _{0,\E} & \text{otherwise}\\
\end{cases} \quad \forall \, \alpha= 1,\dots, \np,\, \forall \, i = 1,\dots, \Ndof.
\]
We firstly deal with the first line and we consider the two cases $\p=1$ and $\p\ge 2$.
\begin{itemize}
\item[($\p=1$)] $P_0(\phibar_i) = \frac{1}{\NVE} \sum_{i=1}^{\NVE} \phibar_i(\nu_i)$. Thus, $\Bbar_{1,i} = \B_{1,i}$ for all $i = 1,\dots, \Ndof$, since the elements of the new basis coincide with the old ones on the skeleton of the mesh.
\item[($\p\ge 2$)] In this case:
\[
\begin{aligned}
P_0(\phibar_i) = 	&\frac{1}{\vert \E \vert} \int_\E \phibar_i = \GS^{-1}_{1,1} \frac{1}{\vert \E \vert} \int_\E \phibar_i \GS_{1,1} = \GS^{-1}_{1,1} \frac{1}{\vert \E \vert} \int_\E \phibar_i \mbar_1 \\
			& = \begin{cases}
				\GS^{-1}_{1,1} & \text{if } \phibar_i \text{ is the first bulk basis element}\\
				0 & \text{else}\\
				\end{cases} \quad \forall \, i = 1, \dots, \Ndof,
\end{aligned}
\]
since $\mbar_1 = \GS_{1,1} \m_1 = \GS_{1,1}$. Thus, we can copy the old line and multiply it for $\GS_{1,1}^{-1}$, i.e. $\Bbar_{1,i} = \GS_{1,1}^{-1} \B_{1,i}$.
\end{itemize}
\medskip \medskip 
Next, we treat the remaining lines. We must compute $(\nabla \mbaralpha, \nabla \phibar_i)_{0,\E}$. We consider three different situations.
\begin{itemize}
\item We assume $\phibar_i$ edge basis element. Then:
\[
\begin{aligned}
& (\nabla \mbaralpha, \nabla \phibar_i)_{0,\E} = (-\Delta \mbaralpha, \phibar_i)_{0,\E} + (\partial _\n \mbaralpha, \phibar_i )_{0,\partial \E} 	\\
& = \sum_{\F \in \partial \E} (\partial _\n \mbaralpha, \phibar_i )_{0,F} \quad \forall \, \alpha= 1,\dots, \np,\, \forall \, i=1,\dots,\NS.
\end{aligned}
\]
Therefore, it suffices to compute integrals over faces. For the purpose, we decompose $\partial _\n \mbaralpha$ on each face $\F$ as a linear combination of elements in the $L^2(\F)$ orthonormal basis $\{\mbarbeta^\F\}_{\vert \boldbeta\vert=0}^{\p-1}$:
\begin{equation} \label{expansion normal derivative}
\partial _\n \mbaralpha | _{\F} = \sum_{\vert \boldbeta \vert=0}^{\p-1} \lambdabar_{\alpha, \beta}^\F \mbarbeta^\F \quad \quad \forall \, \F \in \partial \E, \quad \forall \, \alpha= 1,\dots, \np.
\end{equation}
In order to be able to compute the coefficients $\lambdabar_{\alpha,\beta}^\F$, we test \eqref{expansion normal derivative} with $\mbargamma^\F$ with $\gamma= 1,\dots, \npmu$ and get by orthonormality:
\begin{equation} \label{lambda face}
\begin{aligned}
\lambdabar_{\alpha,\n}^\F = (\partial _\n \mbaralpha, \mbargamma^\F)_{0,\F}
&= \sum_{\beta \le \alpha} \GS_{\alpha, \beta} \left\{ \sum_{ \delta \le \gamma} \GS_{\gamma, \delta} ^\F (\partial _ \n \mbeta, \mdelta^\F)_{0,\F}     \right\}\\
& \quad \forall \, \alpha= 1,\dots, \np, \, \forall \, \gamma= 1,\dots, \npmu,\\
\end{aligned}
\end{equation}
which is easily computable:
\[
\lambdabar_{\alpha,\gamma}^\F = \GS (\alpha,1: \np) \cdot \boldsymbol\Lambda^\F \cdot \GSF(\gamma, 1:\npmu^\F)^T \quad \forall \, \alpha= 1,\dots, \np,\, \forall \, \gamma= 1,\dots, \npmu,
\]
where:
\[
\boldsymbol \Lambda ^\F _{\beta,\delta}= (\partial_\n \mbeta, \mdelta^\F)_{0,\F} \quad \quad \forall \, \beta = 1,\dots, \np,\, \forall \, \delta= 1,\dots, \npmu
\]
is computable by a simple integration of $L^2$ products on monomials.

In short, if we call $\boldsymbol \Lambdabar^\F$ the matrix of the coefficients of expansion \eqref{expansion normal derivative} on face $\F$, then we have:
\[
\boldsymbol \Lambdabar^\F = \GS \cdot \boldsymbol \Lambda^\F \cdot \GSF(1:\npmu^\F, 1:\npmu^\F)^T.
\]

As a consequence, on each face $\F$,  we get by using the definition of the enhancing constraints:
\begin{equation} \label{a bar integral on boundary}
\begin{aligned}
(\partial_\n \mbaralpha, \phibar_i)_{0,\F} 	& = \sum_{\vert \boldbeta \vert=0}^{\p-1} \lambdabar_{\alpha,\beta}^\F (\mbarbeta^\F, \phibar_i)_{0,\F} = \sum_{\vert \boldbeta \vert=\p-1} \lambdabar_{\alpha,\beta}^\F (\mbarbeta^\F, \phibar_i)_{0,\F} \\
							& = \sum_{\vert \boldbeta \vert=\p-1} \lambdabar_{\alpha,\beta}^\F (\mbarbeta^\F, \PinablabarF \phibar_i)_{0,\F} \quad \forall \, \alpha= 1,\dots, \np,\, \forall i=\NS+1,\dots, \NS+\NF,
\end{aligned}
\end{equation}
where $\PinablabarF$ denotes the $H^1$ projection on the polynomial space spanned by the $L^2(\F)$ orthonormal basis $\{\mbaralpha^\F\}_{\alpha=1}^{\np^\F}$, which can be computed on face $\F$ following \cite{fetishVEM}.
The quantity in \eqref{a bar integral on boundary} can be computed; in fact:
\begin{equation} \label{face integrals Pinabla}
\sum_{\vert \boldbeta \vert=\p-1} \lambdabar_{\alpha,\beta}^\F (\mbarbeta^\F, \PinablabarF \phibar_i)_{0,\F} =
\sum_{\vert \boldbeta \vert=\p-1} \lambdabar_{\alpha,\beta}^\F \left(\mbarbeta^\F, \sum_{\vert \boldgamma \vert = 0}^{\p} \chibar_{\gamma, i}^\F \mbargamma ^\F \right)_{0,\F} =
\sum_{\vert \boldbeta \vert = \p-1} \lambdabar_{\alpha,\beta}^\F \chibar_{\beta, i}^\F.
\end{equation}
In order to conclude, one sums \eqref{face integrals Pinabla} on all the faces.

\item Let now $\phibar_i$ be a face basis element. Then:
\begin{equation}\label{grad m grad phi}
\begin{aligned}
(\nabla \mbaralpha, \nabla \phibar_i)_{0,\E} 	& = (\partial_\n \mbaralpha, \phibar_i)_{0,\partial \E} \\
							& = \sum_{\F \in \partial \E}  (\partial_\n \mbaralpha, \phibar_i)_{0,\F} \quad \forall \, \alpha= 1,\dots, \np,\, \forall \, i=\NS + 1,\dots, \NS+\NF.
\end{aligned}
\end{equation}
Plugging expansion \eqref{expansion normal derivative} into \eqref{grad m grad phi} and denoting by $\widetilde \F$ the face associated with basis element $\overline \varphi _i$:
\[
\begin{aligned}
& (\nabla \mbaralpha, \phibar _i)_{0,\E} 	= \sum_{\F \in \partial \E}  \sum_{\vert \boldbeta\vert=0}^{\p-1} \lambda^\F_{\alpha,\beta} (\mbarbeta^\F, \phibar _i)_{0,\F} \\
& = \left( \sum_{\F \in \partial \E} \sum_{\vert \boldbeta\vert=\p-1} \lambda^\F_{\alpha,\beta} (\mbarbeta^\F, \PinablabarF \phibar _i)_{0,\F}   \right) + \vert \widetilde \F \vert \lambdabar_{\alpha,i}^{\widetilde \F}
								\quad \forall\, \alpha= 1,\dots, \np,\, \forall\, i=\NS+1,\dots,  \NS+\NF,
\end{aligned}
\]
where $\lambdabar^\F_{\alpha,i}$ is computed as in \eqref{lambda face} for all faces $\F \in \partial \E$
and where we are assuming that $\phibar_i$ is the face element associated with $\mbarbeta^{\widetilde \F}$.
The integrals over the faces involving $\Pinablap$ can be computed as in \eqref{face integrals Pinabla}.
 	\item Finally, we assume that $\phibar_i$ is a bulk basis element. Then:
\begin{equation} \label{insert laplacian}
(\nabla \mbaralpha, \nabla \phibar_i)_{0,\E} = (-\Delta \mbaralpha, \phibar_i)_{0,\E} \quad \forall \, \alpha= 1,\dots, \np,\, \forall\, i=\NS+\NF+1,\dots, \NS+\NF+\NB.
\end{equation}
As a consequence, we have to expand $-\Delta \mbaralpha$ in terms of the $L^2(\E)$ orthonormal basis $\{\mbaralpha\}_{\alpha=1}^{\npmd}$ as:
\begin{equation} \label{expansion L}
-\Delta \mbaralpha = \sum_{\vert \boldbeta \vert=0}^{\p-2} \mubar_{\alpha, \beta} \mbarbeta \quad \forall \, \alpha= 1,\dots, \np.
\end{equation}
We test \eqref{expansion L} with $\mbargamma$, $\vert \boldgamma \vert = 0,\dots,\p-1$ and get by orthonormality:
\begin{equation} \label{expansion laplacian}
\mubar_{\alpha, \beta} = (-\Delta \mbaralpha, \mbargamma)_{0,\E} = \Lbar_{\alpha,\gamma} \quad \forall \, \alpha= 1,\dots, \np,\, \forall \, \beta= 1,\dots, \npmd,
\end{equation}
which is actually computable. In fact:
\[
(-\Delta \mbaralpha, \mbargamma)_{0,\E} = \GS(\alpha, 1:\np) \cdot \L \cdot \GS(\gamma,1:\npmd)^T \quad \forall \, \alpha= 1,\dots, \np,\, \forall \, \gamma = 1,\dots, \npmd,
\]
where matrix $\L$ can be computed as:
\[
\begin{aligned}
&\L_{\alpha,\beta} = (-\Delta \malpha, \mgamma)_{0,\E}\\
& = -\frac{1}{\hE^2} \left( \alpha_1(\alpha_1-1) (m_{\boldalpha_1}, \mgamma) _{0,\E} +\alpha_2(\alpha_2-1) (m_{\boldalpha_2}, \mgamma) _{0,\E}
+ \alpha_3(\alpha_3-1) (m_{\boldalpha_3}, \mgamma) _{0,\E} \right)\\
& = -\frac{1}{\hE^2} \left( \alpha_1(\alpha_1-1) \H_{\boldalpha_1, \beta} + \alpha_2(\alpha_2-1) \H_{\boldalpha_2, \beta} + \alpha_3(\alpha_3-1) \H_{\boldalpha_3, \beta}   \right) \\
& \forall \, \alpha= 1,\dots, \np,\, \forall \, \beta = 1,\dots, \npmd,
\end{aligned}
\]
having set:
\[
\boldalpha_1 = (\alpha_1-2,\alpha_2,\alpha_3),\quad \boldalpha_2 = (\alpha_1,\alpha_2-2,\alpha_3),\quad \boldalpha_3 = (\alpha_1,\alpha_2,\alpha_3-2).
\]

Inserting \eqref{expansion L} and \eqref{expansion laplacian} in \eqref{insert laplacian}, we obtain:
\[
\begin{aligned}
(\nabla \mbaralpha, \nabla \phibar_i)_{0,\E} & = \sum_{\vert \beta \vert = 0}^{\p-2} \mu_{\alpha, \beta} (\mbarbeta, \phibar_i)_{0,\E}\\
& = \vert \E \vert \mu_{\alpha,i} \quad \forall \, \alpha= 1,\dots, \np,\,
\forall\, i = \NS+\NF+1,\dots, \NS+\NF+\NB.
\end{aligned}
\]
\end{itemize}

%-------------------------------------------------------------------------------------------------------------------------------------------
\subsubsection{Matrix $\Cbar$} \label{subsubsection matrix Cbar}
We define matrix $\Cbar$ as:
\[
\Cbar_{\alpha,i} = (\mbaralpha, \phibar_i)_{0,\E} \quad \forall \, \alpha=1,\dots,\np ,\, \forall \, i = 1,\dots,\Ndof.
\]
It is also clear that, for all $\alpha=1,\dots,\npmd$ and for all $ i = 1,\dots,\Ndof$:
\[
(\mbaralpha, \phibar_i)_{0,\E} =
\begin{cases}
0 & \text{if } \phibar_i \text{ is an edge or face basis element}\\
\vert \E \vert \frac{1}{\vert \E \vert} \int_\E \mbaralpha \phibar_i = \delta_{i,\alpha} \vert \E \vert & \text{otherwise}\\
\end{cases}.
\]
For what concerns the other lines, one employs the so-called enhancing technique \cite{hitchhikersguideVEM}. More precisely, one sets:
\[
\Cbar (\npmd +1:\np, 1:\Ndof) = \left[ \Gbar^{-1} \cdot \Bbar  \right] (\npmd+1: \np, 1:\Ndof).
\]

%%%%%%%%%%%%%%%%%%%%%%%%%%%%%%%%%%%%%%%%%%%%%%%%%
\subsection{A hitchhiker's guide for the ``\hybrid choice''} \label{subappendix hybrid fetish}
The aim of the present appendix, is to give some details for what concerns the computation of the counterpart of the matrices in \eqref{old VEM local matrices} employing the face/bulk polynomial bases of the so-called ``\hybrid choice''
discussed in Section \ref{subsection choice of dofs}. In particular, we fix bases made of $L^2$ orthonormal polynomials in the bulk and standard monomials on faces.

The assembling of the local stiffness matrix boils down to the construction in \cite{hitchhikersguideVEM} and depends on the choice of the local stabilization, see Section \ref{subsection choice stabilizations}.

Henceforth, we denote with a bar at the top and \hyb{}  as subscript the local VEM matrices, the local degrees of freedom and the local canonical basis functions. It is easy to check that:
\[
\Gtildebar_\hyb = \Gtildebar, \quad \quad \Gbar_\hyb = \Gbar, \quad \quad \Hbar_\hyb = \Hbar.
\]
For what concerns matrix $\Dbar_\hyb$ we observe that is coincides with matrix $\Dbar$, with the exception of the entries related to face dofs. In this case, one has:
\[
\Dbar_\hyb (\NS+1: \NS + \NF, 1:\np) = \D (\NS+1: \NS+ \NF, 1:\npmd) \cdot \GS^T.
\]
The treatment of matrix $\Bbar_\hyb$ is rather different. It coincides with $\Bbar$ when considering the first line and when considering all the columns associated with bulk dofs.

Let us fix our attention to the columns associated with skeleton dofs. In this case, we have:
\begin{equation} \label{first equation hybrid Bbar}
(\nabla \mbaralpha, \nabla \phibar_{\hyb,i})_{0,\E} = \sum_{\F \in \partial \E} (\partial_n \mbaralpha, \phibar_{\hyb,i}) \quad \forall \, \alpha=1,\dots,\np,\, \forall \, i=1,\dots, \NS.
\end{equation}
We expand $\partial_n \mbaralpha$ on each face $\F$ into a combination of elements $\{\mbeta^\F\} _{\beta=1}^{\npmu^\F}$, which we recall is the standard (scaled and centered) monomial basis on face $\F$:
\begin{equation} \label{hybrid expansion normal derivative}
\partial_\n \mbaralpha|_\F = \sum_{\beta =1}^{\npmu^\F} \lambdabar^\F_{\hyb, \alpha,\beta} \mbeta^\F \quad \forall \, \F \in \partial \E,\, \forall \, \alpha=1,\dots, \np.
\end{equation}
By testing \eqref{hybrid expansion normal derivative} with $\mgammaF$, for $\gamma=1,\dots,\npmu^\F$, we write for all faces $\F$:
\[
\sum_{\beta = 1}^{\npmu^\F} \lambdabar^\F_{\hyb, \alpha,\beta} (\mbetaF, \mgammaF)_{0,\F} = (\partial_n \mbaralpha, \mgammaF)_{0,\F}.
\]
Therefore, if we want to have an explicit representation of the coefficients $\lambdabar^\F_{\hyb, \alpha, \beta} = \Lambdabar^\F_\hyb$
for $\alpha=1,\dots, \np$, $\beta=1,\dots, \npmu^\F$, we have to  compute:
\[
\Lambdabar^\F_\hyb = \Mbar ^\F_\hyb \cdot \H^\F (1: \npmu^\F, 1:\npmu^\F)^{-1},
\]
where $\H^\F$ denotes the standard 2D VEM matrix on face $\F$ and where:
\[
(\Mbar_\hyb^\F)_{\alpha,\gamma} = (\partial _\n \mbar, \mgammaF)_{0,\F} \quad \forall\, \alpha=1,\dots,\np,\, \forall \, \gamma =1,\dots, \npmu^\F,
\]
which is easily computable.

Having now the coefficients $\lambdabar^\F_{\hyb, \alpha,\beta}$, we obtain from \eqref{first equation hybrid Bbar} and using that $\phibar_{\hyb, i} |_{\partial \E} = \varphi_i|_{\partial \E}$ for all $i=1,\dots, \NS$:
\[
\begin{aligned}
(\partial  \mbaralpha, \nabla \phibar _{\hyb,i})_{0,\E}	& = \sum_{\F \in \partial \E} (\partial _\n \mbaralpha, \phibar_{\hyb, i})_{0,\F} = \sum_{\F \in \partial \E} (\partial _\n \mbaralpha, \varphi_{i})_{0,\F}\\
									& = \sum_{\F \in \partial \E} \left\{ \sum_{\beta=1}^{\npmu^\F} \lambdabar^\F_{\hyb, \alpha, \beta} (\mbeta^\F, \varphi_i)_{0,\E}   \right\} = \sum_{\F \in \partial \E} \left\{ \sum_{\beta=\npmd^\F+1}^{\npmu^\F} \lambdabar^\F_{\hyb, \alpha, \beta} (\mbeta^\F, \varphi_i)_{0,\E}   \right\}\\
									& = \sum_{\F \in \partial \E} \left\{ \sum_{\beta=\npmd^\F+1}^{\npmu^\F} \lambdabar^\F_{\hyb, \alpha, \beta} (\mbeta^\F, \PinablaF \varphi_i)_{0,\E}   \right\}
												\quad \forall \, \alpha=1,\dots, \np,\, \forall \, i=1,\dots,\NS.
\end{aligned}
\]
This is equivalent to say:
\[
(\Bbar_\hyb)_{\alpha,i} = \Lambdabar^\F_\hyb (1:\np, \npmd^\F+1: \npmu^\F) \cdot \H^\F (\npmd^\F+1: \npmu^\F) \cdot \mathbf{\Pi^{\nabla, *,\F}} (1:\np, 1:\NS),
\]
where $\mathbf{\Pi^{\nabla, *,\F}} = \Gbar_\hyb^{-1} \cdot \Bbar_{\hyb}$.

The case $\phibar_{\hyb,i}$ face basis element is dealt with in an utterly analogous way by noting that:
\[
(\nabla \mbaralpha, \nabla \phibar _{\hyb, i})_{0,\E} = \sum_{\beta = \npmd^{\widetilde \F} +1 }^{\npmu^{\widetilde \F}} \lambdabar^{\widetilde \F}_{\hyb, \alpha,\beta} (\mbeta^{\widetilde \F}, \Pi^{\nabla, \widetilde \F} \varphi_i)_{0,\widetilde \F}
+ \vert \widetilde \F \vert,
\]
where $\widetilde F$ is the face associated with $\phibar^\F_{\hyb,i}$.

\end{appendix}

\section*{Acknowledgement}
The first  author was partially supported by the European Research Council through the
H2020 Consolidator Grant (grant no. 681162) CAVE -- Challenges and Advancements in Virtual Elements.

%%%%%%%%%%%%%%%%%%%%%%%%%%%%%%%%%%%%%%%%%%%%%%%%%%%%%%%%%%%%%%%%%%%%%%%%%%%
{\footnotesize
\bibliography{bibliogr}
}
\bibliographystyle{plain}
%%%%%%%%%%%%%%%%%%%%%%%%%%%%%%%%%%%%%%%%%%%%%%%%%%%%%%%%%%%%%%%%%%%%%%%%%%%
\end{document}